\documentclass[reqno]{amsart}
\usepackage{amssymb}
\usepackage{url}
\usepackage{mathrsfs}
\usepackage{graphicx}
\def\qmod#1#2{{\hbox{}^{\displaystyle{#1}}}\!\big/\!\hbox{}_{
\displaystyle{#2}}}

\def\resto#1#2{{
#1\hskip 0.4ex\vline_{\hskip 0.2ex\raisebox{-0.2ex}
{{${\scriptstyle #2}$}}}}}
\def\C{{\mathbb C}}

\def\H{{\mathbb H}}

\def\N{{\mathbb N}}
\def\O{{\mathbb O}}
\def\P{{\mathbb P}}
\def\Q{{\mathbb Q}}
\def\R{{\mathbb R}}

\def\V{{\mathbb V}}
\def\Z{{\mathbb Z}}

\def\union{\mathop{\bigcup}}
\def\qed {\hfill\vrule height6pt width6pt depth0pt \bigskip}
\def\map{\longrightarrow}
\def\textmap#1{\mathop{\vbox{\ialign{
                                 ##\crcr
     ${\scriptstyle\hfil\;\;#1\;\;\hfil}$\crcr
     \noalign{\kern 1pt\nointerlineskip}
     \rightarrowfill\crcr}}\;}}
\newcommand{\cal}{\mathcal}
\def\textlmap#1{\mathop{\vbox{\ialign{
                                 ##\crcr
     ${\scriptstyle\hfil\;\;#1\;\;\hfil}$\crcr
     \noalign{\kern-1pt\nointerlineskip}
     \leftarrowfill\crcr}}\;}}

\newfam\meuffam

\font\tenmeuf=eufm10
\font\sevenmeuf=eufm7
\font\fivemeuf=eufm5

\textfont\meuffam=\tenmeuf
\scriptfont\meuffam=\sevenmeuf
\scriptscriptfont\meuffam=\fivemeuf

\def\germ{\fam\meuffam\tenmeuf}
\def\ag{{\germ a}}

\def\g{{\germ g}}

\def\kg{{\germ k}}

\def\zg{{\germ z}}
\def\Ag{{\germ A}}

\def\Cg{{\germ C}}

\def\Fg{{\germ F}}

\def\Hg{{\germ H}}

\def\Mg{{\germ M}}

\def\Rg{{\germ R}}

\def\Sg{{\germ S}}

\newtheorem{sz}{Satz}[section]
\newtheorem{thry}[sz]{Theorem}
\newtheorem{pr}[sz]{Proposition}
\newtheorem{re}[sz]{Remark}
\newtheorem{co}[sz]{Corollary}

\newtheorem{lm}[sz]{Lemma}

\begin{document}
\def\tr{\mathrm{Tr}}
\def\End{\mathrm{End}}
\def\Aut{{\rm Aut}}
\def\Spin{\mathrm{Spin}}
\def\U{\mathrm{U}}
\def\O{\mathrm{O}}
\def\SU{\mathrm{SU}}
\def\SO{\mathrm{SO}}
\def\PU{\mathrm{PU}}
\def\GL{\mathrm{GL}}
\def\SL{\mathrm{SL}}
\def\gl{\mathrm{gl}}
\def\spin{\mathrm{spin}}
\def\su{\mathrm{su}}
\def\so{\mathrm{so}}
\def\sl{\mathrm{sl}}
\def\ub{\underbar}
\def\proj{\mathrm{pr}}
\def\pu{\mathrm{pu}}
\def\Pic{\mathrm{Pic}}
\def\Iso{\mathrm{Iso}}
\def\NS{\mathrm{NS}}
\def\deg{\mathrm{deg}}
\def\Hom{\mathrm{Hom}}
\def\Aut{\mathrm{Aut}}
\def\Tors{\mathrm{Tors}}
\def\Sph{\mathrm{S}}
\def\Herm{\mathrm{Herm}}
\def\Vol{\mathrm{Vol}}
\def\vol{\mathrm{vol}}

\def\pf{{\bf Proof: }}
\def\id{{\rm id}}
\def\i{{\germ i}}
\def\im{{\rm im}}
\def\rk{{\rm rk}}
\def\ad{{\rm ad}}
\def\coker{{\rm coker}}
%%%%%%%%%%%%%%%%%%%%%%%%%%
\def\dbar{\bar{\partial}}
\def\Lo{{\Lambda_g}}
\def\niq{=\kern-.18cm /\kern.08cm}
%%%%%%%%%%%%%%%%%%%
\def\Ad{{\rm Ad}}
\def\RSU{\R SU}
\def\ad{{\rm ad}}
\def\dva{\bar\partial_A}
\def\da{\partial_A}
\def\p{\partial\bar\partial}
\def\sp{\Sigma^{+}}
\def\sm{\Sigma^{-}}
\def\spm{\Sigma^{\pm}}
\def\smp{\Sigma^{\mp}}
\def\st{{\rm st}}
\def\s{{\rm s}}
\def\oo{{\scriptstyle{\mathcal O}}}
\def\ooo{{\scriptscriptstyle{\mathcal O}}}
\def\sw{Seiberg-Witten }
\def\pa{\partial_A\bar\partial_A}
\def\gr{{\scriptscriptstyle|}\hskip -4pt{\g}}
\def\subsetint{{\  {\subset}\hskip -2.45mm{\raisebox{.28ex}
{$\scriptscriptstyle\subset$}}\ }}
\def\ra{\rightarrow}
\def\pst{{\rm pst}}
\def\sst{{\rm sst}}
\def\td{{\rm td}}
\def\kod{\mathrm{kod}}
\def\degmax{\mathrm{degmax}}
\def\red{{\rm red}}
\def\reg{{\rm reg}}
\def\ASD{{\rm ASD}}

\title{Instantons and curves on class VII surfaces}
\author{Andrei Teleman}
\date{\today}
\maketitle

\centerline{\it Dedicated to the memory of my father Kostake Teleman} 

\begin{abstract}
We develop a general strategy, based on gauge theoretical methods, to
prove existence of curves  on class VII surfaces. We prove that,
  for $b_2=2$, every minimal class VII surface has a cycle of rational curves hence, by a result of Nakamura,  is a global deformation   of a one parameter family of blown up primary Hopf surfaces.     The case $b_2=1$  was solved in \cite{Te2}. The fundamental object intervening in our strategy is the moduli space ${\mathcal M}^{\pst}(0,{\mathcal K})$ of polystable bundles ${\mathcal E}$ with $c_2({\mathcal E})=0$, $\det({\mathcal E})={\mathcal K}$. For large $b_2$ the geometry of this moduli space becomes very complicated. 
The case
$b_2=2$  treated here in detail requires   new ideas and difficult 
techniques   of both complex geometric and
gauge theoretical nature.  
 We   explain   the substantial  obstacles which must be overcome in order
to extend our methods to the case   $b_2\geq 3$. 
\end{abstract}

\tableofcontents

\setcounter{section}{-1}
\section{Introduction}\label{intro}

The classification problem for class VII surfaces is
a very difficult, still unsolved problem.  Solving this problem would
finally fill the most defying gap in the Enriques-Kodaira classification
table. By analogy with the class of {\it algebraic}
surfaces with $\kod=-\infty$, one expects   this class to be actually {\it very
small}. This idea is supported by the  classification in the case
$b_2=0$, {\it which is known}: any class VII surface with $b_2=0$ is
either a Hopf surface or an Inoue surface \cite{Bo1}, \cite{Bo2}, \cite{Te1}, \cite{LY}.  
On the other hand, solving completely the classification
problem for this class of surfaces  has   been considered for a long  time to
be  a  hopeless  goal; the difficulty comes from the lack of lower dimensional complex geometric objects: for instance, it is not known (and there exists no method to decide) whether a minimal class VII surface  with $b_2>0$ possesses a holomorphic curve, a non-constant entire curve, or a holomorphic foliation.

In his remarkable article \cite{Na2} Nakamura, inspired by the previous work  of Kato (\cite{Ka1}, \cite{Ka2}, \cite{Ka3}), and Dloussky \cite{D1}, stated a  
courageous conjecture, which would in principle solve the classification
problem for class VII surfaces, as we explain below:
\\ \\
{\bf  The GSS conjecture:} {\it Any minimal class VII surface with $b_2>0$
contains a global spherical shell.}
\\

We recall that a (bidimensional) spherical shell is an open surface which
is biholomorphic to a standard neighborhood of $S^3$ in $\C^2$. A global
spherical shell  (GSS) in a surface $X$ is an open submanifold $\Sigma$
of $X$  which is a spherical shell and such that $X\setminus\Sigma$ is
connected. Minimal class VII surfaces which allow GSS's (which are usually called GSS
surfaces, or Kato surfaces) are well understood; in particular it is known that any such
surface   is a degeneration of a holomorphic family of blown up primary
Hopf surfaces, in particular it is diffeomorphic to such a blown up Hopf
surface. Moreover, Kato showed \cite{Ka1}  that  any 
GSS surface can be obtained by a very simple construction: First one considers a modification $m=m_b\circ\dots m_1:\hat D_b\to D$ of the standard disk $D\subset \C^2$, where $m_1:\hat D_1\to D$ is the blowing up at $0\in D$, and $m_k$ is obtained inductively by blowing up in  $\hat D_{k-1}$ a point of its   (-1)-curve. Second, one performs a   {\it holomorphic surgery} $S$ to the
resulting manifold $\hat D_b$ in the following way: one removes a closed ball around a point  $p\in \hat
D_b$ belonging to the last exceptional curve of $m$, and then  identifies
holomorphically the two ends of the resulting manifold (which are both
spherical shells). Choosing in a suitable way the identification map
$s$, one gets a minimal surface. The isomorphism class of the resulting surface  is
determined by  two parameters: the modification $m$   and the identification map $s$.    Note however that Kato's 
simple description
of GSS surfaces   does not immediately yield a clear description of  the
moduli space  of   GSS surfaces, because different pairs
$(m,s)$ can produce isomorphic surfaces.  Nevertheless this shows that, in
principle, 
 the complete classification of GSS surfaces can be obtained with
``classical" methods, so the GSS conjecture would solve in principle the
classification problem for the whole class VII.\\

The existence of a GSS reduces to the existence of ``sufficiently many curves".  This is an
 important progress which is due to several mathematicians (Kato, Nakamura, Dloussky, Dloussky-Oeljeklaus-Toma) who worked on the subject in the last decades. 
More precisely one has
\begin{thry} \label{previousresults}
\begin{enumerate}
\item If a minimal class VII surface $X$ with $b_2(X)>0$ admits $b_2(X)$ rational curves, 
then it also has a global spherical shell. 
\item If a minimal class VII surface $X$ admits a numerically
pluri-anticanonical divisor, i.e. a non-empty curve
$C$ such that 
$c_1({\mathcal O}(C))\in \Z_{\leq 0}\ c_1({\mathcal K})$ mod $\Tors$.
then it also has a global spherical shell.
\item \label{cycleimplies} If a minimal class VII surface $X$ admits  a cycle of curves, then  it is a global deformation 
(a degeneration) of a one parameter family of blown up primary Hopf 
surfaces. 
\end{enumerate}
\end{thry}
Here by ``cycle" we mean either a smooth elliptic curve or a cycle of rational curves (which includes a rational curve with an ordinary double point). The first statement  is the remarkable positive solution -- due to  Dloussky-Oeljeklaus-Toma \cite{DOT} -- of  Kato's
conjecture; this conjecture had been solved earlier in the case $b_2=1$   by Nakamura \cite{Na1}. The second  statement is a recent result of G. Dloussky  \cite{D2}, whereas the third is due to   Nakamura \cite{Na2}. This important theorem shows that, as soon as a minimal class VII surface $X$ with $b_2(X)>0$ admits a cycle, it belongs to the ``known component" of the moduli space.  
 \\

In our previous article \cite{Te2} we proved, using techniques from Donaldson
theory, that any class VII surface with $b_2=1$ has curves; using the results
of Nakamura \cite{Na1} or Dloussky-Oeljeklaus-Toma cited above, this implies that the global spherical shell conjecture  holds in the case $b_2=1$. Since the GSS surfaces in the case $b_2=1$ are very well understood, this solves completely the classification  problem in this case.  The method used in \cite{Te2} can be extended to higher $b_2$, and we believe that, at least for small $b_2$, it should give the existence of a cycle. Our general strategy   has two steps:\\ \\
{\bf Claim 1:} {\it If $X$ is a minimal class VII surface {\it with no cycle} and $b_2(X)>0$, then, for suitable Gauduchon metrics, the moduli space ${\mathcal M}^{\pst}(0,{\mathcal K})$  of polystable bundles ${\mathcal E}$ on $X$ with $c_2=0$ and $\det({\mathcal E})={\mathcal K}$  has a smooth  compact connected component $Y\subset {\mathcal M}^{\st}(0,{\mathcal K})$, which contains a non-empty  finite   subset of  filtrable points.}
\\ 

${\mathcal M}^{\pst}(0,{\mathcal K})$ is endowed with the topology induced by the Kobayashi-Hitchin correspondence from the corresponding moduli spaces of instantons. We will see that this moduli space is always compact (see section \ref{topprop}); this is easy to see for $b_2\leq 3$, because in this case the lower strata in the Uhlenbeck compactification are automatically empty \cite{Te2}. This moduli space is not a complex space, but its stable part ${\mathcal M}^{\st}(0,{\mathcal K})$ is an open subset with a natural complex space structure. 
$Y$ will be defined as the connected component   of the {\it canonical extension} ${\mathcal A}$, which, by definition, is the (essentially unique) non-split extension of the form
\begin{equation}\label{stes}
0\map{\cal K}\textmap{i_0}{\cal A}\textmap{p_0} {\cal O}\map 0\ ,
\end{equation}%
and is stable when $\deg_g({\mathcal K})<0$ and $X$ has no cycle of curves.   The condition  $\deg_g({\mathcal K})<0$ is not restrictive; we will show that there exist Gauduchon metrics with this property. 
  \\ \\ 
{\bf Claim 2:} {\it The existence of such  a component $Y$ leads to a contradiction.} \\  

Both claims might be surprising, and one can wonder how we came to these statements.
 The first claim is an obvious  consequence of the following more precise statement, which has been checked for $b_2\in\{1,2\}$: Consider the subspace ${\mathcal M}^{\st}_\emptyset\subset {\mathcal M}^{\pst}(0,{\mathcal K})$ consisting of those stable bundles which can be written as a line bundle extension whose kernel (left  hand term) has torsion Chern class.  %
\\ \\
{\bf Claim 1$'$:}  {\it If $X$ is a minimal class VII surface {\it with no cycle} and $b_2(X) >0$ then, for a Gauduchon metric $g$ with $\deg_g({\mathcal K})<0$,  the closure $\overline{{\mathcal M}^{\st}_\emptyset}$ of ${\mathcal M}^{\st}_\emptyset$ in ${\mathcal M}^{\pst}(0,{\mathcal K})$  is open in ${\mathcal M}^{\pst}(0,{\mathcal K})$ and contains all  filtrable polystable bundles except the bundles of the form ${\mathcal A}\otimes{\mathcal R}$, ${\mathcal R}^{\otimes 2}={\mathcal O}$.  These bundles   are stable   but do not belong to $\overline{{\mathcal M}^{\st}_\emptyset}$.}
\\

This implies that the  connected component   $Y$ of ${\mathcal A}$ in ${\mathcal M}^{\pst}(0,{\mathcal K})$ is contained in the stable part and is a smooth compact manifold which contains a finite non-empty set of filtrable bundles, so it has the properties stated  in Claim 1. \\

 The main purpose of this article is to show that our 2-step strategy works in the case $b_2=2$.  Therefore, we will prove the following   result: 
 \begin{thry} Any minimal class VII surface with $b_2=2$ has a cycle of curves, so it is a global deformation of a family of blown up Hopf surfaces.
\end{thry}

We explain now in a geometric, non-technical way  how Claim 1$'$  will be proved  for $b_2=2$. We suppose  for simplicity that  $\pi_1(X)\simeq\Z$; in this case the cohomology group $H^2(X,\Z)$ is torsion free and,  by Donaldson's  first theorem \cite{Do2},  is isomorphic to $\Z^{\oplus 2}$ endowed with the standard  negative definite intersection form.  Let $(e_1,e_2)$ be an orthonormal basis of $H^2(X,\Z)$ such that
$$e_1+e_2=-c_1(X)=c_1({\mathcal K})\ .$$

Consider first   the spaces ${\mathcal M}^{\st}_\emptyset$, ${\mathcal M}^{\pst}_\emptyset$  of stable, respectively polystable    extensions  ${\mathcal E}$ of the form 
 \begin{equation}\label{firstext}
 0\to {\mathcal L}\to{\mathcal E}\to{\mathcal K}\otimes{\mathcal L}^\vee\to 0
\end{equation}
with $c_1({\mathcal L})=0$.
It is easy to see that, under our assumptions,  ${\mathcal M}^{\st}_\emptyset$  can be identified with $D^\bullet\times\P^1$, where $D^\bullet\subset  \Pic^0(X)\simeq\C^*$ is the  subset of line bundles satisfying the inequality $\deg_g({\mathcal L})<\frac{1}{2}\deg_g({\mathcal K})$. $D^\bullet$ is a punctured disk.   ${\mathcal M}^{\pst}_\emptyset$ can be identified with the space obtained from the product $\bar D^\bullet\times\P^1$ by   collapsing  to a point each fiber over the circle $\partial \bar D$. The circle   $\Rg'$ of collapsed fibers is one of the two components  of the subspace of reductions (split polystable bundles) ${\mathcal M}^{\red}(0,{\mathcal K})$.  \vspace{1mm} 

One also has two 1-dimensional families of extensions corresponding to the cases   $c_1({\mathcal L})=e_i$. When $X$ has no curve in the classes $\pm(e_1-e_2)$ (assume this for simplicity!), the corresponding loci of polystable  bundles ${\mathcal M}^{\pst}_{\{i\}}$ can be identified with two punctured closed disks $\bar D_i^\bullet:=\bar D_i\setminus\{0_i\}$; the subspaces  ${\mathcal M}^{\st}_{\{i\}}$ of stable bundles of these types   are identified with $D_i^\bullet$. There is a natural isomorphism between the two boundaries $\partial \bar D_i$, and  the points which correspond via this isomorphism represent isomorphic split polystable bundles. Therefore, one has to glue the punctured closed disks $\bar D_i^\bullet$ along their boundaries and get a   2-sphere minus two points $S\setminus\{0_1,0_2\}$, which is the second piece of our moduli space. The circle $\Rg''$ given by the identified boundaries  $\partial D_i$ is the second component   of ${\mathcal M}^{\red}(0,{\mathcal K})$.  
There are two more filtrable bundles in our moduli space, namely  the two bundles of the form ${\mathcal A}\otimes{\mathcal R}$ with ${\mathcal R}^{\otimes 2}={\mathcal O}$. These bundles are stable under the assumption that $\deg_g({\mathcal K})<0$ and $X$ has no cycle. Therefore we also have a 0-dimensional subspace ${\mathcal M}_{\{1,2\}}^\st$ of stable bundles. \vspace{1mm} 

Using another classical construction method for bundles, one gets two more points, namely the push-forwards ${\mathcal B}_1$,  ${\mathcal B}_2$ of two line bundles on a double cover $\tilde X$ of $X$ (see section \ref{bicoverings}). These points are   fixed under the   involution $\otimes\rho$ given by tensoring with the flat $\Z_2$-connection  defined by the generator $\rho$ of $H^1(X,\Z_2)$.
\begin{figure}[h]
\centering
\scalebox{0.5}
{\includegraphics{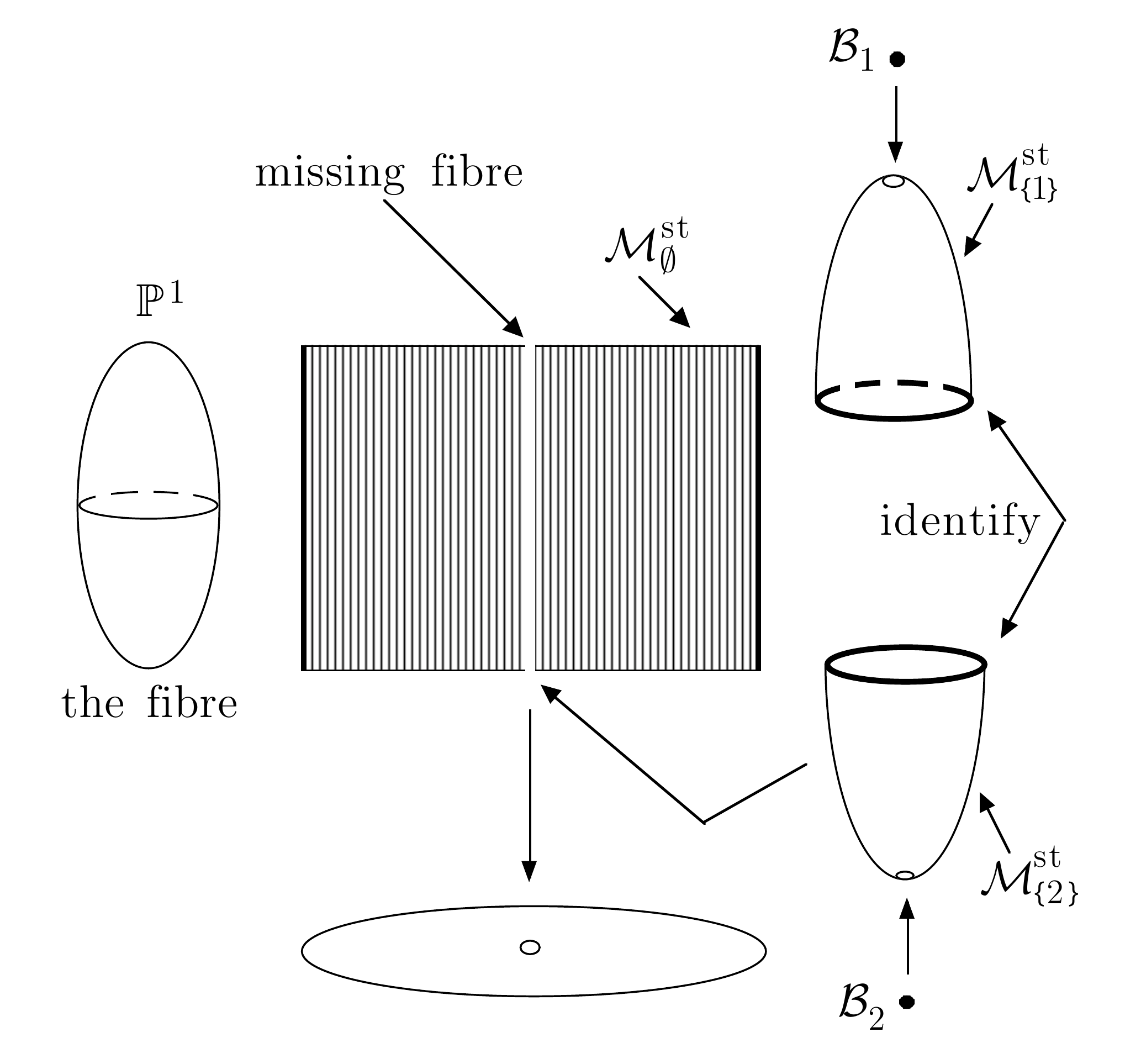}}
\label{picture}
\end{figure}

Under our  assumptions (lack of curves!), the  four pieces  ${\mathcal M}^{\pst}_\emptyset$,  $S\setminus\{0_1,0_2\}$, $\{{\mathcal B}_1\}$,  $\{{\mathcal B}_2\}$ are disjoint. 
 Note now that there is an obvious way to put together these pieces in order to get a compact space (this looks like solving a puzzle!):  one identifies   ${\mathcal B}_1$,  ${\mathcal B}_2$ with the missing points $0_1$, $0_2$ of $S\setminus\{0_1,0_2\}$ and afterwards puts  the obtained sphere $S$ at the place of  the missing fiber of the $\P^1$-fibration ${\mathcal M}^{\st}_\emptyset\to D^\bullet$ (see the {picture} nearby).  The result is a topological space homeomorphic to $S^4$.

Knowing that $\overline{{\mathcal M}^{\st}_\emptyset}$ is obtained in the described way, it will be easy to prove that it is open and that it does not contain any bundle  of the form  ${\mathcal A}\otimes{\mathcal R}$. For the first statement it suffices to compare the local topology of our 4-sphere to the local topology  of ${\mathcal M}^{\pst}(0,{\cal K})$ prescribed by  deformation theory; for the second    it suffices to show that a bundle ${\mathcal A}\otimes{\mathcal R}$ does not belong to any of the four pieces; this follows easily  using again our assumption concerning non-existence of curves.

Therefore, the idea of proving Claim 1$'$ (so also Claim 1) is very clear: solving our puzzle game yields a compact component of the moduli space; the two elements of ${\mathcal M}_{\{1,2\}}$ {\it are not needed} in the construction of this compact component, so they  belong to a new component (or to two new components). However,  the fact  that  our 4-sphere (the  space obtained solving the puzzle game in the most natural way)   is indeed the closure of   ${\mathcal M}^{\st}_\emptyset$ is  difficult to prove. The point is that one has absolutely no control on extensions of the form  (\ref{firstext}) when $\deg_g({\mathcal L})\to -\infty$, because the volume of the section defined by ${\mathcal L}$ in the projective bundle $\P({\mathcal E})$ tends  to $\infty$ as $\deg_g({\mathcal L})\to -\infty$.  In other words, there exists no method to prove that a certain family of extensions  is contained in the closure of another family, so {\it incidence relations between families of extensions are difficult to understand and prove}. This is one of the major difficulties in understanding the global geometry of moduli spaces of bundles on non-K\"ahlerian surfaces. The fact that the above construction gives indeed the closure of ${\mathcal M}^{\st}_\emptyset$ will follow from:
\begin{enumerate}
\item \label{itext} The holomorphic structure of   ${\mathcal M}^{\st}(0,{\mathcal K})$ extends across $\Rg''$ (see sect. \ref{general}) and, with respect to the extended holomorphic structure, the complement ${\mathcal D}$ of ${\mathcal M}^{\pst}_\emptyset$ in ${\mathcal M}^{\pst}(0,{\mathcal K})$ is a divisor. 
\item The circles $\Rg'$ and $\Rg''$ belong to the same component of the moduli space. This important result will be obtained using the Donaldson $\mu$ class associated with a generator of $H_1(X,\Z)/\Tors$ and  a gauge theoretical cobordism argument. 
 \end{enumerate} 
 
 We still have an important detail to explain: Why did we use   the two $\otimes\rho$-fixed points  ${\mathcal B}_1$, ${\mathcal B}_2$ in solving our puzzle (which should produce the closure $\overline{{\mathcal M}^{\st}_\emptyset}$) and  not, for instance, the two filtrable elements of ${\mathcal M}^\st_{\{1,2\}}$ or two non-filtrable bundles? The point is that the involution $\otimes\rho$ acts non-trivially on the divisor ${\mathcal D}$, which  in our simplified case is a projective line. Therefore this divisor  contains the two fixed points of this involution.  \\
 
Note that without assuming $\pi_1(X)\simeq\Z$ the  structure of our moduli space  will a priori be  slightly more complicated. Two complications arise:
\begin{enumerate}
\item    $\overline{{\mathcal M}^{\st}_\emptyset}$ contains a disjoint union of components ${\cal M}_c$ associated with classes $c\in\Tors(H^2(X,\Z))$.  
\item It is difficult to prescribe the   number of fixed points of the involution induced by $\otimes\rho$ on each component (see section \ref{bicoverings}).
\end{enumerate}
However one can easily prove that each component ${\cal M}_c$ can be obtained from the product $\bar D\times\P^1$ by collapsing the fibers over $\partial\bar D$ to  points as above, and applying a (possibly empty) sequence of blow ups in points  lying above the center of the disk. The fiber over this center will be a tree ${\cal D}_c$ of rational curves, one of which, say ${\cal D}^0_c$, contains a circle of reductions and consists of extensions of types $\{1\}$ or $\{2\}$ and two fixed points of $\otimes\rho$.  This will slightly complicate our argument, because we will also have to rule out the case when the bundle ${\cal A}$ belongs to an irreducible component ${\cal D}^1_c\ne {\cal D}^0_c$ of  the tree ${\cal D}_c$. But we will see that such a component (if it exists) consists generically of non-filtrable bundles, so it suffices to apply Corollary 5.3 in \cite{Te2}, which  shows   that there cannot exist a family of rank 2-simple bundles over $X$ parameterized by a closed Riemann surface which contains both filtrable and non-filtrable bundles.
Note that the existence  of a cycle implies $\pi_1(X)=\Z$ by Nakamura's theorem (see Theorem \ref{previousresults}, \ref{cycleimplies}. above)  hence, a posteriori, the moduli space ${\cal M}^\pst(0,{\cal K})$ is  just a 4-sphere, as  explained above.

The first section contains general results concerning moduli spaces of polystable bundles (projectively ASD connections) on non-K\"ahlerian surfaces. These results  will play an important role in the future attempts to solve the GSS conjecture in its full generality.  Particular attention will be given to the  topological properties of moduli spaces of projectively ASD connections on general Riemannian 4-manifolds with $b_+=0$, and to the structure of these moduli spaces around the reduction loci.

The second section is dedicated to the geometry of the moduli space ${\mathcal M}^{\st}(0,{\mathcal K})$ in the case $b_2=2$. In this section we will prove Claim 1$'$ in full generality (without any assumption on $\pi_1(X)$) following the geometric ideas explained above. 
The following sections are dedicated to Claim 2: the appearance of a smooth compact  component in the moduli space leads to a contradiction. This contradiction will be obtained in several steps as follows:
In the third section we will show that the embedding $Y\subset {\cal M}^\st(0,{\mathcal K})$ has a  universal family ${\cal F}\to Y\times X$. This result will enable us in the fourth section  to apply the Grothendieck-Riemann-Roch theorem to the sheaves ${\cal F}$, ${\cal E}nd_0({\cal F})$ and the projection $Y\times X\to Y$; this will give us important information about the Chern classes of the family ${\cal F}$ and  about the Chern classes of  $Y$ itself. The most important result is a parity theorem: the first Chern class of $Y$ is even modulo torsion, i.e. its image in $H^2(Y,\Q)$ belongs to the image of $2H^2(Y,\Z)$. This is a very restrictive condition; it implies for instance that $Y$ is minimal. On the other hand, using the results in \cite{Te2}, we see that $Y$  cannot be covered by curves, so $a(Y)=0$. Therefore, we are left with very few possibilities: a class VII surface with $b_2=0$, a K3 surface, or a torus.  The case when $Y$ is a class VII surface requires a careful examination. This case will be  treated in the fifth section, which contains the final arguments. The other two cases (a K3 surface or a torus) are K\"ahlerian, so they can be ruled out using the results in \cite{Te5}; we will explain briefly the arguments used in \cite{Te5} for completeness.

Therefore we make use essentially of the theory
of surfaces, so it is not clear yet how to  generalize our arguments to larger $b_2$. On the other hand, by the results in \cite{LT1}, the regular part of any  moduli space of stable bundles over a Gauduchon compact manifold is a strong KT (K\"ahler with torsion) manifold.   Therefore, future progress in the classification of this class of manifolds will be very useful  for extending our program to class VII surfaces with  arbitrary $b_2$.
\\
\\
{\bf Acknowledgements:} I have benefited from useful discussions with  many mathematicians, who took their time trying to answer my questions and to follow  my arguments.   I am especially indebted  to Nicholas Buchdahl for his careful and  professional  comments. He  kindly pointed  out to me that the compactness theorem  stated in \cite{Te3} holds for arbitrary $b_2$ (see Theorem \ref{cp}) and came with a short proof of Lemma \ref{degree}.   I learnt a lot about the properties of the ``known" class VII surfaces from Georges Dloussky,  Karl Oeljeklaus and Matei  Toma, who also explained me their recent  result about surfaces with $b_2$ curves.        I also had extensive discussions with Martin L\"ubke  about moduli spaces of holomorphic bundles on non-K\"ahlerian surfaces and their properties.  

I  thank   Simon Donaldson, Richard Thomas, Stefan Bauer and Kim Froyshov  for their interest in my work, their encouragements,  and for giving me the opportunity to give talks about my results on class VII surfaces at Imperial College and Bielefeld University.

\section{General results}\label{strategy}

\subsection{Holomorphic bundles with $c_2=0$, $\det={\cal K}$ on class VII surfaces}
\label{firstresults}

In this section we prove several general results concerning 2-bundles with $c_2=0$, $\det={\cal K}$ on class VII surfaces. 

 Let $X$ be a class VII surface with second Betti number $b$.  Since
$b_2^+(X)=0$, the intersection form
$q_X:H^2(X,\Z)/{\rm Tors}\times H^2(X,\Z)/{\rm Tors}\to \Z$ is  
definite  so, by Donaldson's first  theorem \cite{Do2}, it is trivial over $\Z$.    Put
$$k:=c_1({\mathcal K})=-c_1(X)\ .
$$
 Since $\bar k:=k$ mod Tors is a characteristic element for $q_X$ and $\bar k^2=-b$, it follows easily that there exists a  unique  (up to
order) basis
$(e_1,\dots, e_b)$   in the free $\Z$-module
$H^2(X,\Z)/{\rm Tors}$ such that 
$$e_i\cdot e_j=-\delta_{ij}\ ,\ \bar k=\sum_{i=1}^b e_i\  
 \ .
$$

For instance, when $X$ is a primary Hopf surface blown up at $b$ simple
points, $e_i$ are just the Poincar\'e duals of the exceptional divisors mod $\Tors$. For a
subset  
$$I\subset \{1,\dots,b\}=:I_0\ ,$$
 we put
$$e_I:=\sum_{i\in I} e_i\ ,\ \bar I:=I_0\setminus I\ .
$$

The connected components $\Pic^c$, $c\in H^2(X,\Z)=NS(X,\Z)$ of
the Picard group $\Pic$ of $X$ are isomorphic to $\C^*$. We put  
$$\Pic^T:=\union_{c\in \Tors}\Pic^c\ ,\ \Pic^e:=\union_{c\in
e}\Pic^c\ ,$$
for a class
$e\in H^2(X,\Z)/\Tors$. Let $g$ be a Gauduchon metric  on $X$. We will
use the notations $\Pic^c_{< d}$, $\Pic^e_{<  d}$, $\Pic^T_{< d}$ etc.
for the subspaces of $\Pic^c$, $\Pic^T$, $\Pic^e$ defined by the
inequality $\deg_g({\mathcal L})< d$. Similarly for the subscripts $\leq d$,
$=d$.
\vspace{1mm}

We recall that a holomorphic vector bundle on a compact Gauduchon manifold is called polystable if it is either stable, or it decomposes as a direct sum of stable bundles of the same slope \cite{LT1}. 

Consider the moduli space
${\mathcal M}^{\rm pst}(0,{\mathcal K})$ of holomorphic, $g$-polystable rank 2-bundles ${\mathcal E}$ on $X$ with $c_2({\mathcal E})=0$ and $\det({\mathcal E})={\mathcal K}$.
 The geometry of this moduli space plays a fundamental
role in our arguments. The idea to use this moduli space is surprising      
and might look  artificial; the point is that, whereas for a class VII surface
with no curves the ``classical" complex geometric methods fail,   a lot can be said about the 
corresponding moduli space
${\mathcal M}^{\rm pst}(0,{\mathcal K})$, and the geometry of this space carries important information about the base surface.     
 
The characteristic number $\Delta({\mathcal E}):=
4c_2({\mathcal E})-c_1({\mathcal E})^2$ of a  bundle ${\mathcal E}$ with these
invariants is
$b_2(X)$  and, by the
Riemann Roch theorem, it follows that the expected complex dimension of
the moduli space is also
$b_2(X)$. As explained in \cite{Te2}, this moduli space can be identified with a moduli space of oriented projectively ASD unitary connections via the
Kobayashi-Hitchin correspondence. {\it We will endow this moduli space
with the topology induced by this identification} (see section \ref{topprop} for the main properties of this topology).

One should {\it not} expect this moduli space to be a
complex space: in the non-K\"ahlerian framework, moduli spaces of
instantons have  complicated singularities around the reductions, and
these singularities are {\it not} of a complex geometric nature (see \cite{Te2} and section \ref{topprop} in this article).     
Denote by ${\mathcal M}^{\rm red}(0,{\mathcal K})$ the subspace of reductions (of
split poystable  bundles) in ${\mathcal M}^{\rm pst}(0,{\mathcal
K})$. The open
subspace
$${\mathcal M}^{\rm st}(0,{\mathcal K})={\mathcal M}^{\rm pst}(0,{\mathcal
K})\setminus {\mathcal M}^{\rm red}(0,{\mathcal K})$$ is a complex space
\cite{LT1}. 

It is important to note that ${\mathcal M}^{\rm pst}(0,{\mathcal
K})$ comes with a natural involution.
Indeed, the group
$H^1(X,\Z)/2H^1(X,\Z)\simeq\Z_2$ is a subgroup of
$H^1(X,\Z_2)$, and the latter can be identified  with the group  of flat line
bundles with structure group $\{\pm 1\}\subset S^1$ (see \cite{Te2}).
We denote by $\rho$ the generator of
$H^1(X,\Z)/2H^1(X,\Z)$,   by $\otimes \rho$ the corresponding
involution on ${\mathcal M}^{\rm pst}(0,{\mathcal
K})$, and by ${\mathcal M}^\rho(0,{\mathcal K})$ the fixed point set of this involution. We will see that its points correspond to stable bundles whose pull-back to the double cover $\tilde X_\rho$ associated with $\rho$ are split (see section \ref{bicoverings}).

The   filtrable bundles ${\mathcal E}$ with $c_2=0$, $\det({\mathcal
E})\simeq{\mathcal K}$ can be easily described as extensions. More
precisely, as in
 Proposition 3.2 \cite{Te2} one can show  that 
 \begin{pr}\label{extypes}
Let  ${\mathcal E}$ be a    rank 2-bundle  on $X$ with   $c_2({\mathcal E})=0$ and $\det({\mathcal E})={\mathcal K}$. Then any rank 1 subsheaf ${\mathcal L}$ of ${\mathcal E}$ with torsion free quotient  is a line  {subbundle}  of ${\mathcal E}$ and has $c_1({\mathcal L})\in e_I$   for some
$I\subset I_0$.  In particular, if ${\mathcal E}$ is filtrable,  it is  the central term   of an
extension  
\begin{equation}\label{ext}0\to {\mathcal L}\to {\mathcal E}\to {\mathcal K}\otimes{\mathcal L}^{\vee}\to
0\ , 
\end{equation} 
where $c_1({\mathcal L})\in e_I$   for some
$I\subset I_0$.
\end{pr}

  We recall that an {\it Enoki surface} is a minimal class $VII$ surface with $b_2>0$ which has a   non-trivial effective divisor $D>0$ with  $D\cdot D=0$ (or, equivalently, with $c_1({\cal O}(D))\in\Tors$). By the ``Main Theorem" of \cite{E} it is known that any Enoki surface is an exceptional compactification of an affine line bundle over an elliptic curve and contains a global spherical shell (so also a cycle). Therefore, these surfaces belong to the ``known list", so they are not interesting for our purposes. Recall  also an important vanishing result (see Lemma 1.1.3 \cite{Na3}):  
\begin{lm} \label{nak} On  a minimal class VII surface one has
\begin{equation}
H^0({\mathcal U})=0\ \forall {\mathcal U}\in\Pic(X)\hbox{ \it with }   k \cdot c_1({\mathcal U})<0\ .
\end{equation}
\end{lm}
\pf It suffices to note that   $k \cdot c_1({\cal O}(C))\geq 0$ for every irreducible curve $C$. This follows easily from the genus formula (see \cite{BHPV} p. 85) taking into account that the intersection form of $X$ is negative definite.
\qed

Using   Proposition \ref{extypes} and the vanishing lemma stated above, one gets
 easily the following important regularity result:
 \begin{pr}\label{reg} Let $X$ be a minimal class VII surface with $b_2(X)>0$ which is not an Enoki surface, and let ${\mathcal E}$ be a rank 2-holomorphic bundle on $X$ with $c_2({\mathcal E})=0$, $\det({\mathcal E})={\mathcal K}$. Then $H^2({\mathcal E}nd_0({\mathcal E}))=0$ except when ${\mathcal E}$ is an extension of ${\mathcal K}\otimes{\mathcal R}$ by ${\mathcal R}$, where ${\mathcal R}^{\otimes 2}\simeq{\mathcal O}$.
\end{pr}
\pf  An element $\varphi\in H^0({\mathcal E}nd_0({\mathcal E})\otimes{\mathcal K})\setminus\{0\}$ defines a section $\det(\varphi)\in H^0({\mathcal K}^{\otimes 2})$, and this space   vanishes for class VII surfaces. Therefore $\ker(\varphi)$ is a rank 1 subsheaf of ${\mathcal E}$, so ${\mathcal E}$ is filtrable. By Proposition \ref{extypes}, ${\mathcal E}$ fits into an exact sequence of type (\ref{ext})  with ${\mathcal L}\in\Pic^{e_I}$ for some $I\subset I_0$.

Consider the diagram 
\begin{equation}
\label{morextk}
\begin{array}{ccccccccc}
0&\map&{\mathcal L}&\textmap{\alpha}& {\mathcal E}&\textmap{\beta}&{\mathcal K}\otimes{\mathcal L}^\vee&\map& 0\phantom{\ .}
\\
&&&&\downarrow\varphi
\\
0&\map&{\mathcal K}\otimes{\mathcal L}&\textmap{\id\otimes \alpha}&{\mathcal K}\otimes{\mathcal E} &\textmap{\id\otimes \beta}&{\mathcal K}^{\otimes 2}\otimes{\mathcal L}^\vee&\map& 0\ .
\end{array}
\end{equation}
{\it Case 1.} $(\id\otimes \beta)\circ\varphi\circ\alpha\ne 0$.

The morphism $(\id\otimes \beta)\circ\varphi\circ\alpha$ can be regarded as an element of  $H^0({\mathcal K}^{\otimes 2}\otimes{\mathcal L}^{\otimes-2})$.   Using Lemma \ref{nak} we obtain $H^0({\mathcal K}^{\otimes 2}\otimes{\mathcal L}^{\otimes-2})=0$, except perhaps when $I=I_0$, in which case the Chern class of ${\mathcal K}^{\otimes 2}\otimes{\mathcal L}^{\otimes-2}$ is torsion. But $X$ is not an Enoki surface, so we conclude that in fact  $H^0({\mathcal K}^{\otimes 2}\otimes{\mathcal L}^{\otimes-2})=0$, except  only when ${\mathcal L}^{\otimes 2}={\mathcal K}^{\otimes 2}$. But in this case, any non-trivial morphism ${\mathcal L}\to {\mathcal K}^{\otimes 2}\otimes{\mathcal L}^\vee$ is an isomorphism. Therefore, if   $(\beta\otimes\id)\circ\varphi\circ\alpha$ did not vanish, it would split the second exact sequence, so the first would be also split. This gives ${\mathcal E}={\mathcal L}\oplus ({\mathcal K}\otimes{\mathcal L}^\vee)$, where the second summand is a square root of ${\mathcal O}$.\\
 \\ 
{\it Case 2.} $(\id\otimes \beta)\circ\varphi\circ\alpha= 0$.
  
In this case $\varphi$ maps $\ker(\beta)\simeq{\cal L} $ into $\ker(\id\otimes \beta)\simeq {\cal K}\otimes {\cal L}$. Since $H^0({\mathcal K})=0$, the induced morphism between the two kernels will vanish, hence there exists a well-defined morphism $\psi:{\cal K}\otimes{\cal L}^\vee\to {\cal K}\otimes{\cal E}$ such that  $\varphi=\psi\circ\beta$. The composition $(\id\otimes\beta)\circ\psi$ vanishes (again because $H^0({\mathcal K})=0$), so $\psi$ factorizes as $\psi=(\id\otimes\alpha)\circ\chi$ for a morphism $\chi: {\cal K}\otimes{\cal L}^\vee\to {\cal K}\otimes{\cal L}$, which can be regarded as a section in $H^0({\mathcal L}^{\otimes 2})$. By the same vanishing Lemma \ref{nak}, one has $H^0({\mathcal L}^{\otimes 2})=0$ except when $I=\emptyset$. Since $X$ is not an Enoki surface, $\varphi$ can be non-zero only when ${\mathcal L}^{\otimes 2}\simeq{\mathcal O}$.
\qed

Finally we recall a result proved in \cite{Te3}. This result answers the  question  whether the canonical extension ${\mathcal A}$ can be written as an extension in a different way, and shows that the answer to this question is related to the existence of a cycle in $X$. The result is (see \cite{Te3} Corollary  4.10, Proposition 4.11).
\begin{pr}\label{prevresult} If the bundle ${\mathcal A}$ can be written as   an extension
\begin{equation}\label{diffext}
0\map {\mathcal M}\textmap{i}{\mathcal A}\textmap{p} {\mathcal K}\otimes   {\mathcal M}^{-1}\map 0
\end{equation}
in which the kernel  $\ker(p)\subset {\mathcal A}$ does not coincide with the standard kernel $\ker(p_0)$ of the canonical extension (\ref{stes}), then there exists a non-empty effective divisor $D$ such that:
\begin{enumerate}
\item ${\mathcal M}\simeq {\mathcal O}(-D)$,
\item ${\mathcal K}\otimes {\mathcal O}_D(D)\simeq {\mathcal O}_D$.
\item $c_1({\cal O}(-D))=e_I$  mod Tors  for a subset $I\subset I_0$.
\item $h^0({\mathcal K}\otimes {\mathcal O}_D(D))-h^0({\mathcal K}\otimes {\mathcal O}(D))=1$
\end{enumerate}
Moreover, one of the following holds
\begin{enumerate}
\item $D$ is a cycle,
\item ${\mathcal O}(-D)\simeq{\mathcal K}$ (i.e. $D$ is an anti-canonical divisor).  
 \end{enumerate}
\end{pr}
Note that any anti-canonical divisor contains a cycle (see \cite{Na1} Lemma 12.4). We include a self-contained proof of Proposition \ref{prevresult}  for completeness:\\
\pf 
Since $\ker(p)\subset {\cal A}$ does not coincide with the  kernel $\ker(p_0)\subset {\cal A}$ of the standard exact sequence (\ref{stes})
the composition $p_0\circ i:{\mathcal M}\to{\cal O}$ is non-zero. On the other hand it cannot be an isomorphism because, if it was, $i$ would define a right splitting of (\ref{stes}). Therefore the image of  $p_0\circ i$ is the ideal sheaf of a non-empty effective divisor $D$ and  $p_0\circ i$ defines an isomorphism ${\cal M}\textmap{\simeq} {\cal O}(-D)\subset{\cal O}$, which proves the first  statement. Since all the sheaves in  (\ref{stes}) are locally free, we obtain an exact sequence 
$$0\map {\cal K}_D\textmap{i_0^D} {\cal A}_D\textmap {p_0^D} {\cal O}_D\map 0\ .
$$
of locally free sheaves on $D$. Since $p_0^D\circ i^D=0$,    the restriction $i^D:{\cal M}_D\to{\cal A}_D$ factorizes as $i^D=i_0^D\circ j^D$ for a morphism  $j^D:{\cal M}_D\to {\cal K}_D$.  But $i$ is a bundle embedding, so the induced maps between fibers ${\cal M}(x)\to {\cal A}(x)$ are all injective. This shows that  $j^D(x)\ne 0$ for every $x\in D$, so $j^D$ defines a trivialization of the line bundle $({\cal K}\otimes {\cal M}^\vee)_D\simeq {\cal K}\otimes{\cal O}_D(D)$.  Note that the argument is also valid when $D$ is non-reduced.  This proves  2. 3. follows directly from  Proposition \ref{extypes}. Consider now the commutative diagram
$$
\begin{array}{cccccc}
&&&&&H^0({\cal K}\otimes{\cal O}(D))\\
&&&&&\downarrow\\
&&&&&H^0({\cal K}\otimes{\cal O}_D(D))\\
&&&&&\ \ \downarrow u \\
\map&H^0({\cal A})&\map& H^0({\mathcal O})&\textmap{\partial}& H^1({\cal K})\\
&\downarrow&&\ \ \downarrow a&&\ \ \downarrow v\\
\map&H^0({\cal A}(D))&\map&H^0( {\mathcal O}(D))&\textmap {\partial_D}&H^1({\cal K}\otimes{\cal O}(D))
\end{array}
$$
where the horizontal exact sequences are associated with the short exact sequence (\ref{stes}) and its tensor product  with ${\cal O}(D)$, whereas the vertical exact sequence is associated with the short exact sequence $0\to {\cal K}\to {\cal K}\otimes{\cal O}(D)\to{\cal K}\otimes{\cal O}_D(D)\to 0$.  The morphism $i:{\cal O}(-D)\to {\cal A}$ can be regarded as a lift of the canonical section $s=a(1)$ to $H^0({\cal A}(D))$. Therefore $\partial_D(a(1))=0$, so $v(\partial(1))=0$. But $\partial(1)\in H^1({\cal K})\simeq\C$ is the extension invariant  of the extension (\ref{stes}), which  is non-zero by the definition of this extension. The vertical exact sequence yields an exact sequence
$$0=H^0({\cal K})\map H^0({\cal K}\otimes{\cal O}(D))\map H^0({\cal K}\otimes{\cal O}_D(D))\map H^1({\cal K})\simeq\C\map 0\ ,
$$
which proves 4. Taking into account that $a(X)=0$ we have $h^0({\cal K}\otimes{\cal O}(D))\leq 1$, so there are two possibilities:\\ \\
1. $h^0({\cal K}\otimes{\cal O}(D))=0$, $h^0({\cal K}\otimes{\cal O}_D(D))=h^1({\cal O}_D)=1$.

Using Lemma 2.7 in \cite{Na1} and the obvious inequality $h^1({\cal O}_{D_{\rm red}})\leq h^1({\cal O}_{D})$ (the canonical map $H^1({\cal O}_D)\to H^1({\cal O}_{D_{\rm red}})$ is surjective),  we get $h^1({\cal O}_{D_{\rm red}})=1$.  Let $0<C\leq D_{\rm red}$ be a minimal divisor such that $h^1({\cal O}_{C})=1$. By Lemma 2.3 in \cite{Na1} the divisor $C$ is a cycle. Write $D=C+E$ with $E\geq 0$. Denoting ${\cal N}:={\cal K}\otimes{\cal O}(D)$ and noting that $h^2({\cal N})=h^0({\cal O}(-D))=0$ we get an exact sequence
$$0\to H^0({\cal N}(-C))\to H^0({\cal N})\to H^0({\cal N}_C)\to H^1({\cal N}(-C))\to H^1({\cal N})\to $$
$$\to H^1({\cal N}_C)\to H^2({\cal N}(-C))\to 0
$$
We know that ${\cal N}$ is trivial on $D$, so it is also trivial on $C\subset D$. Therefore $h^1({\cal N}_C)=h^1({\cal O}_C)=1$.  Note that $c_1({\cal N})=e_{\bar I}$ (where $\bar I:=I_0\setminus I$), so $\chi({\cal N})=0$ by the Riemann-Roch theorem. We have assumed  $h^0({\cal N})=0$ and we know $h^2({\cal N})=0$, so $h^1({\cal N})=0$. Therefore $h^2({\cal N}(-C))=h^1({\cal O}_C)=1$. But $h^2({\cal N}(-C))=h^0({\cal K}\otimes{\cal N}^\vee(C))=h^0(-E)$. This shows that $E=0$, so $D$ coincides with the cycle $C$.\\

2. $h^0({\cal K}\otimes{\cal O}(D))=1$, $h^0({\cal K}\otimes{\cal O}_D(D))=h^1({\cal O}_D)=2$.

Since $c_1({\cal K}\otimes{\cal O}(D))=e_{\bar I}$ and   $H^0({\cal K}\otimes{\cal O}(D))\ne 0$, we get by the vanishing Lemma \ref{nak} that $\bar I=\emptyset$, so ${\cal K}\otimes{\cal O}(D)$ is a flat line bundle. When $X$ is not an Enoki surface it follows already that ${\cal O}(-D)\simeq{\cal K}$. When $X$ is an Enoki surface, it follows easily that in fact $X$ is a parabolic Inoue surface and $D=E+C$, where $E$ is  the elliptic curve of $X$ and $C$ is its numerically trivial cycle; it is well-known that this sum is anti-canonical (see \cite{Te3} for details). 
\qed

For every $I\subset I_0$ we have a family  ${\mathcal
F}_I$ of extensions; the elements of ${\mathcal F}_I$ -- which will be called {\it extensions of type $I$} --  are in 1-1 
 correspondence with pairs
$({\mathcal L},\varepsilon)$, where  ${\mathcal L}\in\Pic^{e_I}$   and
$\varepsilon\in H^1({\mathcal L}^{\otimes 2}\otimes{\mathcal K}^{\vee})$, so
${\mathcal F}_I$ is naturally a linear space over $\Pic^{e_I}$. We will denote by ${\mathcal E}({\mathcal L},\varepsilon)$ the central term of the extension associated with the pair  $({\mathcal L},\varepsilon)$. For ${\mathcal L}\in\Pic^{e_I}$ one has
$$\chi({\mathcal L}^{\otimes 2}\otimes{\mathcal K}^{\vee})=\frac{1}{2}(e_I-
e_{\bar I})(-2 e_{\bar I})=-(b-| I|)\ (b:=b_2(X)),
$$
so the dimension of the generic fiber of this linear space is $b-| I|$; the
dimension of the fiber $H^1({\mathcal L}^{\otimes 2}\otimes{\mathcal K}^{\vee})$ 
jumps  when  $h^0({\mathcal L}^{\otimes 2}\otimes{\mathcal K}^{\vee})$ or
$h^0({\mathcal L}^{\otimes -2}\otimes{\mathcal K}^{\otimes 2})>0$.
It might happen that the same bundle can be written as extension in
many ways, so in general the loci ${\mathcal M}_I^{\rm st}$ of stable bundles     
defined by the elements of ${\mathcal F}_I$ might have intersection points. 

Therefore, it is important to have general rules to decide whether two different extensions have
 isomorphic central terms. This problem will be addressed in the following section in a general framework.

\subsection{Morphisms of extensions}
\label{morextsect}
  
  In this section we will address the following questions:
  \begin{itemize}
   \item Under which conditions are the central terms of two different  line bundle extensions isomorphic?
   \item Is the central term of a given non-trivial line bundle  extension simple?
   \end{itemize}
  
Let ${\mathcal L}'$, ${\mathcal L}''$, ${\mathcal M}'$, ${\mathcal M}''$ line bundles on a compact manifold $X$. Consider a diagram of the form
\begin{equation}
\label{morext}
\begin{array}{ccccccccc}
0&\map&{\mathcal L'}&\textmap{\alpha'}& {\mathcal E}'&\textmap{\beta'}&{\mathcal M}'&\map& 0
\\
&&&&\downarrow\varphi
\\
0&\map&{\mathcal L''}&\textmap{\alpha''}& {\mathcal E}''
&\textmap{\beta''}&{\mathcal M}''&\map& 0
\end{array}
\end{equation}
with exact lines.

\begin{pr} \label{prext} Suppose that  $\beta''\circ\varphi\circ  \alpha'=0$. Then  
\begin{enumerate}
\item    There exist  morphisms $u:{\mathcal L}'\to{\mathcal L}''$, $v:{\mathcal M}'\to{\mathcal M}''$ making commutative the diagram (\ref{morext}).
\item If $\varphi:{\mathcal E}'\to {\mathcal E}''$ is an isomorphism, then $u$ and $v$ are   isomorphisms.
\item If $H^0({{\mathcal L}'}^\vee\otimes {\mathcal L''})=0$  then   $\varphi:{\mathcal E}'\to {\mathcal E}''$   is induced by a morphism ${\mathcal M}'\to {\mathcal E}''$, so it cannot be an isomorphism.  
\item If  $H^0({{\mathcal M}'}^\vee\otimes {\mathcal M''})=0$, then
  $\varphi:{\mathcal E}'\to {\mathcal E}''$   is induced by a morphism  ${\mathcal E}'\to {\mathcal L}''$, so it cannot be an isomorphism. 
\item If $H^0({{\mathcal L}'}^\vee\otimes {\mathcal L''})=0$ and $H^0({{\mathcal M}'}^\vee\otimes {\mathcal M''})=0$,  then any morphism $\varphi:{\mathcal E}'\to {\mathcal E}''$ is induced by a morphism ${\mathcal M}'\to {\mathcal L''}$. 
 \end{enumerate}
\end{pr}
\pf 

1.  Since
$\beta''\circ\varphi\circ  \alpha'=0$,
$\varphi$ maps
${\mathcal L}'$ to
${\mathcal L}''$ (defining a morphism $u:{\mathcal L}'\to{\mathcal L}''$) and induces a morphism $v:{\mathcal M}'\to {\mathcal M''}$. 
\\

2. It is easy to show that, when $\varphi$ is an isomorphism, $u$ will be a monomorphism and $v$ an epimorphism. But any epimorphism of locally free rank 1 sheaves is an isomorphism.  Diagram chasing shows that $u$ is also surjective.\\ 

3.  Suppose that $H^0({{\mathcal L}'}^\vee\otimes {\mathcal L''})=0$. In this case $u=0$, so $\varphi$
vanishes on ${\mathcal L}'$, hence it is induced by a morphism $\nu:{\mathcal M}'\to {\mathcal E}''$.\\

4.  Suppose that $H^0({{\mathcal M}'}^\vee\otimes {\mathcal M}'')=0$. In this case $v=0$, so the image of $\varphi$ is contained in $\ker(\beta'')={\mathcal L}''$, hence $\varphi$ is induced
 by a morphism $\mu:{\mathcal E}'\to {\mathcal L}''$.\\

5. Suppose   that  	$H^0({{\mathcal L}'}^\vee\otimes {\mathcal L''})=H^0({{\mathcal M}'}^\vee\otimes {\mathcal M}'')=0$. Then $\beta''\circ\nu=0$, hence $\im(\nu)$ is contained in ${\mathcal L}''$, proving that
$\varphi$ is induced by a morphism ${\mathcal M}'\to {\mathcal L''}$.
\qed
\begin{lm}\label{nontrext} If ${\mathcal L}'\simeq{\mathcal M}''$ and  the second exact sequence
is non-trivial, then one has always $\beta''\circ\varphi\circ  \alpha'=0$ hence the conclusions of Proposition \ref{prext} hold.
\end{lm}
\pf We can suppose ${\mathcal L}'={\mathcal M}''$. If $\beta''\circ\varphi\circ  \alpha'\ne 0$, the composition  $\beta''\circ\varphi\circ  \alpha'$ would be an isomorphism, hence a suitable scalar multiple of  $\varphi\circ\alpha'$   would split the second exact sequence.
\qed
\begin{co}\label{injectivity} [extensions with isomorphic determinants  and  central terms]
Suppose that  in (\ref{morext}) one has ${\mathcal L}'\otimes{\mathcal M}'\simeq{\mathcal L}''\otimes{\mathcal M}''\simeq{\mathcal D}$ and $\varphi:{\cal E}'\to {\cal E}''$ is an isomorphism. Denote by  
 $\varepsilon'\in H^1({{\mathcal M}'}^\vee\otimes{\mathcal L}')$, $\varepsilon''\in H^1({{\mathcal M}''}^\vee\otimes{{\mathcal L}''})$ the invariants associated with the two extensions.  Then one of the following holds:
\begin{enumerate}
 \item There exists isomorphisms $u:{\cal L}'\to {\cal L}''$, $v:{\cal M}'\to {\cal M}''$ making (\ref{morext}) commutative. In this case $\varepsilon''=H^1(w)(\varepsilon')$ where $w:{{\mathcal M}'}^\vee\otimes{\mathcal L}'\to {{\mathcal M}''}^\vee\otimes{\mathcal L}''$ is the isomorphism induced by the pair $(u,v)$.
 \item   ${\mathcal D}\otimes{{\mathcal L}'}^\vee\otimes{{\mathcal L}''}^\vee$ is not trivial  and $H^0({\mathcal D}\otimes{{\mathcal L}'}^\vee\otimes{{\mathcal L}''}^\vee)\ne 0$.
 \item  ${\mathcal D}\otimes{{\mathcal L}'}^\vee\otimes{{\mathcal L}''}^\vee$ is trivial and $\varepsilon'=\varepsilon''=0$. 
\end{enumerate} 
\end{co}
\pf  When $\beta''\circ\varphi\circ\alpha'=0$, the statement 1. will hold by  Proposition \ref{prext}. When $\beta''\circ\varphi\circ\alpha'\ne 0$ and ${\cal M}''\otimes {{\cal L}'}^\vee\simeq {\cal D}\otimes {{\cal L}'}^\vee\otimes {{\cal L}''}^\vee$ is not trivial, 2. holds obviously. Finally, when  $\beta''\circ\varphi\circ\alpha'\ne 0$ and   ${\cal M}''\otimes {{\cal L}'}^\vee$ is trivial, we see using Lemma \ref{nontrext} that the second exact sequence  splits, hence $\varepsilon''=0$. Changing the roles  and noting that ${\cal M}''\otimes {{\cal L}'}^\vee\simeq{\cal M}'\otimes {{\cal L}''}^\vee$ we obtain $\varepsilon'=0$.
\qed
\begin{co} \label{inject} Let ${\mathcal L}$, ${\mathcal M}$ be two line bundles on $X$.
\begin{enumerate}
\item Denote by ${\mathcal E}'$, ${\mathcal E}''$ the middle terms of the extensions associated with $\varepsilon'$, $\varepsilon''\in H^1({\mathcal M}^\vee\otimes {\mathcal L})$.
Suppose $H^0({\mathcal L}^\vee\otimes{\mathcal M})=0$. Then
  ${\mathcal E}'  \simeq{\mathcal E}''$ if and only if $\varepsilon',\ \varepsilon'$  are conjugate modulo
$\C^*$.
\item Suppose   $H^0({\mathcal L}^\vee\otimes{\mathcal M})=0$, $H^0({\mathcal M}^\vee\otimes{\mathcal L})=0$. Then  any nontrivial extension of ${\mathcal M}$ by ${\mathcal L}$ is simple.
 
\end{enumerate}
\end{co}
\pf  
The first statement is a particular case of Corollary \ref{injectivity}, 1. 
For the second, let 
$$0\map{\cal L}\textmap{\alpha}{\mathcal E}\textmap{\beta}{\cal M}\map 0
$$
be a nontrivial extension, and $\varphi:{\mathcal E}\to{\mathcal E}$   a bundle morphism. By  Proposition \ref{prext} (1) and the first assumption there exist morphisms $u:{\cal L}\to{\cal L}$,  $v:{\mathcal M}\to{\mathcal M}$ such that $\varphi\circ\alpha=\alpha\circ u$, $v\circ \beta=\beta\circ\varphi$. Since ${\cal M}$ is a line bundle, we can write $v=\zeta\id_{{\mathcal M}}$, for a constant $\zeta\in\C$. The endomorphism  $\psi:=\varphi- \zeta\id_{\cal E}$ has the property
$$\beta\circ \psi=\beta \circ(\varphi- \zeta\id_{\cal E})=v\circ\beta- \zeta\beta=0\ ,
$$
so it factorizes as $\psi=\alpha\circ\mu$ for a morphism $\mu:{\mathcal E}\to {\mathcal L}$. Note that $\mu\circ \alpha:{\mathcal
L}\to {\mathcal L}$ cannot be an isomorphism  because, if it was, a suitable scalar multiple of $\mu$ would split  the first exact sequence. Therefore $\mu\circ\alpha=0$, so $\mu$ vanishes on $\alpha({\cal L})$, which shows that it is induced by a morphism  ${\mathcal M} \to {\mathcal L}$. This implies $\mu=0$  by our second assumption.
\qed
\subsection{Stable bundles defined by line bundles on double covers}\label{bicoverings}
Let $(X,g)$ be a compact Gauduchon manifold and $\rho\in H^1(X,\Z_2)\setminus\{0\}$.
Let
 $\pi_\rho:\tilde X_\rho\ra X$ 
be  the corresponding double cover,  $\iota$ its canonical
involution, and   
${\cal L}_\rho$ the  holomorphic line bundle on $X$ defined by $\rho$ (regarded as a
representation
$\pi_1(X)\to\{\pm 1\}\subset\C^*$). A push-forward 2-bundle (i.e. a bundle of the form
$(\pi_\rho)_*({\cal M})$, where ${\cal M}$ is a holomorphic line bundle on $\tilde X_\rho$)
is always polystable. Indeed, choosing a Hermitian-Einstein metric $h$ on ${\cal M}$ with
respect to  the $\iota$-invariant metric $\pi_\rho^*(g)$, one gets a Hermitian-Einstein
metric $\iota^*(h)$ on $\iota^*({\cal M})$ (with the same Einstein constant as $h$), hence a Hermitian-Einstein metric $h\oplus\iota^*(h)$ on ${\cal M}\oplus \iota^*({\cal
M})=\pi_\rho^*((\pi_\rho)_*({\cal M}))$  which descends to
$(\pi_\rho)_*({\cal M})$. Conversely, one has the well known

\begin{thry}\label{push} Let ${\cal E}$  be  a  stable rank 2-bundle on $X$, such
that
${\cal E}\otimes{\cal L}_\rho\simeq{\cal E}$. 
Then there exists a line bundle ${\cal M}$ on $\tilde X_\rho$ such that
\begin{enumerate}
\item $[\pi_\rho]_*({\cal M})\simeq {\cal E}$.
\item $\pi_\rho^*({\cal E})={\cal M}\oplus
\iota^*({\cal M})$.
\end{enumerate}
\end{thry}
\pf 
Consider the sheaf of ${\cal O}_X$-algebras
${\cal A}:=[\pi_\rho]_*({\cal O}_{\tilde X_\rho})={\cal O}_X\oplus{\cal
L}_\rho$.  The category of  ${\cal O}_{\tilde X_\rho}$-modules  is equivalent (via the functor 
 $[\pi_\rho]_*$) with the category of 
${\cal A}$-modules. Composing with the obvious automorphism $(u,v)\mapsto (u,-v)$
of ${\cal A}$ one gets a functor $(\cdot)'$ on the category of ${\cal
A}$-modules, which corresponds to the functor   $\iota^*(\cdot)$ on the category of ${\cal
O}_{\tilde X_\rho}$-modules.  Let now $f:{\cal E}\otimes{\cal L}_\rho\to{\cal E}$ be an
isomorphism. The composition 
$$f\circ (f\otimes\id_{{\cal L}_\rho}):{\cal E}\otimes{\cal L}_\rho^{\otimes 2}={\cal
E}\map{\cal E}$$
is an automorphism of ${\cal E}$, which is stable, hence simple. Multiplying $f$ by a
constant if necessary, we can assume that this composition is $\id_{\cal E}$. Therefore,
$f$ defines the structure of  a locally free rank 1 ${\cal A}$-module ${\cal M}$ on  the
sheaf associated with 
${\cal E}$. This proves 1. The bundle $\pi_\rho^*({\cal E})$ corresponds to the ${\cal A}$-module ${\cal
E}\otimes_{{\cal O}_X}{\cal A}$, which can be easily identified with ${\cal M}\oplus {\cal
M}'$. So $\pi_\rho^*({\cal E})\simeq {\cal M}\oplus
\iota^*({\cal M})$.
\qed

Suppose now that $X$ is a class VII surface with $\pi_1(X)\simeq \Z$, $b_2(X)=b$, let $\rho$ be  the non-trivial element of $H^1(X,\Z_2)\simeq\Z_2$, and $\pi:\tilde X\to X$ the corresponding double cover.   $\tilde X$ will also be a class VII surface with $\pi_1(\tilde X)\simeq \Z$ and, comparing the Euler characteristics, we get $b_2(\tilde X)=2b$. We denote by ${\cal K}$, $\tilde{\cal K}$ the corresponding canonical line bundles and by $k$, $\tilde k=\pi^*_\rho(k)$ their Chern classes. By the coefficients formula the groups $H^2(X,\Z)$, $H^2(\tilde X,\Z)$ are both torsion free. As explained  in section \ref{firstresults} we can write $k=\sum_{i=1}^{b} e_i$, $\tilde k=\sum_{s=1}^{2b} f_s$ where $e_i\cdot e_j=-\delta_{ij}$,  $f_s\cdot f_t=-\delta_{st}$. Put $\tilde e_i:=\pi^*(e_i)$.    Using the formulae 

$$\tilde k=\sum_i\tilde e_i=\sum_s f_s\ ,\ \tilde e_i\cdot \tilde e_j=-2\delta_{i,j}\ ,\  \tilde e_i\cdot \tilde k=-2\ ,$$
 we see that there exists a partition $\{1,\dots,2b\}=\cup_{i=1}^n  J_i$ in subsets $J_i$ with $|J_i|=2$ such that
$\tilde e_i=\sum_{s\in J_i} f_s$. Writing $J_i=\{e_i',e_i''\}$ we obtain an orthonormal basis   $(e_1',e_1'',\dots,e_b',e_b'')$ of  $H^2(\tilde X,\Z)$ such that $\tilde e_i=e_i'+e_i''$. The morphism $\iota^*$ leaves invariant each set $\{e_i',e_i''\}$. The  Lefschetz fixed point theorem gives $\tr (\iota^*:H^2(\tilde X)\to H^2(\tilde X))=0$, which implies $\iota^*(e_i')=e_i''$ for any $i$.
\begin{pr} \label{fixedpoints} Let $X$ be a class VII surface with $\pi_1(X)\simeq\Z$. Then
\begin{enumerate}
\item A line bundle ${\cal M}\in\Pic(\tilde X)$ satisfies $c_1([\pi_\rho]_*({\cal M}))=k$, $c_2([\pi_\rho]_*({\cal M}))=0$ if and only if there exists $I\subset\{1,\dots,b\}$ such that $c_1({\cal M})=e'_I+e''_{\bar I}$.
\item For every $I\subset \{1,\dots,b\}$ there exists a unique line bundle ${\cal M}_I\in \Pic^{e'_I+e''_{\bar I}}(\tilde X)$ such that $\det([\pi_\rho]_*({\cal M}_I))\simeq{\cal K}$, $c_2([\pi_\rho]_*({\cal M}_I))=0$.
\item  The involution $\otimes\rho$ on the moduli space ${\cal M}^\st(0,{\cal K})$ has $2^{b-1}$  fixed points.
\end{enumerate}
\end{pr}
\pf  Write $c_1({\cal M})=:c=\sum_i n'_i e'_i+\sum_j n''_j e''_j$. ${\cal M}$ satisfies the conditions in 1. if and only if  $c\cup(\iota^*c)=0$ and $c+\iota^*(c)=\tilde k$. This is equivalent to the system
$$n'_i+n''_i=1\ \forall i\in\{1,\dots,n\},\ \sum_in'_in''_i=0\ . $$
We get $\sum_in'_i(1-n'_i)=0$, hence $n_i'\in\{0,1\}$.   This proves the first statement. For the second,  the pull-back morphism $\Pic^0(X)\to\Pic^0(\tilde X)$ is surjective and its kernel is generated by the non-trivial square root ${\cal R}$ of the trivial line bundle ${\cal O}$. Fix   ${\cal M}_0\in \Pic^{e'_I+e''_{\bar I}}(\tilde X)$. Since
$$\det([\pi_\rho]_*({\cal M}_0\otimes\pi^*({\cal L})))\simeq \det((\pi_\rho)_*({\cal M}_0)\otimes {\cal L})\simeq \det([\pi_\rho]_*({\cal M}_0))\otimes {\cal L}^2\ ,  $$
we find two flat line bundles ${\cal L}_1$, ${\cal L}_2={\cal L}_1\otimes{\cal R}\in\Pic^0(X)$ for which  it holds  $\det([\pi_\rho]_*({\cal M}_0\otimes\pi^*({\cal L_i})))\simeq{\cal K}$. Since $\pi^*({\cal R})$ is trivial, the claim is proved. For the third statement, note that $[\pi_\rho]_*({\cal M}_I)$ is stable and that $[\pi_\rho]_*({\cal M}_I)\simeq[\pi_\rho]_*({\cal M}_J)$ if and only if either $J=I$ or $J=\bar I$. This follows easily by taking the pull-back of the bundles $[\pi_\rho]_*({\cal M}_I)$ to $\tilde X$.
\qed

Note that without assuming $\pi_1(X)\simeq\Z$, counting the fixed points of the involution $\otimes\rho:{\cal M}^\st(0,{\cal K})\to{\cal M}^\st(0,{\cal K})$ becomes very difficult.

\subsection{The moduli spaces ${\cal M}^\pst_{\cal D}(E)$, ${\cal M}_a(E)^\ASD$. Topological properties}\label{topprop}

The results obtained so far   allow us to describe certain pieces of our moduli space: spaces of reductions, spaces of extensions and spaces of fixed points under the natural involution $\otimes\rho$. In order to understand how these pieces fit together, we need several important general results about the topology of the moduli spaces of polystable rank 2-bundles on non-K\"ahlerian complex surfaces.  
 
Let $(E,h)$ be a Hermitian  rank 2-bundle with $c_2(E)=c$ on a compact Gauduchon surface 
$(X,g)$. We fix a holomorphic structure ${\cal D}$ on the determinant line bundle $D:=\det(E)$ and denote by ${\cal M}^\pst_{\cal D}(E)$ the moduli space of polystable holomorphic structures ${\cal E}$ on $E$ with $\det({\cal E})={\cal D}$ (of polystable ${\cal D}$-{\it oriented holomorphic structures}). Two such structures are considered equivalent if they are equivalent  modulo the action of  complex gauge group  ${\cal G}^\C:=\Gamma(X, \SL(E))$.  This moduli space plays a  fundamental role in this article. Its points correspond bijectively to   isomorphism classes of polystable holomorphic bundles ${\cal F}$ with $c_2({\cal F})=c$, $\det({\cal F})\simeq{\cal D}$, so  we will also use the notation ${\cal M}^\pst(c,{\cal D})$ (used more frequently in the complex geometric literature) for this space.

Denote by $a\in{\cal A}(\det(E))$ the Chern connection of the pair $({\cal D},\det(h))$. As in \cite{Te2} we denote by ${\cal A}_a(E)$ (respectively ${\cal A}_a(E)^*$) the space of (irreducible) $a$-oriented Hermitian connections, i.e. the space of (irreducible) Hermitian connections $A$ on $E$ with $\det(A)=a$.  We denote as usually
$${\cal B}_a(E):=\qmod{{\cal A}_a(E)}{{\cal G}}\ ,\ {\cal B}_a(E)^*:=\qmod{{\cal A}_a(E)^*}{{\cal G}}
$$
the infinite dimensional quotients of these spaces by the (real) gauge group ${\cal G}:=\Gamma(X,\SU(E))$ of $(E,h)$. The latter quotient becomes a Banach manifold after suitable Sobolev completions \cite{DK}.  It is convenient to use the  action of the quotient group $\bar {\cal G}:={\cal G}/\{\pm \id_E\}$  which is effective on ${\cal A}_a(E)$ and has trivial stabilizers at the irreducible connections.  

Let ${\cal A}^{\rm ASD}_a(E)$ (respectively ${\cal A}^{\rm ASD}_a(E)^*$) be the subspace of (irreducible) solutions of the {\it projectively ASD equation}

$$(F_A^0)^+=0\ , \eqno{(\ASD)}
$$
and  ${\cal M}^{\rm ASD}_a(E):={{\cal A}^{\rm ASD}_a(E)}/{{\cal G}}$, ${\cal M}^{\rm ASD}_a(E)^*:={{\cal A}^{\rm ASD}_a(E)^*}/{{\cal G}}$ the corresponding moduli spaces. In this formula ${\cal G}:=\Gamma(X,\SU(E))$ denotes  the (real) gauge group of $(E,h)$.  

As explained in the introduction and in \cite{Te2}, {\it ${\cal M}^\pst_{\cal D}(E)$ is endowed with the topology which makes the bijection ${\cal M}^{\rm ASD}_a(E)\to  {\cal M}^\pst_{\cal D}(E)$ given by the Kobayashi-Hitchin correspondence (\cite{Do1}, \cite{Te2} \cite{Bu1}, \cite{LY},  \cite{LT1}) a homeomorphism. } Via this correspondence ${\cal M}^\st_{\cal D}(E)$ corresponds to ${\cal M}^{\rm ASD}_a(E)^*$, so it is open in ${\cal M}^\pst_{\cal D}(E)$.  ${\cal M}^\st_{\cal D}(E)$ can also be identified with an open subspace of the moduli space ${\cal M}^\s_{\cal D}(E)$ of simple ${\cal D}$-oriented holomorphic structures on $E$, so it has a natural complex space structure  inherited from ${\cal M}^\s_{\cal D}(E)$. In general this structure does not extend to  ${\cal M}^\pst_{\cal D}(E)$. Understanding the local structure of ${\cal M}^\pst_{\cal D}(E)$ around the split polystable bundles in the non-K\"ahlerian framework is a difficult task.   The difficulty is the following: whereas the germ of ${\cal M}^\st_{\cal D}(E)$ at a {\it stable} point ${\cal E}\in {\cal M}^\st_{\cal D}(E)$ can be identified with the  universal deformation of ${\cal E}$ in the sense of holomorphic deformation theory, the local  structure at a split polystable bundle (a reduction) cannot be understood using only complex geometric methods.  For instance, the moduli space described in \cite{Te2} contains a finite union of circles of split polystable bundles; although any such bundle ${\cal E}$ has $H^2({\cal E}nd_0({\cal E}))=0$ (so it satisfies the naive complex geometric ``regularity condition") it gives a boundary point in the moduli space (which is a union of compact disks bounded by circles of ``regular" reductions).  

In this section we will discuss topological properties of the moduli spaces ${\cal M}^\pst_{\cal D}(E)$ on non-K\"ahlerian surfaces,  namely: 
\begin{enumerate}
\item compactness properties,
\item  the structure of such moduli spaces around reduction loci,
\item  the restriction of the Donaldson $\mu$-classes to a boundary of a standard neighborhood of a reduction locus.
\end{enumerate}

Compactness has  already been discussed in \cite{Te2}. Here we will give a more general result due to Nicholas Buchdahl. The structure of the ASD instanton moduli spaces around reduction loci and the behavior of the Donaldson $\mu$-classes around these loci have been extensively studied in \cite{Te4} with gauge theoretical methods. For completeness we will give here short, self-contained  proofs of the results we need.  We will not make use  of the general (but difficult) results about normal neighborhoods of reduction loci in ${\cal B}_a$ \cite{Te4}; instead, in sections \ref{localstructure}, \ref{donclass} we will ``blow up" the reduction loci in a Donaldson moduli space, and we will show that the Donaldson $\mu$-classes extend to the blow up. As explained in the introduction, we are  interested in the classes  $\mu(\gamma)$ associated with elements $\gamma\in H_1(X,\Z)/\Tors$. The result we need states that, if $X$ is a class VII surface with $b_2(X)=2$ and $\gamma$ is a  generator of $H_1(X,\Z)/\Tors$, then the restriction of $\mu(\gamma)$ to the (suitably oriented) boundary of a standard neighborhood of a circle of reductions is the fundamental class of this boundary. This result will give us important information about the position of the circles of reductions in the moduli space (see Proposition \ref{firststr}).

The results in sections \ref{localstructure}, \ref{donclass} hold under the assumption that the reductions are regular, i.e. that the second cohomology spaces of their deformation elliptic complexes vanish. It is important to have a complex geometric criterion for this regularity condition, which can be checked for instantons associated with split polystable bundles.  This problem will be addressed in section \ref{compell}, in which we will compare the deformation  elliptic complex  of a split polystable bundle to the deformation elliptic complex of the corresponding reducible instanton. The main result states that, on non-K\"ahlerian surfaces, the second cohomology spaces of the two complexes can be identified. Note that this result {\it does not hold in the K\"ahlerian case}: on a K\"ahlerian surface a reducible instanton has always non-vanishing second cohomology (so it can never be regular), see Corollary \ref{comcoh}.  

In section \ref{general} we will apply our general results  to the moduli space ${\cal M}^\pst(0,{\cal K})$ on a class VII surface with $b_2=2$, and we will show that (for suitable Gauduchon metrics) this moduli space is a topological 4-manifold and that the complex structure of ${\cal M}^\st(0,{\cal K})$ is smooth and extends smoothly across a part of the reduction locus.

\subsubsection{Compactness properties}
Endowed with the topology induced by the Koba\-ya\-shi-Hitchin correspondence ${\cal M}^\pst_{\cal D}(E)$ has the following important properties  inherited from ${\cal M}^{\rm ASD}_a(E)$ (see \cite{DK}):
\begin{enumerate}
\item It is Hausdorff in all cases,
\item It is compact when $\Delta(E):=4c_2(E)-c_1(E)^2\leq 3$.
\end{enumerate}
The second statement follows easily using the general properties of the Uhlenbeck compactification  of an instanton moduli space, and the well-known Chern class inequality for bundles admitting projectively ASD connections \cite{DK}, \cite{Te2}: when $\Delta(E)\leq 3$ all lower strata in the Uhlenbeck compactification are automatically empty, so ${\cal M}^{\rm ASD}_a(E)$ will be compact.

For the class of moduli spaces on which we will focus in this article, one has the following general result\footnote{In the first version of this article we stated this compactness result only for $b_2\leq 3$; this case is sufficient for the purposes of this article and has already been explained in \cite{Te2}.  The fact that this result can be extended to the case $b_2(X)>3$ has been noticed by Nicholas Buchdahl.}:
\begin{thry}\label{cp} Let $X$ be a class VII surface, ${\cal K}$ its canonical line bundle,  and $K$ the underlying ${\cal C}^\infty$-line bundle of ${\cal K}$. Let $(E,h)$ be a Hermitian rank 2-bundle   on $X$ with $c_2(E)=0$,  $\det(E)=K$. Then ${\cal M}^\pst_{\cal K}(E)$ is compact.
\end{thry}
\pf We have to prove that ${\cal M}^{\rm ASD}_a(E)$ is compact, where $a$ is the Chern connection of the pair $({\cal K},\det(h))$. One has $\Delta(E)=c_1(K)^2=-b_2(X)$, so the statement is obvious for $b_2(X)\leq 3$. For the general case, note that a stratum in the Uhlenbeck compactification of ${\cal M}^{\rm ASD}_a(E)$ has the form $S^k(X)\times {\cal M}^{\rm ASD}_a(E_k)$, where $S^k(X)$ stands for the $k$-th symmetric power of $X$, $\det(E_k)=\det(E)=K$, and $c_2(E_k)=c_2(E)-k=-k$. We claim that ${\cal M}^{\rm ASD}_a(E_k)=\emptyset$ for any  $k>0$. Using again the Kobayashi-Hitchin correspondence, it suffices to prove that there does not exist any holomorphic 2-bundle ${\cal E}$ on $X$ with   $\det({\cal E})={\cal K}$, and $c_2({\cal E})<0$. Identifying ${\cal E}$ with the corresponding locally free coherent sheaf, and using the Riemann-Roch theorem, one obtains easily $\chi({\cal E})>0$, so $h^0({\cal E})>0$ or   $h^0({\cal K}\otimes {\cal E}^\vee)>0$. In both cases  ${\cal E}$ would be filtrable, so it would fit  into an exact sequence of the form 
$$0\to {\mathcal L}\to {\mathcal E}\to {\mathcal K}\otimes{\mathcal L}^{\vee}\otimes {\cal I}_Z\to 0\ ,
$$
for a holomorphic line bundle ${\mathcal L}$ on $X$ and a  codimension 2  locally complete intersection $Z\subset X$. We use now the same method as in the proof of  Proposition 3.2 \cite{Te2}: Writing $c_1({\mathcal L})=\sum l_i e_i$ mod Tors (where $(e_i)_i$ is the basis considered in section \ref{firstresults} and $l_i\in\Z$), one obtains 
$0>c_2({\mathcal E})=|Z|+\sum_{i=1}^{b_2(X)} l_i(l_i-1)$,
which is obviously a contradiction.
\qed

 \subsubsection{The structure of ${\cal M}^{\rm ASD}_a(E)$  around the reductions}\label{localstructure}
 Let $(X,g)$ be a closed, connected, oriented  Riemannian 4-manifold with $b_+(X)=0$,   $(E,h)$   a Hermitian rank 2-bundle  on $X$, $L\hookrightarrow E$ a line subbundle of $E$ (endowed with the induced metric).  A connection $A\in{\cal A}_a(E)$ will be called $L$-{\it reducible} if it admits  a parallel line subbundle $L'\subset E$ isomorphic to $L$. The goal of this section is to describe the moduli space ${\cal M}^{\rm ASD}_a(E)$ around the subspace ${\cal M}^{\rm ASD}_a(E)^L$ of $L$-reducible instantons.  In \cite{Te4} we have shown that\\ 

{\it A neighborhood of ${\cal M}^{\rm ASD}_a(E)$ around the reduction locus ${\cal M}^{\rm ASD}_a(E)^L$   can be identified with the moduli space associated with an abelian moduli problem of Seiberg-Witten type}. \\

This identification has two important consequences:
\begin{itemize}
\item It describes explicitly a fundamental system of open neighborhoods of  the reduction locus ${\cal M}^{\rm ASD}_a(E)^L$ in ${\cal M}^{\rm ASD}_a(E)$, called standard neighborhoods. Under suitable regularity conditions any such neighborhood  can be identified with the total space of a fiber bundle whose  basis is a $b_1$-dimensional torus,  and whose fiber is a cone over a complex projective space. This generalizes the well-known theorem concerning the structure of a Donaldson moduli space around an {\it isolated reduction} \cite{FU}, \cite{DK}.
\item  It allows to compute explicitly the restriction of the Donaldson $\mu$-classes to the boundary of a standard neighborhood of  ${\cal M}^{\rm ASD}_a(E)^L$ in ${\cal M}^{\rm ASD}_a(E)$.
\end{itemize} 
\vspace{2mm}
We explain this formalism briefly.
Put  $M:=L^\bot$ and  $S:=L\otimes M^\vee\simeq L^{\otimes 2}\otimes D^\vee$, where $D:=\det(E)$.  The gauge group ${\cal G}_L:=\Aut(L)\simeq{\cal C}^\infty(X,S^1)$ acts on the space of Hermitian connections ${\cal A}(L)$ in the usual way 
$$u(d_b)=u\circ d_b\circ u^{-1}=d_b-u^{-1} du,$$
and on the space $A^1(S)$ of $S$-valued 1-forms via  the map 
$$\Aut(L)={\cal C}^\infty(X,S^1)\ni u\mapsto u^2\in {\cal C}^\infty(X,S^1)=\Aut(S)$$
induced by the identification $S=L^{\otimes 2}\otimes D^\vee$. 

The affine map $\Psi:{\cal A}(L)\times A^1(S)\to{\cal A}_a(E)$ defined by
 $$\Phi(b,\alpha)=A_{b,\alpha}:=\left(
 \begin{array}{cc}
 b&\alpha\\
 -\alpha^*& a\otimes b^\vee
 \end{array}\right) 
 $$
is equivariant with respect to the group morphism 
$$\psi:{\cal G}_L:=\Aut(L)\simeq{\cal C}^\infty(X,S^1)\to{\cal G}\ ,\ u\mapsto\left(\begin{matrix}u&0\cr 0&u^{-1}\end{matrix}\right)\ .$$
The connection  $A_{b,\alpha}$ is projectively ASD if and only if 
$$F_b^+-\frac{1}{2} F_a^+=(\alpha\wedge\alpha^*)^+\ ,\ d_{b^2\otimes a^\vee}^+\alpha=0\ .$$
Consider the moduli space ${\cal M}_a(L)$ of pairs $(b,\alpha)\in  {\cal A}(L)\times A^1(S)$ solving the system
$$
  \left\{
\begin{array}{ccc}
(d_{b^2\otimes a^\vee}^*,d_{b^2\otimes a^\vee}^+)\alpha&=&0\\\\
F_b^+ - \frac{1}{2} F_a^+&=&(\alpha\wedge\alpha^*)^+\ ,
\end{array}\right.   \eqno{(\Ag)}
$$
modulo the gauge group ${\cal G}_L$.   This abelian moduli problem   has been introduced in  \cite{Te2} and studied extensively in \cite{Te4}.  The factor group  $\bar {\cal G}_L:={\cal G}_L/\{\pm 1\}$ acts freely on the subspace ${\cal A}(L)\times (A^1(S)\setminus\{0\})$ of irreducible pairs, and acts with stabilizer $\bar S^1:=S^1/\{\pm 1\}\simeq S^1$ on the complement. The subspace of  ${\cal M}_a(L)$ consisting of solutions with $\alpha=0$ is a $b_1(X)$-dimensional torus. Indeed, this subspace can be identified with the  quotient 
$$T_a(L):=\qmod{{\cal T}_a(L)}{{\cal G}_L}$$
 of  the space ${\cal T}_a(L)$ of connections $b\in{\cal A}(L)$ solving the equation 
$$F_b^+ - \frac{1}{2} F_a^+=0\ .$$
 As explained in \cite{Te2}, \cite{Te4} this equation is equivalent to the condition 
$$F_b=\frac{1}{2} (-2\pi i\sigma+F_a)\ ,$$
 where $\sigma$ denotes the $g$-harmonic representative of the Chern class $c_1(S)$. Since we  assumed $b_+(X)=0$, this form is ASD. It is well-known that the moduli space of Hermitian connections with fixed curvature on a Hermitian line bundle is a $iH^1(X,\R)/2\pi i H^1(X,\Z)$-torsor, so it can be identified with this $b_1$-dimensional torus. The identification is well-defined, up to translations. We will identify  $T_a(L)$ with its image  in ${\cal M}_a(L)$ via the map $[b]\mapsto [b,0]$.

The following important lemma states that  a neighborhood of the reduction locus ${\cal M}^{\rm ASD}_a(E)^L$ in ${\cal M}_a^\ASD(E)$ can be identified with a neighborhood of $T_a(L)$ in ${\cal M}_a(L)$, so studying the Donaldson space  ${\cal M}_a^\ASD(E)$ around the reduction locus ${\cal M}^{\rm ASD}_a(E)^L$ reduces to studying ${\cal M}_a(L)$ around the torus $T_a(L)$.
\begin{lm}\label{isoms} Suppose $c_1(M)\ne c_1(L)$.  The map $\bar\Psi:{\cal M}_a(L)\to {\cal M}^{\rm ASD}_a(E)$ induced by $\Psi$ maps homeomorphically
\begin{enumerate}
\item the torus $T_a(L)$  onto ${\cal M}^{\rm ASD}_a(E)^L$,
\item  a sufficiently small open neighborhood ${\cal V}_L$  of $T_a(L)$ in ${\cal M}_a(L)$ onto an open neighborhood ${\cal U}_L$ of    ${\cal M}^{\rm ASD}_a(E)^L$ in ${\cal M}_a(L)$ such that
\begin{equation}\label{irrim}\bar\Psi({\cal V}_L\setminus T_a(L))\subset {\cal M}_a^\ASD(E)^*\ .
\end{equation}
\end{enumerate}
\end{lm}

\pf  1. {\it Surjectivity:} By definition, a connection $A\in {\cal M}^{\rm ASD}_a(E)^L$   splits as a direct sum $A=b'\oplus (a\otimes {b'}^\vee)$ with respect to an $A$-parallel orthogonal decomposition $E=L'\oplus  M'$, where $L'\simeq L$ (which implies $M'\simeq M$). Since $c_1(L)\ne c_1(M)$, this decomposition is unique (including the order).  Consider two Hermitian isomorphisms  $\lambda:L'\to L$, $\mu:M'\to M$ such that the induced Hermitian automorphism $v:E\to E$ belongs to ${\cal G}$ (i.e. it has determinant $\equiv 1$). Then $v(A)=\Psi(b,0)$, where $b:=\lambda (b')$, which shows that $[A]\in \bar\Psi(T_a(L))$.  

{\it Injectivity:} If $[\Psi(b)]=[\Psi(b_1)]$ there exists   $u\in {\cal G}$ such that $u\circ d_{\Psi(b)}=d_{\Psi(b_1)}\circ u$. Writing $u$ as a matrix with respect to the splitting $E=L\oplus M$, we see that $u^1_2\in A^0(\Hom(M,L))$ is $b_1\otimes (b\otimes a^{\vee})$-parallel, so $u^1_2=0$ because $L$ and $M$ are not isomorphic. Similarly one obtains $u^2_1=0$, hence $u^1_1\in{\cal G}_L$ and $u^1_1(d_b)=d_{b_1}$.\vspace{3mm}\\
2.\\ {\it Step 1.}  $\bar\Psi$ is a local homeomorphism at any point $[(b,0)]\in  T_a(L)$.

We use the standard procedure used in gauge theory to construct  locals models   (see \cite{DK} section 4.2.2).      Put $A_b=\Psi(b,0)=b\oplus(a^\vee\otimes b)$. The stabilizer $\bar{\cal G}_{A_b}$ of $A_b$ is the circle $\psi(\bar S^1)$. Consider the affine subspaces $S_b\subset {\cal A}(L)\times A^1(S)$, $\Sigma_{A_b}\subset {\cal A}_a(E)$ 
$$S_b:=(b,0)+\{(\beta,\alpha)\in  iA^1(X)\times A^1(S)|\ d^*\beta=0\} \ ,\  \Sigma_{A_b}:=A_b+i\ker d_{A_b}^*=$$
$$=A_b+\left\{\left(\begin{matrix} \beta& \alpha\cr-\alpha^*&-\beta\end{matrix}\right)\vline\ (\beta,\alpha)\in iA^1(X)\times A^1(S),\ d^* \beta=0,\ d_{b^2\otimes a^\vee}^*\alpha=0\right\}\ .
$$
$S_b$ (respectively $\Sigma_{A_b}$) is  $L^2$-orthogonal in $(b,0)$ (respectively in $A_b$) to the orbit  ${\cal G}_L\cdot (b,0)$ (respectively to the orbit ${\cal G}\cdot A_b$).   We denote by  $S^\Ag_b$, $\Sigma^{\rm ASD}_{A_b}$ the subspaces of   $S_b$, $\Sigma_{A_b}$ of points solving the equations $(\Ag)$ and $({\rm ASD})$ respectively, and  by $S^\Ag_{b,\eta}$, $\Sigma^{\rm ASD}_{A_b,\eta}$ the open subspaces defined by the inequality $\|(\alpha,\beta)\|_{L^2_k}<\eta$ (for a fixed sufficiently large Sobolev index $k$). Using standard gauge theoretical arguments (\cite{DK} section 4.2.2, Proposition 4.2.9) we see that for sufficiently small $\eta=\eta(b)>0$, the natural maps 
 $$\qmod{S^\Ag_{b,\eta}}{\bar S^1}\to {\cal M}_a(L)\ ,\ \qmod{\Sigma^{\rm ASD}_{A_b,\eta}}{\psi(\bar S^1)}\to {\cal M}_a^\ASD(E)$$
  are homeomorphisms onto open neighborhoods $V_b$, $U_b$ of  the orbits $[b,0]$ and $[A_b]$ in ${\cal M}_a(L)$ and ${\cal M}_a^\ASD(E)$ respectively.  But $\Psi$ maps $S^\Ag_{b,\eta}$ isomorphically onto $\Sigma^{\rm ASD}_{A_b,\eta}$, so $\bar\Psi$ defines a homeomorphism $V_b\to U_b$. 
\vspace{2mm}\\ 
{\it Step 2.}  $T_a(L)$ has a fundamental system of compact neighborhoods in ${\cal M}_a(L)$. 

Indeed, by the well-known bootstrapping procedure used to prove compactness in Seiberg-Witten theory \cite{KM}, it follows that
$$W^\varepsilon:=\{[b,\alpha]\in {\cal M}_a(L)|\ \|\alpha\|_\infty\leq \varepsilon\}
$$
is compact for any $\varepsilon>0$.  Using the compactness of the $W^\varepsilon$s, and the fact that $\bar\Psi$ is injective on $T_a(L)$ and local homeomorphic around  $T_a(L)$,  it follows easily (by reductio ad absurdum) that for sufficiently small $\varepsilon_0>0$, $\bar\Psi$ remains injective on $W^{\varepsilon_0}$. It suffices now to take
${\cal V}_L\subset W^{\varepsilon_0}\cap(\union_{b\in {\cal T}_a(L)} V_b)$.

  The equality (\ref{irrim}) follows directly from the second part of Proposition 4.2.9 \cite{DK}.
\qed
\begin{re} The proof of Lemma \ref{isoms} shows that the Kuranishi local model of a point $[b]\in T_a(L)$  in the moduli space ${\cal M}_a(L)$ is identified via $\bar\Psi$ with the   Kuranishi local model of $[A_b]$ in the moduli space  ${\cal M}^\ASD_a(E)$. Therefore, $\bar\Psi$ defines  an analytic isomorphism of real analytic spaces around  $T_a(L)$.
\end{re}
 
 Lemma  \ref{isoms} can be interpreted in the following way: {\it in a neighborhood of the reduction locus ${\cal M}^\ASD_a(E)^L\subset {\cal M}^\ASD_a(E)$, the ASD-moduli problem is equivalent to the abelian moduli problem (\Ag).} Therefore, it suffices to describe a neighborhood of the torus $T_a(L)$ in  ${\cal M}_a(L)$.  For this purpose it is convenient to introduce polar coordinates in the $\alpha$-direction. In other words, consider  the surjection 
 $$q: {\cal A}(L)\times {\mathbb S}(A^1(S))\times\R\to {\cal A}(L)\times A^1(S)\ ,\ q(b,\ag,\rho):=(b,\rho\ag)\ ,$$
where ${\mathbb S}(A^1(S))\subset A^1(S)$ denotes the $L^p$-unit  sphere of $A^1(S)$, and $p\in 2\N_{>0}$. Note that the  $L^p$-unit  sphere ${\mathbb S}(A^1(S)_k)$ in the $L^2_k$-Sobolev completion $A^1(S)_{k}$ is a real analytic hypersurface of this Hilbert space (if $k$ is   sufficiently large), because it is the fiber over the regular value 1 of  the real analytic  map $\alpha\mapsto \int_X(\alpha,\alpha)^{p/2}\vol_g$.

The pull-back PDE system $q^*(\Ag)$ on the new configuration space is  
 $$
  \left\{
\begin{array}{ccc}
(d_{b^2\otimes a^\vee}^*,d_{b^2\otimes a^\vee}^+)\ag&=&0\\ \\
F_b^+ - \frac{1}{2} F_a^+&=&\rho^2(\ag\wedge\ag^*)^+\ ,
\end{array}\right.   \eqno{(\tilde\Ag)}
$$
for triples $(b,\ag,\rho)\in {\cal A}(L)\times {\mathbb S}(A^1(S))\times\R$. 

In order to understand the idea behind this construction in an elementary example, consider the standard action of $S^1$ on $\C^n$.  A simple way to  describe the structure of the quotient $\C^n/S^1$ around  the singular orbit $\theta=\{0\}$ is to introduce polar coordinates writing $z=\rho \zg$, with $\zg\in S^{2n+1}$, $\rho\in \R$. $S^1$ acts freely on the set of pairs $(\rho,\zg)$, so the difficulty caused by the appearance of non-trivial stabilizers disappears. In this way we see immediately that $\C^n/S^1$ can be obtained from the cylinder  $[0,\infty)\times \P^{n-1}_\C$ (which is a smooth manifold with boundary) by collapsing the boundary $\{0\}\times \P^{n-1}_\C$ to a point.\\

The gauge group $\bar{\cal G}_L$ acts {\it freely} on the  whole space of triples $(b,\ag,\rho)$. Put
$$\tilde{\cal B}(L):=\qmod{{\cal A}(L)\times {\mathbb S}(A^1(S))\times\R}{\bar{\cal G}_L}\ ,$$
denote by  $\tilde{\cal M}_a(L)\subset \tilde{\cal B}(L)$ the moduli space of solutions of ($\tilde\Ag$), and consider its subspaces 
$$\hat{\cal M}_a(L):=\{[b,\ag,\rho]\in \hat{\cal M}_a(L)|\ \rho\geq 0\}\ ,\ \hat T_a(L):=\{[b,\ag,\rho]\in \hat{\cal M}_a(L)|\ \rho=0\}\ , $$

For $b\in{\cal T}_a(L)$ we denote   by ${\cal C}^+_{b^2\otimes a^\vee}$ the elliptic complex associated with the ASD connection $b^2\otimes a^\vee$ on $S$. With this notations, we can state the following lemma, which describes the geometry of ${\cal M}_a(L)$ around the reduction torus  $T_a(L)$.
 \begin{lm} \label{blowup} With the assumptions and notations of Lemma \ref{isoms} suppose also that, for every $b\in  {\cal T}_a(L)$, one has $H^2({\cal C}^+_{b^2\otimes a^\vee})=0$.  Then
 \begin{enumerate}
 \item   The assignment ${\cal T}_a(L)\ni b\mapsto \H^1_b:= \H^1({\cal C}^+_{b^2\otimes a^\vee})$ descends to a complex vector bundle $H_L$ of rank $-c_1(S)^2+(b_1-1)$  over the torus $T_a(L)$.
 \item $\hat T_a(L)$ can be identified  with the projectivization $\P(H_L)$.
  \item After suitable Sobolev completions the map $\tilde\Phi$ defined by the left hand side  of $(\tilde\Ag)$ is a submersion at any  solution $(b_0,\ag_0,0)$, so $\tilde{\cal M}_a(L)$ is a smooth manifold at any point  $[(b_0,\ag_0,0)]\in \hat T_a(L)$.
 \item The map  $r:\tilde{\cal M}_a(L)\to \R$ defined by the projection on the $\rho$-component is a submersion at any point $[(b,\ag,0)]\in \hat T_a(L)$, so $\hat {\cal M}_a(L)$ has the structure of a smooth manifold with boundary $\partial\hat {\cal M}_a(L)=\hat T_a(L)$ around this subspace.
 \item ${\cal M}_a(L)$ is obtained from $\hat {\cal M}_a(L)$ by collapsing to a point every fiber of the projective bundle $\partial\hat {\cal M}_a(L)=\P(H_L)\to T_a(L)$. 
 \item  $T_a(L)$ admits a fundamental system of neighborhoods  which can be identified with   cone bundles over  $\P(H_L)$.
 \item When $-c_1(S)^2+(b_1-1)=2$ (respectively $-c_1(S)^2+(b_1-1)=1$) the moduli space ${\cal M}_a(L)$ has the structure of a topological manifold  (respectively with boundary) around the torus of reductions $T_a(L)$.
 \end{enumerate}
 \end{lm}
 \pf 1. Since $c_1(L)\ne c_1(M)$, the line bundle $S$ is not  topologically trivial, so $H^0({\cal C}^+_{b^2\otimes a^\vee})=\ker d_{b^2\otimes a^\vee} =\{0\}$ for any connection  $b$.  Taking into account the hypothesis, we  see that  $(\H^1_b)_{b\in{\cal T}_a(L)}$  is the family of kernels associated with a smooth family of surjective elliptic operators
$$\delta_{b}:=(d_{b^2\otimes a^\vee}^*,d_{b^2\otimes a^\vee}^+):A^1(S)\map A^0(S)\oplus A^2_+(S)\ ,
$$
so the union ${\cal H}_L:=\union_{b\in {\cal T}_a(L)}\{b\}\times \H^1_b$ is the total space of a complex vector bundle over ${\cal T}_a(L)$.  Using the identity $\delta_{u(b)}=u^2\circ \delta_b\circ u^{-2}$, we see that ${\cal H}_L$ is gauge invariant. We   define $H_L$ to be the quotient of ${\cal H}_L$  by the based gauge group 
 $${\cal G}_{L,x_0}:=\{u\in{\cal G}_L\ |\ u(x_0)=1\}\subset {\cal G}_{L}$$
associated with a  point $x_0\in X$. The space ${\cal T}_a(L)$ is a ${\cal G}_{L,x_0}$-principal bundle over $T_a(L)$, so this quotient is a vector bundle over the basis $T_a(L)$ (see \cite{DK} p. 195). 

2. The cokernel of ${\cal G}_{L,x_0}\hookrightarrow \bar{\cal G}_L$ is $\bar S^1\simeq S^1$, which acts trivially on $T_a(L)$ and acts in the standard way on the fibers of $H_L$. 

3. This follows using standard transversality arguments, which we explain briefly: Let $\tau=(b_0,\ag_0,0)$ be a solution of $(\tilde \Ag)$ and  $(u,v)\in [A^0(S)\times A^2_+(S)]\times iA^2_+(X)$   a pair which is $L^2$-orthogonal to   $\im(d_\tau\tilde\Phi)$. Using variations $\dot\ag\in T_{\ag_0}:=T_{\ag_0}({\mathbb S}(A^1(S)))$ of $\ag_0$ and the surjectivity of the operator $\delta_{b_0}$ we obtain $u=0$. The important point here is that $\delta_{b_0}(T_{\ag_0})= \delta_{b_0}(A^1(S))$ because $\delta_{b_0}$ vanishes on the complement $\R\ag_0$ of $T_{\ag_0}$ in $A^1(S)$. Using variations $\dot b$ of $b_0$ and the surjectivity of $d^+$ we obtain $v=0$. 

4. The tangent space $T_\tau(\tilde\Phi^{-1}(0))=\ker(d_\tau\tilde\Phi)$   contains the line $\{0\}\times\{0\}\times \R$, which is obviously mapped surjectively onto $\R$ via the projection on the third factor.  

5. is obvious. To prove 6. note that, by 2., 3.  and 4., $\hat T_a(L)$ has a fundamental system of neighborhoods in  $\hat {\cal M}_a(L)$ which can be identified $\P(H_L)\times [0,\epsilon]$, such that the projection on the second factor corresponds to $r$. It suffices to apply 5.

  7.  The condition $-c_1(S)^2+(b_1-1)=2$ implies $\rk(H_L)=2$. In this case the cone over $\P^1_\C\simeq S^2$ is $D^3$, so the cone bundle of $\P(H_L)$ is a $(b_1+3)$-dimensional topological manifold around its tautological section. When $-c_1(S)^2+(b_1-1)=1$, the cone bundle over $ \P(H_L)$ can be identified with $T_a(L)\times [0,\varepsilon]$. \qed

Combining Lemmata \ref{blowup}, \ref{isoms} we obtain the following result concerning the structure of the Donaldson moduli space around the reduction locus ${\cal M}_a^\ASD(E)^L$. Geometrically, the meaning of this result is very simple:  replacing in ${\cal M}_a^\ASD(E)$ the reduction  locus ${\cal M}_a^\ASD(E)^L$ (which is a torus) with the projective bundle $\P(H_L)$, we get a space which has a natural structure of manifold with boundary around $\P(H_L)$.
\begin{co} (Blowing up reduction loci) \label{localred} Let $L\hookrightarrow E$ be a line subbundle of $E$ such that $2c_1(L)\ne c_1(E)$  and $H^2({\cal C}^+_{\su(E),A})=0$ for any $L$-reducible instanton $A$. 
\begin{enumerate}
\item The reduction locus ${\cal M}_a^\ASD(E)^L$ can be identified with the torus $T_a(L)$.
\item The union $\hat {\cal M}_a^\ASD(E)_L:=[{\cal M}_a^\ASD(E) \setminus {\cal M}_a^\ASD(E)^L]\cup \P(H_L)$ 
has a natural structure of   manifold with boundary $\partial\hat {\cal M}_a^\ASD(E)_L=\P(H_L)$ around $\P(H_L)$, and comes with a  continuous map   $\hat {\cal M}_a^\ASD(E)_L\to {\cal M}_a^\ASD(E)$ which is the identity on the complement of $\P(H_L)$, and whose restriction to this subspace is the bundle projection $\P(H_L)\to T_a(L)$. 
\item The reduction locus ${\cal M}_a^\ASD(E)^L$ has a fundamental system of neighborhoods ${\cal U}_L^\varepsilon$ which can be identified with   cone bundles over $\P(H_L)$ and have the property ${\cal U}_L^\varepsilon\setminus {\cal M}_a^\ASD(E)^L\subset{\cal M}_a^\ASD(E)^*$. 
\item The union 
$$\hat {\cal M}_a^\ASD(E)_L^*:={\cal M}_a^\ASD(E)^* \cup \P(H_L)$$  has  the structure of manifold with boundary $\P(H_L)$ around $\P(H_L)$.  
\item When $-c_1(S)^2+(b_1-1)=2$ (respectively $-c_1(S)^2+(b_1-1)=1$) the moduli space ${\cal M}_a^\ASD(E)$ is a topological manifold (respectively with boundary) around the reduction locus ${\cal M}_a^\ASD(E)^L$.
\end{enumerate}
 \end{co}
 
   The Chern class $c_1(D)$ is an integral lift of the Stiefel-Whitney  class $w:=w_2(\su(E))$. Choosing a subbundle $L\hookrightarrow E$ defines a new integral lift of $w$, namely $\hat w_L=c_1(L^{\otimes 2}\otimes D^\vee)$. Together with a fixed orientation  $\oo_1$ of $H^1(X,\R)$ this lift defines an orientation $\oo_L$ of the regular part of the moduli space ${\cal M}^*_a(E)$ (see \cite{DK} section 7.1 p. 283). The restriction of $\oo_L$  to $[{\cal U}^\varepsilon_L]^*:={\cal U}_L^\varepsilon\setminus {\cal M}_a^\ASD(E)^L$ is just the obvious orientation induced by the the natural complex orientations of the fibers $\P(H_{L,[b]})$ and the orientation of the base $T_a(L)$ induced by $\oo_1$. The rule explained in   \cite{DK} loc. cit. states that
\begin{lm}\label{orientation}  The orientations $\oo_{L_1}$,  $\oo_{L_2}$ associated with two line subbundles $L_i\hookrightarrow E$ compare according to the  parity of $(c_1(L_2)-c_1(L_1))^2$.
\end{lm} 

\subsubsection{Donaldson classes around the reductions} \label{donclass}

We recall (see \cite{DK}) that the quotient $\V:=[{\cal A}^*_a(E)\times \su(E))]/\bar {\cal G}$ can be regarded as an $\SO(3)$ vector bundle over ${\cal B}^*_a(E)\times X$; this is the Donaldson  universal vector bundle. The Donaldson map $\mu:H_i(X,\Q)\to H^{4-i}({\cal B}^*_a(E),\Q)$ (see \cite{DK} Definition 5.1.11) is defined by
$$\mu(\theta):=-\frac{1}{4} p_1(\V)/\theta\ .
$$
In general the $\mu$-classes do not extend across the reduction loci ${\cal M}_a^\ASD(E)^L$.  However  these classes do extend in a canonical way to the blown up moduli spaces $\hat {\cal M}_a^\ASD(E)_L^*$  and, in the conditions of Corollary \ref{localred}, the obtained classes in the cohomology of the boundaries $\P(H_L)$ can be computed explicitly. The point is that   $\V$ extends   in a natural way across the boundary $\P(H_L)\times X$ of  $\hat {\cal M}_a^\ASD(E)_L^* \times X$, and this extension has an obvious  $S^1$-reduction  (i.e. an obvious splitting as the sum of a complex line bundle and a trivial $\R$-bundle). To see this it suffices to note that the pull-back of $\V$ via the natural map  $[\tilde{\cal B}(L)_{\rho>0}]\times X\to {\cal B}^*(E)\times X$ splits as   $i\R\times [\resto{\Sg}{[\hat{\cal B}(L)_{\rho>0}]\times X}]$, where $\Sg$ is the universal complex line bundle
$$\Sg:=\qmod{{\cal A}(L)\times {\mathbb S}(A^1(S))\times \R \times S}{\bar {\cal G}_L}\ ,
$$
which is defined on the whole of $\tilde {\cal B}(L)\times X$. The restriction $\resto{\Sg}{\P(H_L)\times X}$ is the line bundle $\left\{ {\mathbb S}({\cal H}_L)\times S\right\}/\bar {\cal G}_L$ over $\P(H_L)\times X$.   The Chern class of this restriction  can be computed easily  as in Seiberg-Witten theory  (see \cite{Te4} Corollary 2.5). Using these facts we obtain (see \cite{Te4} Corollary 2.6):
\begin{lm}\label{muclass} Let $\gamma\in H_1(X,\Q)$ and $\delta(\gamma)\in H^1(T_a(L),\Q)$ the associated class via the obvious isomorphism. The Donaldson class $\mu(\gamma)$ extends to   $\hat {\cal M}_a^\ASD(E)_L^*$ and the restriction of this extension  to  the boundary $\P(H_L)$ is $-p^*(\delta(\gamma))\cup h$ where $p:\P(H_L)\to T_a(L)$ is the bundle projection and $h\in H^2( \P(H_L),\Q)$ is the Chern class of the $S^1$-bundle ${\mathbb S}(H_L)\to  \P(H_L)$. When $b_1=1$, $c_1(S)^2=-2$ and $\gamma$ is a generator of $H_1(X,\Z)/\Tors$,  this class is (up to a universal sign) the fundamental class of the boundary $\partial \hat {\cal M}_a^\ASD(E)_L$, with respect to the orientation  $\oo_L$ defined by $L$  and the orientation $\oo_1$ of $H^1(X,\R)$ associated with $\gamma$.
\end{lm}

In order to explain in a geometric way the role played by these gauge theoretical results in this article, consider again a minimal class  VII  surface with $b_2(X)=2$, $\pi_1(X)\simeq\Z$, and the corresponding  moduli space ${\cal M}^\pst(0,{\cal K})$  described in the introduction. We know now that (under suitable regularity conditions) this space is a compact topological 4-manifold containing two circles of reductions $\Rg'=T_a(L')$, $\Rg''=T_a(L'')$. In the introduction we claimed that {\it the circles $\Rg'$, $\Rg''$ must belong to the same connected component of the moduli space.} Indeed, if not, one would get two compact oriented  manifolds with boundary $\hat {\cal M}'$, $\hat {\cal M}''$ having the closed 3-manifolds $\P(H_{L'})$,  $\P(H_{L''})\simeq S^1\times S^2$  as boundaries. But this is impossible  because, for a generator $\gamma$ of $H_1(X,\Z)/\Tors$, the  Donaldson class $\mu(\gamma)$  extends to $\hat {\cal M}'$, $\hat {\cal M}''$, and its restriction to the two boundaries is non-trivial. The same idea will be used in the proof of  Proposition \ref{firststr} concerning the structure of ${\cal M}^\pst(0,{\cal K})$ on a general minimal class VII surface with $b_2=2$.
 
 \subsubsection{Comparing deformation elliptic complexes}\label{compell}

Let $(X,g)$ be a Gauduchon surface, $F$ an Euclidean rank $r$ bundle and $B$ an ASD connection on $F$. Our goal now is to compare the cohomology of the elliptic complex 
$$0 \map A^0(F) \textmap{d_B} A^1(F) \textmap{d^+_B} A^2_+(F) \map  0 \eqno{({\cal C}^+_B)}$$
of $B$ with the Dolbeault elliptic complex associated with the  operator $\bar\partial_B$ on its complexification $F^\C$.  The complexified bundle $F^\C$ comes with an induced Hermitian structure and an induced Hermitian connection  (denoted  by the same symbol $B$) which will also be ASD, because the two curvature forms coincide via the obvious embedding $A^2(\so(F))\hookrightarrow A^2(\su(F^\C))$. The holomorphic bundle $(F^\C,\bar\partial_B)$ is a polystable bundle of degree 0. We have an obvious morphism of elliptic complexes
\begin{equation}
\label{morcmp}
\begin{array}{ccccccccc}
0&\map&A^0(F)&\textmap{d_B}&A^1(F)&\textmap{d^+_B}&A^2_+(F)&\map& 0\ \phantom{.}
\\  \\
&&\phantom{j_0}\downarrow j_0&&\simeq\downarrow j_1 &&\phantom{j_2=p^{02}}\downarrow j_2:=p^{02}
\\  \\
0&\map&A^{0}(F^\C)&\textmap{\bar\partial_B}& A^{01}(F^\C)
&\textmap{\bar\partial_B}&A^{02}(F^\C)&\map& 0\ .
\end{array}
\end{equation}
Denote by  $H^j(B)$, $H^j(\bar\partial_B)$ the corresponding cohomology spaces.
\begin{lm} Let $B$ be an ASD connection and suppose that $F$ has no  non-trivial $B$-parallel section. Then 
\begin{enumerate} 
\item The operators $\Lambda_g d^c_B d_B:A^0(F)\to A^0(F)$, $i\Lambda_g\bar\partial_B\partial_B:A^0(F^\C)\to A^0(F^\C)$ are isomorphisms.
\item The natural morphisms
$$\Hg^1_B:=\{a\in A^1(F)|\ d^+_B(a)=0\ ,\ \Lambda_g d_B^c a =0\}\to H^1(B)\ ,
$$
$$\Hg^1_{\bar\partial_B}:=\{\alpha\in A^{01}(F^\C)|\ \bar\partial_B(\alpha)=0\ ,\ \Lambda_g \partial_B \alpha =0\}\to H^1(\bar\partial_B)
$$
are isomorphisms.
\end{enumerate}
\end{lm}
\pf 1. Note first that, since $B$ is ASD,  one has the identities
$$\bar\partial_B \bar\partial_B=0\ ,\  \partial_B\partial_B=0\ ,\ \Lambda_g\bar\partial_B \partial_B=-\Lambda_g\partial_B \bar\partial_B\ ,\ \Lambda_g d_B^c d_B=-\Lambda_g d_B d_B^c\ .
$$
 Suppose that $\varphi\in \ker (i\Lambda_g\bar\partial_B \partial_B)$, with $\varphi\in A^0(F^\C)$. Using the same method  as in the proof of Proposition 4.3 \cite{Te2} one gets $i\Lambda_g\bar\partial\partial|\varphi|^2=-|d_B\varphi|^2$,
which implies that $|\varphi|^2$ is constant and $d_B\varphi=0$ by the maximum principle. Therefore $\varphi=0$, since $F$ (so also $F^\C$) has no  non-trivial $B$-parallel section by assumption.
The same argument proves the injectivity of $\Lambda d^c_B d_B$.  It suffices to note that the two operators have vanishing index (because they have self-adjoint symbols). \\
2. For the first isomorphism, one has to prove that for every  $F$-valued 1-form $a\in \ker(d^+_B)$ there exists a unique $\beta\in A^0(F)$ such that $i\Lambda_g d^c_B(a+d_B\beta)=0$. This follows from (1). The same argument applies for the second morphism.
\qed

The following comparison theorem is  known in the K\"ahlerian framework \cite{K}.
\begin{pr}\label{compcplx} Suppose that $F$ has no  non-trivial $B$-parallel section.  The diagram (\ref{morcmp}) induces isomorphisms $H^0(B)=H^0(\bar\partial_B)=0$,  $H^1(B)=H^1(\bar\partial_B)$, $H^2(B)=H^2(\bar\partial_B)$.
\end{pr}
 \pf $H^0(\bar\partial_B)=0$: Since $(F^\C,\bar\partial_B)$ is a polystable bundle of degree 0, any $\bar\partial_B$-holomorphic section is parallel \cite{K}. \vspace{2mm}\\
 $(j_1)_*$ is an isomorphism: it suffices to note that the map $\Hg^1_B\to \Hg^1_{\bar\partial_B}$ given by $a\mapsto a^{01}$ is an isomorphism. \vspace{2mm}\\
 $(j_2)_*$ is an isomorphism: the surjectivity is obvious. For the injectivity, consider an element $(a^{20}+a^{02}+ u\omega_g)\in\ker (d_B^+)^*$. This means $\partial_B a^{02}+ \bar \partial_B (u\omega_g)=0$, which implies $\partial_B\bar \partial_B (u\omega_g)=0$. Using   the properties of the  operator $i\Lambda_g\bar\partial_B  \partial_B$ and its adjoint (\cite{Te2}), it follows that $u$ is $B$-parallel so, under our assumption,  $u=0$. Here the Gauduchon condition $\bar\partial\partial\omega_g=0$ plays a crucial role. This shows
 $$\ker (d_B^+)^*=\{a^{20}+a^{02}|\ a^{0,2}\in A^{0,2}(F^\C), \ a^{2,0}=\bar a^{0,2},\ \partial_B a^{02}=0\}
 $$
 which is obviously identified with $\ker(\bar\partial_B^*:A^{02}(F^\C)\to A^{01}(F^\C))$ via $j_2$.
 \qed

 Let now $B$ be any ASD connection on $F$. Consider the $B$-parallel decomposition $F=[X\times H^0(B)]\oplus F^\bot$. Let $B^\bot$ be the connection induced on $F^\bot$. One has obvious  isomorphisms:
 $$  H^1(B)=H^1(B^\bot)\oplus [H^0(B)\otimes H^1(X,\R)],\ H^2(B)=H^2(B^\bot)\oplus  [H^0(B)\otimes H^2_+(X,\R)]$$
 $$ H^0(\bar\partial_B)=H^0(B)^\C,\  H^1(\bar\partial_B)=H^1(\bar\partial_{B^\bot})\oplus [H^0(\bar\partial_B)\otimes H^{01}(X)]\ ,$$
 $$ H^2(\bar\partial_B)=H^2(\bar\partial_{B^\bot})\oplus  [H^0(\bar\partial_B)\otimes H^{02}(X)] \ , \
 $$
Applying Proposition \ref{compcplx} to $B^\bot$, one gets the following important
 \begin{co}  \label{comcoh}
\begin{enumerate}
\item If $g$ is K\"ahlerian, $H^1(B)=H^1(\bar\partial_B)$, and   $H^2(\bar\partial_B)$ is a subspace of real codimension $h^0(B)$ in  $H^2(B)$.
\item If $b_1(X)$ is odd, $H^1(B)$ is a subspace of real codimension $h^0(B)$ in $H^1(\bar\partial_B)$, and $H^2(B)=H^2(\bar\partial_B)$.
\end{enumerate}
  \end{co}
 We apply these results to the Euclidean bundle $F:=\su(E)$ associated with a Hermitian bundle $E$ endowed with an {\it oriented} integrable projectively ASD connection $A\in{\mathcal A}_a(E)$.  The deformation complex of   $A$ is just the  complex ${\cal C}^+_{\su(E),A}$  of the pair $(\su(E),A)$, whereas the deformation elliptic complex of an {\it oriented} holomorphic structure ${\cal E}$ on $E$ is the 
 Dolbeault complex of the pair $(\End_0(E),\bar\partial_{\cal E})$. 
  \begin{co}\label{compconcl} Let $(E,h)$ be a Hermitian rank $r$-bundle on a Gauduchon surface $(X,g)$ with $b_1(X)$ odd, and ${\cal E}$ the  oriented polystable holomorphic structure on $E$ associated to $A$ via the Kobayashi-Hitchin correspondence. The second cohomology of the deformation elliptic complex of  $A$ can be identified with the second cohomology of the  deformation elliptic complex of   ${\cal E}$.
  \end{co}

 \subsubsection{The structure of ${\cal M}^\pst(0,{\cal K})$ around the reductions}\label{general}

Let $X$ be a class VII surface with $b_2(X)=2$ and $g$ a Gauduchon metric on $X$.  We denote 
$$\kg:=\frac{1}{2}\deg_g({\cal K})$$
 (the slope of the bundles in our moduli space). The subspace   ${\mathcal M}^{\red}(0,{\mathcal K})$  of reductions has two disjoint parts
$$\Rg':=\left\{{\mathcal L}\oplus({\mathcal K}\otimes{\mathcal L}^{\vee})\ \vline\ {\mathcal L}\in\Pic^T_{=\kg}\right\}\ ,\ \Rg'':=\left\{{\mathcal L}\oplus({\mathcal K}\otimes{\mathcal L}^{\vee})\                  
\vline\ {\mathcal L}\in\Pic^{e_1}_{=\kg}\right\}
$$
which are disjoint unions of $\tau:=|\Tors(H^2(X,\Z))|$ circles.  

It is certainly impossible to extend the complex space structure of ${\mathcal M}^{\st}(0,{\mathcal K})$ across $\Rg'$, because any neighborhood of a point ${\mathcal L}\oplus({\mathcal K}\otimes{\mathcal L}^{\vee})\in \Rg'$ contains holomorphic curves (namely projective lines of stable nontrivial extensions of ${\mathcal K}\otimes{\mathcal L'}^{\vee}$ by ${\cal L}'$, where ${\cal L}'$ is close to ${\cal L}$ and has $\deg_g({\cal L}')<\kg$).

The purpose of this section is to show that, if $X$ is minimal,  ${\mathcal M}^{\pst}(0,{\mathcal K})$ has a natural smooth holomorphic structure around $\Rg''$, extending the complex structure of   ${\mathcal M}^{\st}(0,{\mathcal K})$. The idea is very simple. Perturbing the metric in a convenient way, one gets a new moduli space ${\mathcal M}^{\pst}_{g_t}(0,{\mathcal K})$ which can be homeomorphically identified with the old ${\mathcal M}^{\pst}(0,{\mathcal K})$.  This homeomorphism maps $\Rg''$ into the smooth locus of ${\mathcal M}^{\st}_{g_t}(0,{\mathcal K})$ and its inverse maps $\Rg''_{g_t}$ into the smooth locus of ${\mathcal M}^{\st}(0,{\mathcal K})$. Our main tool will be the following easy
\begin{lm}\label{newlemma} Let $f:U_1\to U_2$  be a homeomorphism between Hausdorff paracompact topological spaces and $F_i\subset U_i$ closed sets such $f(F_1)\cap F_2=\emptyset$. Suppose that $U_i\setminus F_i$ are endowed with structures of complex manifolds such that 
$$\resto{f}{U_1\setminus(F_1\cup f^{-1}(F_2))}:U_1\setminus(F_1\cup f^{-1}(F_2))\map U_2\setminus(F_2\cup f (F_1))$$
is biholomorphic. Then there exist unique structures of complex manifolds   on $U_i$ for which $f$ is  biholomorphic.
\end{lm}
\pf We define a holomorphic atlas ${\cal A}_1$ on $U_1$ by ${\cal A}_1={\cal B}_1\cup{\cal C}_1$ where ${\cal B_1}$ is a holomorphic  atlas of $U_1\setminus F_1$ and ${\cal C}_1$ is the set of  maps of the form  $h\circ f:f^{-1}(U)\to V$, where $h:U\to V\subset\C^n$  is a holomorphic chart of $U_2\setminus F_2$. Such a map is holomorphically compatible with ${\cal B}_1$ and the domains of these maps cover $F_1$. Similarly for $U_2$. \qed

Suppose now that $X$ is a minimal class VII surface with $b_2=2$. By the vanishing  Lemma \ref{nak}, we know that $h^2({\mathcal K}^\vee \otimes{\mathcal L}^{\otimes 2})=0$, for every ${\mathcal L}\in\Pic^{e_i}$, $i=1,\ 2$. Note that $\deg_g({\mathcal K}^\vee \otimes{\mathcal L}^{\otimes 2})=2(\deg_g({\cal L})-\kg)$. Since $X$ cannot contain curves with arbitrary small volume, and ${\mathcal K}^\vee \otimes{\mathcal L}^{\otimes 2}$ is not topologically trivial, it follows  that there exist $\varepsilon>0$ such that $h^0({\mathcal K}^\vee \otimes{\mathcal L}^{\otimes 2})=0$ for every  ${\mathcal L}\in\Pic^{e_i}_{<\kg+\varepsilon}$.  Therefore, by the Riemann-Roch theorem, one gets $h^1({\mathcal K}^\vee \otimes{\mathcal L}^{\otimes 2})=1$ so, for every  ${\mathcal L}\in\Pic^{e_i}_{<\kg+\varepsilon}$,  there exists an essentially unique non-trivial extension  ${\mathcal E}({\mathcal L})$ of ${\mathcal K}\otimes{\mathcal L}^\vee$ by ${\mathcal L}$. 

For every $\eta\in  (\kg-\varepsilon,\kg+\varepsilon)$ we define 
$$\varphi_\eta: [\Pic^{e_1}]_{<\kg+\varepsilon}^{>\kg-\varepsilon}\map \{\hbox{Isomorphism classes of bundles on }X\}
$$
by
$${\mathcal L}\mapsto\left\{\begin{array}{ccc} {\mathcal E}({\mathcal L})&{\rm when} &\deg_g({\mathcal L})\in(\kg-\varepsilon,\eta)\\
{\mathcal L}\oplus({\mathcal K}\otimes{\mathcal L}^{\vee})&{\rm when} &\deg_g({\mathcal L})=\eta\\
{\mathcal E}({\mathcal K}\otimes{\mathcal L}^\vee)&{\rm when} &\deg_g({\mathcal L})\in(\eta,\kg+\varepsilon)\ .
   \end{array}\right.
 $$
 \begin{lm} If $\varepsilon$ is sufficiently small, then $\varphi_\eta$ is injective for every $\eta\in(-\varepsilon,\varepsilon)$.
 \end{lm}
 \pf The bundles ${\cal E}_i:=\varphi_\eta({\cal L}_i)$ associated to two line bundles ${\cal L}_i\in  [\Pic^{e_1}]_{<\kg+\varepsilon}^{>\kg-\varepsilon}$ are given by two extensions
 $$0\map {\cal M}_i^{(1)}\map {\cal E}_i\map {\cal M}_i^{(2)}\map 0\ .
 $$
The $g$-degree of   ${\cal M}^{(j)}_i$ belongs to $(\kg-\varepsilon,\kg+\varepsilon)$. Since $X$ cannot contain   curves of arbitrary small volume,  we see that, for sufficiently small $\varepsilon$,   $H^0(({\cal M}_1^{(j)})^\vee\otimes {\cal M}_2^{(k)})=0$ except when ${\cal M}_1^{(j)}\simeq {\cal M}_2^{(k)}$. The result follows now from Corollary \ref{injectivity}.
 \qed

Let $h_i$ be the harmonic representative of the de Rham class $e_i$. For any sufficiently small $t>0$, the form $\omega_g+t(h_1-h_2)$ is the K\"ahler form of a Gauduchon metric $g_t$ on $X$. Note that 
\begin{equation}\label{degchange}
\deg_{g_t}({\mathcal L})=\left\{
\begin{array}{cc}
\deg_{g}({\mathcal L}) &\forall {\mathcal L}\in\Pic^T\cup\Pic^{\bar k},\\ \deg_{g}({\mathcal L})+(-1)^i t & \forall {\mathcal L}\in\Pic^{e_i} ,\ i\in\{1,2\}\ .
\end{array}\right.
\end{equation}
When we pass from $g$ to $g_t$  the stability properties of all bundles are preserved, except certain type $\{1\}$   and type $\{2\}$ extensions. More precisely:
\begin{lm} \label{defmet} For sufficiently small $\varepsilon>0$, the following holds:
 For every $t \in (0,\varepsilon)$ one has 
\begin{equation}\label{compmoduli}
\begin{split}
\im(\varphi_{\kg+t})\subset {\mathcal M}^{\pst}_{g_t}(0,{\mathcal K})\ ,\ \varphi_{\kg+t}\left([\Pic^{e_1}]_{<\kg+\varepsilon}^{>\kg-\varepsilon}\setminus [\Pic^{e_1}]_{=\kg+t}\right)\subset {\mathcal M}^{\st}_{g_t}(0,{\mathcal K})\ ,
\\
{\mathcal M}^{\pst}_{g_t}(0,{\mathcal K})\setminus  \varphi_{\kg+t}\left([\Pic^{e_1}]_{\leq\kg+t}^{\geq\kg}\right)= {\mathcal M}^{\pst}_{g}(0,{\mathcal K})\setminus \varphi_{\kg}\left([\Pic^{e_1}]_{\leq\kg+t}^{\geq\kg}\right)  \ .
\end{split}
\end{equation}
\end{lm} 
 \pf When $\deg({\cal L})=\kg+t$, the bundle $\varphi_{\kg+t}({\cal L})$ is obviously a split $g_t$-polystable bundle. Suppose now $\deg({\cal L})\ne\kg+t$. The bundle ${\cal E}= \varphi_{\kg+t}({\cal L})$ is defined as the central term of a non-trivial extension $0\to{\cal M}\textmap{u} {\cal E}\textmap{v} {\cal N}\to 0$ with $\deg_g({\cal M}),\  \deg_g({\cal N})\in (\kg-\varepsilon,\kg+\varepsilon)$. By the definition  of $\varphi_{\kg+t}$, the monomorphism $u$ does not $g_t$-destabilize ${\cal E}$. Suppose there existed a $g_t$-destabilizing morphism $w:{\cal S}\to {\cal E}$ with torsion free quotient. The destabilizing condition reads $\deg_{g_t}({\cal S})\geq \frac{1}{2}\deg_{g_t}({\cal K})=\kg$. By Proposition \ref{extypes} we may assume $c_1({\cal S})\in e_I$ with $I\subset\{1,2\}$. Therefore $\deg_{g_t}({\cal S})\leq \deg({\cal S})+t$ by (\ref{degchange}). One has 
 $$\deg_{g}({\cal S}^\vee\otimes{\cal N})=-\deg_{g}({\cal S})+\deg_{g}({\cal N})\leq -\deg_{g_t}({\cal S})+t+\deg_{g}({\cal N})\leq-\kg+t+(\kg+\varepsilon)\ .$$
 
 Therefore $H^0({\cal S}^\vee\otimes{\cal N})=0$ if $\varepsilon$ is sufficiently small and $S\not\simeq {\cal N}$.  But for  ${\cal S}\simeq {\cal N}$ we have $v\circ w=0$, because otherwise the extension defined by $(u,v)$ would split. In all cases we get $v\circ w=0$, so $w$ would factorize through a morphism ${\cal S}\to{\cal M}$, so it cannot $g_t$-destabilize ${\cal E}$. This proves the two inclusions in (\ref{compmoduli}).
 
 For the third formula in (\ref{compmoduli}), note first that the subsets of split polystable bundles in the two sets coincide. Suppose now that ${\cal E}$   belongs to the left hand set, but is non-split. Therefore it is $g_t$-stable. If ${\cal E}$ was 
  not $g$-stable there would exist a monomorphism  $u:{\cal S}\to {\cal E}$ with $\deg_g({\cal S})\geq \kg$ and torsion free quotient. By  Proposition \ref{extypes}   $u$ is a bundle embedding and $c_1({\cal S})\in e_I$ where $I\subset\{1,2\}$.  Since  ${\cal E}$ is $g_t$-stable we get by (\ref{degchange}) that $c_1({\cal S})\in e_1\cup e_2$   and $\deg_{g}({\cal S})+(-1)^i t <\kg$ if ${\cal S}\in\Pic^{e_i}$.    This would imply ${\cal E}\in  \varphi_{\kg+t}\left([\Pic^{e_1}]_{\leq\kg+t}^{\geq\kg}\right)$, which contradicts the choice of ${\cal E}$. This proves that ${\cal E}$ is $g$-stable. ${\cal E}$ cannot belong to $\varphi_{\kg}\left([\Pic^{e_1}]_{\leq\kg+t}^{\geq\kg}\right)$, because the bundles   in  this set are not $g_t$-polystable. The other inclusion is proved similarly.\qed
\begin{pr}\label{extholstr} Let $X$ be a minimal class VII surface with $b_2=2$. Then    
\begin{enumerate}
\item  For a sufficiently small neighborhood ${\mathcal U}''$ of $\Rg''$, ${\mathcal U}''\setminus \Rg''$ is contained in ${\mathcal M}^{\st}(0,{\mathcal K})$ and is a smooth complex manifold.

\item  The holomorphic structure of ${\mathcal U}''\cap{\mathcal M}^{\st}(0,{\mathcal K})$ extends across $\Rg''$ such that $\im(\varphi_\kg)$ is a holomorphic curve.

\end{enumerate}  
\end{pr}
\pf The proof of Proposition \ref{reg} shows that, if ${\mathcal E}$ is any  extension  of type $\{1\}$ or $\{2\}$ and $X$ is minimal, then $H^2({\mathcal E}nd_0({\mathcal E}))=0$.  Therefore, for any ${\mathcal E}$ in a sufficiently small open neighborhood ${\mathcal U}''$ of $\Rg''$, one will still  have $H^2({\mathcal E}nd_0({\mathcal E}))=0$. We can choose this neighborhood such that ${\mathcal U}''\cap \Rg'=\emptyset$. This proves 1.

For 2. choose $t\in (0,\varepsilon)$, and consider the following symmetric  relation between the   moduli spaces of polystable bundles associated with $g$ and $g_t$:   
$$R_t=\left\{({\mathcal E},{\mathcal E}')\in {\mathcal M}^{\pst}_{g}(0,{\mathcal K})\times {\mathcal M}^{\pst}_{g_t}(0,{\mathcal K}) \vline\  h^0{\mathcal H}om({\mathcal E},{\mathcal E}')\ne 0,\ h^0 {\mathcal H}om({\mathcal E}',{\mathcal E}'') \ne 0\right\}.
$$
When ${\mathcal E}$ is polystable with respect to both metrics, it is in relation only with itself. Using Lemma 
\ref{defmet} it is easy to check that $R_t$ is in fact one-to-one. The corresponding  bijections are continuous by elliptic semicontinuity, so $R_t$ defines a homeomorphism $r_t$. Put ${\cal U}''_t:=r_t({\cal U}'')\subset{\mathcal M}^{\pst}_{g_t}(0,{\mathcal K})$.    ${\cal U}''_t\setminus  {\cal R}''_{g_t}$ is also a smooth complex manifold, because a point in this set is either an element of ${\cal U}''\setminus\Rg''$, or an extension of type $\{1\}$ or $\{2\}$. The claim follows now directly from  Lemma \ref{newlemma}.
 \qed

 Taking into account the proofs of Lemma \ref{newlemma} and Proposition \ref{extholstr}, it follows
 \begin{lm}\label{holmap} Let $Y$ be a complex manifold,  $\tau:Y\map {\cal M}^\pst(0,{\cal K})$ and $y_0\in Y$ a point such that $\tau(y_0)\in{\cal R}''$. The map $\tau$ is holomorphic at $y_0$ with respect to the holomorphic structure given by Proposition \ref{extholstr}  if and only if $r_t\circ \tau$ (which maps $y_0$ to a $g_t$-stable extension) is holomorphic at $y_0$.   
 \end{lm}
 Around the other part $\Rg'$ of the reduction locus, the holomorphic structure does not extend. However, using Corollary \ref{compconcl}  and Corollary \ref{localred}, we get
\begin{pr}\label{loc} Suppose that $H^2({\cal E}nd_0({\cal E}))=0$ for all ${\cal E}\in \Rg'$. Then ${\cal M}^\pst(0,{\cal K})$ has the structure of a topological manifold around $\Rg'$.
\end{pr}

\section{The moduli space ${\mathcal M}^{\rm pst}(0,{\mathcal K})$ in the case $b_2=2$}\label{b2=2}

Let $X$ be class VII surface with $b_2(X)=2$, $g$ a Gauduchon metric on $X$ and $\kg:=\frac{1}{2}\deg_g({\cal K})$.   The subspace ${\mathcal M}^\red(0,{\mathcal K})$ of reductions (split  polystable bundles) in the moduli space is a finite union of circles. For every $c\in\Tors$, $d\in e_1$ we denote
$$\Rg'_c:=\left\{{\mathcal L}\oplus({\mathcal K}\otimes{\mathcal L}^{\vee})\ \vline\ {\mathcal L}\in\Pic^c_{=\kg}\right\}\ ,\ \Rg''_d:=\left\{{\mathcal L}\oplus({\mathcal K}\otimes{\mathcal L}^{\vee})\ \vline\ {\mathcal L}\in\Pic^d_{=\kg}\right\}\ .$$
and we put $\Rg':=\union_{c\in\Tors}\Rg'_c$, $\Rg''=\union_{d\in e_1}\Rg'_d$. One has ${\mathcal M}^\red(0,{\mathcal K})=\Rg'\cup\Rg''$.
The filtrable bundles can be also easily classified. 
  In our case we have only four extension types: $\emptyset$, $\{1\}$, $\{2\}$, $I_0=\{1,2\}$. Using the Riemann Roch theorem and  Serre duality we obtain  the following formulae for the
dimension of the extension space $\mathrm{Ext}^1({\mathcal K}\otimes{\cal L}^\vee,{\cal L})=H^1({\mathcal K}^{\vee}\otimes {\mathcal L}^{\otimes 2})$:
\begin{equation} \label{dimensions}
h^1({\mathcal K}^{\vee}\otimes {\mathcal L}^{\otimes 2})=
\left\{
\begin{array}{cl}
2+h^0({\mathcal K}^{\vee}\otimes {\mathcal
L}^{\otimes 2})+h^0({\mathcal K}^{\otimes 2}\otimes{\mathcal L}^{\otimes-2})&\hbox{for type }\emptyset\ ,\\
 1+h^0({\mathcal K}^{\vee}\otimes {\mathcal
L}^{\otimes 2})+h^0({\mathcal K}^{\otimes 2}\otimes{\mathcal L}^{\otimes-2})&\hbox{for type }\{1\}\ ,\ \{2\}\ ,\\
h^0({\mathcal K}^{\vee}\otimes {\mathcal
L}^{\otimes 2})+h^0({\mathcal K}^{\otimes 2}\otimes{\mathcal L}^{\otimes-2})&\hbox{for type } I_0 \ .
\end{array}\right.
\end{equation}
Using the vanishing  Lemma \ref{nak}, one gets
\begin{lm}  \label{vanishing} Suppose that $X$ is   minimal. Then $h^0({\mathcal K}^{\otimes 2}\otimes{\mathcal L}^{\otimes-2})=0$   for type $\emptyset$, type $\{1\}$ and type $\{2\}$ extensions, whereas $h^0({\mathcal K}^{\vee}\otimes {\mathcal
L}^{\otimes 2})=0$   for type $I_0$ extensions. 
Moreover, when $X$ is not an Enoki surface, $h^0({\mathcal K}^{\otimes 2}\otimes{\mathcal L}^{\otimes-2})=0$   for type $I_0$ extensions, except  the case when ${\mathcal L}^{\otimes 2}={\mathcal K}^{\otimes 2}$,  i.e. when ${\mathcal L}$ has the form ${\mathcal R}\otimes{\mathcal K}$ for  a square root ${\mathcal R}$ of ${\mathcal O}$. 
\end{lm}
Therefore, under these assumptions, for any square root ${\mathcal R}$  we get an (essentially unique)   nontrivial extension of type $I_0$ 
$$0\map {\mathcal K}\otimes{\mathcal R}\map{\mathcal A}_{\mathcal R}\map{\mathcal R}\map 0\ .
$$
One has ${\mathcal A}_{\mathcal R}={\mathcal A}\otimes{\mathcal R}$, where ${\mathcal A}:={\mathcal A}_{{\mathcal O}}$ is the ``canonical extension" introduced in the introduction.  Recall from \cite{Te2} that, for every $c\in\Tors_2(H^2(X,\Z))$, the space $\Pic^c$ contains   two square roots of ${\cal O}$, which are conjugate by the  square root ${\cal R}_0\in\Pic^0(X)$ associated with the standard representation $\rho:\pi_1(X)\to \Z_2$. 

For the other right hand terms  in (\ref{dimensions}), wee see that the problem simplifies further  as soon as ${\mathcal L}$ defines  an extension  with semistable   middle term ${\mathcal E}$. More precisely:
\begin{lm}\label{ineq}   Let $g$ be a Gauduchon metric on $X$. There exists $\varepsilon>0$ such that  for any line bundle  ${\mathcal L}$  on $X$ with $\deg_g({\mathcal L})<  \frac{1}{2}\deg_g({\mathcal K})+\varepsilon$ one has $H^0({\mathcal K}^{\vee}\otimes {\mathcal
L}^{\otimes 2})=0$.
\end{lm}
\pf Indeed,  if $H^0({\mathcal K}^{\vee}\otimes {\mathcal
L}^{\otimes 2})\ne 0$, the vanishing locus of a non-trivial section would be an effective  divisor of volume $<2\varepsilon$. Therefore, if $\varepsilon$ is sufficiently small, this divisor will be empty, which would imply  ${\mathcal L}^{\otimes 2}\simeq {\mathcal K}$. But $\bar k$ is not divisible by 2 in $H^2(X,\Z)/\Tors$.  
\qed

{\it From now on we will always suppose that $X$ is minimal and is not an Enoki surface. }
\begin{lm} \label{degree} There exists a Gauduchon metric $g$ on $X$ such that $\deg_g({\cal K})<0$.
\end{lm}
The following short proof is due to Nicholas Buchdahl\footnote{Already in 2004  Nicholas Buchdahl pointed out to me that in \cite{Te2} one can assume that $\deg_g({\cal K})<0$ without loss of generality.}.   A different proof can be obtained  using   Lamari's description of the pseudo-effective cone of a non-K\"ahlerian surface \cite{Te3}, \cite{La}.

Suppose   that $g$ is a Gauduchon metric with $d:=\deg_g({\cal K})>0$ and let $\kappa$ be a real (1,1)-form representing the Chern class of ${\cal K}$ in Bott-Chern cohomology. One has $\int_X\kappa\wedge\kappa=k^2=-b$ (where again $b:=b_2(X)$). We claim that the $dd^c$-closed form $\eta:=\omega_g+\frac{d}{b}\kappa$ satisfy all  the conditions in  Buchdahl's positivity criterion (see \cite{Bu2} p. 1533), namely:
\begin{enumerate}
\item $\int_X \eta\wedge \eta> 0$,  $\int_X \eta\wedge\omega>0$, where $\omega$ is a strictly positive $dd^c$-closed (1,1)-form.
\item $\int_{C} \eta>0$ for every irreducible curve $C\subset X$ with $C^2<0$.
\end{enumerate}
Indeed,    $\int_X \eta\wedge \eta=\int_X \eta\wedge\omega=\int_X\omega_g^2+\frac{d^2}{b}>0$,   
and for every curve $C\subset X$ one has $\int_{C} \eta=\int_{C}\omega+\frac{d}{b}  c_1({\cal K})\cdot [C] \geq  \int_{C}\omega$, because $c_1({\cal K})\cdot [C]\geq 0$ by Lemma \ref{nak}. Using Buchdahl's positivity criterion, it follows that there exists a smooth real function $\varphi$ such that $\eta+dd^c\varphi$ is strictly positive.  The (1,1)-form $\eta+dd^c\varphi +\varepsilon\kappa$  is  $dd^c$-closed and   strictly positive for sufficiently small $\varepsilon >0$; but $\int_X(\eta+dd^c\varphi +\varepsilon\kappa)\wedge\kappa=-\varepsilon b$, which shows that, with respect to the corresponding Gauduchon metric, ${\cal K}$ has negative degree.\qed

We are now  able to prove  the following important:
\begin{pr}\label{topman} If $\deg_g({\mathcal K})<0$ then the following holds:
\begin{enumerate}
\item The moduli space ${\mathcal M}^{\pst}(0,{\mathcal K})$ is a compact topological manifold which has a natural smooth holomorphic structure on  ${\mathcal M}^{\pst}(0,{\mathcal K})\setminus\Rg'$, extending the  standard holomorphic structure on ${\mathcal M}^{\st}(0,{\mathcal K})$.
\item If none of the classes $-e_1$, $-e_2$ is represented by a cycle, then the bundles ${\mathcal A}_{\mathcal R}$ (where ${\mathcal R}^{\otimes 2}={\mathcal O}$) are stable.
\end{enumerate}
\end{pr} 
\pf 1. By Proposition \ref{reg}, one has $H^2({\mathcal E}nd_0({\mathcal E}))=0$ for every $g$-polystable bundle ${\mathcal E}$ with $c_2({\mathcal E})=0$, $\det({\mathcal E})={\mathcal K}$. Indeed, since $\deg_g({\mathcal K})<0$,  a  bundle  ${\cal E}$ with these invariants and non-vanishing  $H^2({\mathcal E}nd_0({\mathcal E}))$ cannot be $g$-semistable. The claim follows now from Theorem \ref{cp},    Proposition \ref{extholstr} and  Proposition \ref{loc}.\\

2.  Using Corollary \ref{injectivity} we see that  ${\mathcal A}_{\mathcal R}$ cannot be written   as an extension of type $\emptyset$ (when $X$ is not an Enoki surface). By Proposition  \ref{prevresult} it cannot be  written as an extension of type $\{i\}$ (when $-e_i$  is not represented by a cycle). Therefore,    under our assumptions, it could only be destabilized by a line bundle of the form ${\mathcal K}\otimes{\mathcal R}'$, with ${\mathcal R'}^{\otimes 2}={\mathcal O}$. But $\deg_g({\mathcal K}\otimes{\mathcal R}')=\deg_g{\mathcal K}<\kg$ when $\deg_g({\mathcal K})<0$.
\qed
\\
Consider the  three vector bundles 
 $${\mathcal F}_I^{\leq\kg}:=\coprod_{{\mathcal L}\in\Pic^{e_I}_{\leq \kg}} H^1({\mathcal K}^{\vee}\otimes{\mathcal L}^{\otimes 2})\ ,\ I=\emptyset,\ \{1\},\ \{2\}.
 $$
${\mathcal F}_I^{\leq\kg}$ is a bundle  over $\Pic^{e_I}_{\leq\kg}$ (which is a finite union of punctured closed disks). According to (\ref{dimensions}), Lemma \ref{vanishing}, and Lemma \ref{ineq}, ${\mathcal F}_\emptyset^{\leq\kg}$ is a rank 2-bundle, whereas ${\mathcal F}_{\{1\}}^{\leq\kg}$, ${\mathcal F}_{\{2\}}^{\leq\kg}$ are line bundles.  We denote by  $\P({\mathcal F}_I^{\leq\kg})$ the projectivizations of these bundles.   
 $\P({\mathcal F}_\emptyset^{\leq\kg})$   is a $\P^1$-bundle over $\Pic^{T}_{\leq\kg}$, whereas 
 $\P({\mathcal F}_{\{1\}}^{\leq\kg})$, $\P({\mathcal F}_{\{2\}}^{\leq\kg})$ can be identified with 
 $\Pic^{e_1}_{\leq\kg}$, $\Pic^{e_2}_{\leq\kg}$ respectively, so they are finite unions of punctured closed disks.\\
 
 Let $\Pi_\emptyset$ be the space obtained by collapsing to points the fibers of  $\P({\mathcal F}_\emptyset^{\leq\kg})$ over the boundary $\Pic^{T}_{=\kg}$ of the base $\Pic^{T}_{\leq\kg}$. $\Pi_\emptyset$ is a topological manifold, because collapsing the fibers over $\Pic^{T}_{=\kg}$ is equivalent to gluing  a finite union of copies of $(S^1\times D^3)$  by identifying in the obvious way the boundary of this union with the boundary of   $\P({\mathcal F}_\emptyset^{\leq\kg})$. Let $\Cg'\subset\Pi_\emptyset$ be the subset formed by the points corresponding to  collapsed fibers (which is just a copy of the finite union of circles $\Pic^{T}_{=\kg}$), and denote by $\Pi_\emptyset^c$, $\Cg'_c$ the component of $\Pi_\emptyset$ (respectively $\Cg'$) corresponding to $c\in\Tors$.

 We assign a bundle  ${\mathcal E}(p)$ to each  point $p\in\Pi_\emptyset$ in the following way:  ${\mathcal E}(p)$ is the split polystable bundle ${\mathcal L}\oplus[{\mathcal K}\otimes{\mathcal L}^\vee]$ if $p$ is the collapsed fiber over $ {\mathcal L}$, and  ${\mathcal E}(p)$  is the bundle ${\mathcal E}({\mathcal L},\varepsilon)$ if $\deg({\mathcal L})<\kg$ and $p=({\cal L},[\varepsilon])$ with $\varepsilon\in H^1({\cal K}^\vee\otimes{\cal L}^{\otimes 2})$.
\begin{pr}\label{phi}  \begin{enumerate}
\item ${\mathcal E}(p)$ is polystable for every  $p\in \Pi_\emptyset$, and is stable when $p\in \Pi_\emptyset\setminus \Cg'$. The map $\varphi_\emptyset:\Pi_\emptyset\to {\mathcal M}^{\pst}(0,{\mathcal K})$ given by $\varphi_\emptyset(p)={\cal E}(p)$ has the properties:
\begin{enumerate}
\item Identifies homeomorphically $\Cg'_c$ with the circle of reductions  $\Rg'_c$ for  every class $c\in\Tors$.
\item Is injective. 
\end{enumerate}
\item   Assuming  $\deg_g({\mathcal K})<0$, it holds
\begin{enumerate}
\item $\varphi_\emptyset$  is a  holomorphic open embedding on $\Pi_\emptyset\setminus \Cg'$.
\item 
$\varphi_\emptyset$ maps homeomorphically $\Pi_\emptyset$ onto an open subspace of ${\mathcal M}^{\pst}(0,{\mathcal K})$.

\end{enumerate}
\end{enumerate}
\end{pr} 
\pf 1. It is clear that ${\mathcal E}(p)$ is a split polystable bundle when $p\in\Cg'$, and that the induced map  $\Cg'_c\to \Rg'_c$ is a homeomorphism for every $c\in\Tors$. When $p\in
\Pi_\emptyset\setminus\Cg'$, then ${\mathcal E}(p)$ is the non-trivial   extension ${\mathcal E}({\mathcal L},\varepsilon)$ with ${\mathcal L}\in  \Pic^{T}_{<\kg}$.  ${\mathcal E}(p)$ is not destabilized by ${\mathcal L}$, but a priori it could be destabilized by another line bundle. This would imply that ${\mathcal E}(p)$ can be written as an extension in a different way. Using Proposition \ref{extypes}, we see that this new extension is  of one of the four types $\emptyset$, $\{1\}$, $\{2\}$, $I_0$. By Corollary \ref{injectivity}, this would imply that $X$ has a curve in one of the classes $0$, $e_1$, $e_2$, $e_1+e_2$.   This is not possible, because $X$ is minimal (so the vanishing  Lemma \ref{nak} applies) and is not an Enoki surface.   A similar argument, based on Corollaries \ref{injectivity}, \ref{inject} proves injectivity. 
\\
  
2. (a) It is easy to construct a classifying holomorphic bundle on   $\P({\mathcal F}_\emptyset^{<\kg})\times X$ which induces $\varphi_\emptyset$. This proves that  $\varphi_\emptyset$ is holomorphic on $\Pi_\emptyset\setminus\Cg'$. But, when $\deg_g({\mathcal K})<0$, the space  ${\mathcal M}^{\st}(0,{\mathcal K})$ is a smooth complex surface, by Proposition \ref{reg}. Since $\varphi_\emptyset$ is injective, it is a holomorphic  open embedding on $\Pi_\emptyset\setminus\Cg'$.
\\

2. (b) This is a delicate point, because we   endowed ${\mathcal M}^{\pst}(0,{\mathcal K})$ with the topology induced by the Kobayashi-Hitchin correspondence from the corresponding moduli space of oriented ASD connections (see the comments at the beginning of section \ref{topprop}). By Proposition \ref{topman}, ${\mathcal M}^{\pst}(0,{\mathcal K})$ is a topological manifold. This substantial simplification of the problem is specific to the case $b_2=2$.   

We prove first continuity: continuity on $\Pi_\emptyset\setminus\Cg'$ is clear by 2(a). Continuity at the points of $\Cg'$ follows easily by elliptic semicontinuity taking into account the compactness of  ${\mathcal M}^{\pst}(0,{\mathcal K})$. But   $\Pi_\emptyset$ and ${\mathcal M}^{\st}(0,{\mathcal K})$ are both  topological manifolds, so $\varphi_\emptyset(\Pi_\emptyset)$ is  open in ${\mathcal M}^{\st}(0,{\mathcal K})$ and $\varphi_\emptyset:\Pi_\emptyset\to \varphi_\emptyset(\Pi_\emptyset)$ is a homeomorphism, by the  Brouwer invariance of domain theorem. 
 \qed 

 Let ${\mathcal M}^0$ be the  union of the connected components of ${\mathcal M}^{\pst}(0,{\mathcal K})$ which intersect $\varphi_\emptyset(\Pi_\emptyset)$.  For $c\in\Tors$ denote by ${\mathcal M}_c$ the connected component of ${\mathcal M}^{\pst}(0,{\mathcal K})$ (or, equivalently, of ${\mathcal M}^0$) which contains the circle $\Rg'_{c}$.  One can write
 \begin{equation}\label{unionofcomponents}
 {\mathcal M}^0=\union_{c\in\Tors} {\mathcal M}_c\ .
\end{equation}
Our next goals are
\begin{itemize}
\item  to prove that   ${\mathcal M}_c\ne {\mathcal M}_{c'}$ for $c\ne c'$, \item to describe geometrically the  components ${\mathcal M}_c$.
\end{itemize}
 
The main point is that any ${\mathcal M}_c$ is compact and contains a family of projective lines parameterized by a punctured disk. The idea is to prove that ${\mathcal M}_c\setminus\Rg'_c$ can be identified with an open subset of a (possibly blown up) ruled surface.
   
  We will need the following simple
 \begin{lm}\label{ruled} Let $g:(\P^1\setminus\{0\})\times\P^1\to   Z$ be a holomorphic open embedding into a connected compact complex surface $Z$.    
 There   exists a ruled  surface $r:S\to \P^1$ over $\P^1$ and a modification $\mu:Z\to S$ (which is either an isomorphism or obtained by  iterated blowing-ups at points lying  above $0\in \P^1$) such that $g$ factorizes as $g=j\circ g_0$, where $g_0:(\P^1\setminus\{0\})\times\P^1\stackrel{\simeq}{\to} r^{-1}(\P^1\setminus\{0\})$
 is a trivialization of $r$ over $\P^1\setminus\{0\}$, and $j:r^{-1}(\P^1\setminus\{0\})\to Z$ is the obvious open embedding inverting $\mu$ on $r^{-1}(\P^1\setminus\{0\})$. In particular $U:=\im(g)$ is  Zariski open.
 \end{lm} 
 \pf  Fix $p_0\in \P^1\setminus\{0\}$ an put $C_0:=g(\{p_0\}\times\P^1)$.   By Proposition 4.3 p. 192 in \cite{BHPV} we get a (locally trivial) ruled surface $r:S\to B$ over a compact curve $B$ and  a modification $\mu:Z\to S$,   such that $C_0$ does not meet any exceptional curve of $\mu$, and $\mu(C_0)$ is a fiber of $r$. It is easy to see that any other projective line $g(\{p\}\times\P^1)$ also has these two properties. Therefore we have an induced injective holomorphic map $u:\P^1\setminus\{0\}\simeq\C\to B$.  Such a map can only exist when $B$ is rational, and in this case it extends to an isomorphism $v:\P^1\to B$. The exceptional divisor of $\mu$ lies over the singleton $B\setminus\im(u)$.  
 \qed

Lemma \ref{ruled} has a generalization for non-connected, non-compact surfaces. For  a finite set  $A$ put $\P_A:=A\times\P^1$ (the disjoint union of $|A|$ copies of $\P^1$), and for $a\in A$ put $\P_a:=\{a\}\times\P^1$. We denote by $D$, $D_r$ the  open disk of radius 1, respectively $r$ (where $r\in(0,1)$),  by $D^\bullet$, $D^\bullet_r$ the corresponding punctured disks, and by $\Omega(r,1)$ the annulus of biradius  $(r,1)$. 
   \begin{co}\label{newruled} Let $Y$ be an arbitrary (not necessarily connected or compact) complex surface  and let $g:D^\bullet\times\P_A\to Y$  a  holomorphic open embedding  such that
\begin{enumerate}
\item $g((D\setminus D_r)\times\P_a)$ is closed in $Y$ for every $a\in A$,
 \item $Y\setminus g(\Omega(r,1)\times\P_A)$ is compact.
\end{enumerate} 

 Then 
\begin{enumerate}
\item Any component  of $Y_0\in\pi_0(Y)$ with  $Y_0\cap\im(g)=\emptyset$ is compact. 
\item $\pi_0(g):A\to \pi_0(Y)$ is injective.
\item Let $Y_a$ be the connected component  of $Y$ which contains  $g(D^\bullet\times\P_a)$. There exists a modification $\mu_a:Y_a\to D\times\P_a$ (which is either an isomorphism or obtained by  iterated blowing-ups at points lying  above $0\in D$) such that $\resto{g}{D^\bullet\times\P_a}=j_a\circ g_a $ where $g_a:D^\bullet\times\P_a\to D^\bullet\times\P_a$ is a biholomorphism over $D^\bullet$, and $j_a:D^\bullet\times\P_a\to Y_a$ is  the obvious open embedding inverting $\mu_a$ on $D^\bullet\times\P_a$.
\end{enumerate}
  \end{co}
 \pf Let $Y_0$ be a connected component of $Y$. If $Y_0\cap \im(g)=\emptyset$, then $Y_0$ is compact by condition  2. When $Y_0\cap \im(g)\ne\emptyset$, there exists a subset $A_0\subset A$ and a holomorphic  open embedding $g_0:D^\bullet\times\P_{A_0}\to  Y_0$ still satisfying  properties 1., 2. in the hypothesis.  Let $\Delta$ be the complement of $\bar D_{r}$ in $\P^1$. We glue $\Delta\times\P_{A_0}$ to $Y_0$ via the restriction $\resto{g_0}{\Omega(r,1)\times\P_{A_0}}$ and we get a {\it connected}, compact, complex surface $Z_0$ with an open embedding $f_0:(\P^1\setminus\{0\})\times\P_{A_0}\to Z_0$. Condition 1. assures that $Z_0$ is Hausdorff. By Lemma \ref{ruled} the image $f_0((\P^1\setminus\{0\})\times\P_a)$ is Zariski open in $Z_0$, for every $a\in A_0$.  Since a connected  surface cannot contain two disjoint non-empty Zariski open sets, we get $|A_0|=1$.  The third statement follows now directly from  Lemma \ref{ruled}. 
 \qed
 \begin{pr}\label{firststr} Suppose $\deg_g({\mathcal K})<0$. Then
 \begin{enumerate} 
 \item ${\mathcal M}_{c_1}\ne {\mathcal M}_{c_2}$ for $c_1\ne c_2$.
 \item  \label{compldiv} ${\mathcal M}_{c}\setminus \Rg'_c$ has a natural smooth holomorphic structure and ${\mathcal M}_{c}\setminus \varphi_\emptyset(\Pi_\emptyset^c)$ is a divisor ${\mathcal D}_c$. 
 \item  ${\mathcal D}_c$   is a smooth rational curve or a tree of smooth rational curves.
  \item \label{bij} For every $c\in \Tors$ there exists a unique $d(c)\in e_1$ such that $\Rg''_{d(c)}\subset {\mathcal M}_c$;  one has  $\Rg''_{d(c)}\subset{\mathcal D}_c$ and the assignment $\Tors\ni c\mapsto d(c)\in e_1$ is a bijection.
 \end{enumerate} 
 \end{pr}
 \pf  We know by Proposition \ref{topman}, that ${\mathcal M}^{\pst}(0,{\cal K})\setminus\Rg'$ is a  smooth open complex surface with  $|\Tors(H^2(X,\Z)|$ ends (towards the removed circles $\Rg'_c$). We apply Corollary \ref{newruled} to the open embedding
 $$\varphi_\emptyset:\P({\mathcal F}_\emptyset^{<\kg})\to {\mathcal M}^{\pst}(0,{\cal K})\setminus\Rg'\ .
 $$
 Corollary \ref{newruled} applies because the bundle ${\mathcal F}_\emptyset^{<\kg}$ is trivial; indeed, by a theorem of Grauert \cite{Gr}, on a Stein manifold  the  classification of holomorphic bundles coincides with the  classification of topological bundles.  This proves 1., 2., and 3.\\

For 4.  we need results from Donaldson theory (see sections \ref{donclass}, \ref{localstructure}). Every component ${\mathcal M}_c$ contains a single component $\Rg'_c$ of $\Rg'$. Suppose that it also contains $j$ components $\Rg''_{d_1},\dots \Rg''_{d_j}$ of  $\Rg''$, with $j\geq 0$. Removing standard  neighborhoods $\bar U'_c$, $\bar U''_{d_i}$ of the $j+1$ circles of reductions (or equivalently blowing up these reduction loci as in section \ref{localstructure}) we see that the sum  of the fundamental classes of the resulting (oriented) boundary components is homologically trivial in the moduli space ${\mathcal B}^*_a$ of irreducible oriented connections. Here we use  the boundary orientations induced by the orientation $\oo_{L_0}$ of the moduli space (see section \ref{localstructure}) which corresponds to a line bundle $L_0$ with torsion Chern class. Let $\mu(\gamma)$ be  the Donaldson $\mu$-class   associated with a generator $\gamma\in H_1(X,\Z)/\Tors$. Using Corollary \ref{localred}, and Lemmas \ref{orientation}, \ref{muclass} we obtain 
$$\langle \mu(\gamma), [\partial\bar U'_c]\rangle=\pm 1\ ,\ \langle \mu(\gamma), [\partial\bar U''_d]\rangle=\mp 1
$$
for any  $c\in\Tors$   and $d\in e_1$.   The sign difference comes from the relation between the orientations $\oo_{L_0}$, $\oo_{L_1}$ associated with line bundles $L_i$ having $c_1(L_0)\in\Tors$, respectively $c_1(L_1)\in e_1\cup e_2$ (see Lemma \ref{orientation}). Therefore  $j=1$, so ${\mathcal M}_c=\varphi_\emptyset(\Pi_ \emptyset ^c)\cup{\mathcal D}_c$ contains a single component $\Rg''_{d(c)}$ of $\Rg''$. This component is contained in ${\mathcal D}_c$, because, by Proposition \ref{phi}, $\im(\varphi_\emptyset)$ cannot intersect $\Rg''$. Finally, we  have $d(c_1)\ne d(c_2)$ for $c_1\ne c_2$, because $\Rg''_{d(c_i)}$ belong to different connected components of ${\cal M}^\pst(0,{\cal K})$.
\qed
\begin{re}\label{invinv}    The involution $\otimes\rho$ leaves invariant every divisor ${\mathcal D}_c$ as well as the irreducible component ${\mathcal D}_c^0$ which contains the circle $\Rg''_{d(c)}$.  
\end{re}
Indeed, $\otimes\rho$ obviously leaves invariant $\varphi_\emptyset(\Pi^c_\emptyset)$ and ${\cal R}''_{d(c)}$.
\qed

A posteriori (after proving the existence of a cycle) we obtain $\pi_1(X)=\Z$ by Theorem \ref{previousresults} (\ref{cycleimplies}). This gives $\Tors=\{0\}$, so by (\ref{unionofcomponents}) the space ${\cal M}^0$ has a single component ${\cal M}_0=\varphi_\emptyset(\Pi_\emptyset^0)\cup{\cal D}_0$.  We will see that $\pi_1(X)=\Z$ also implies  ${\mathcal D}_0= {\mathcal D}_0^0$ (see Proposition \ref{aposteriori} below), hence the divisor  ${\mathcal D}_0$ is in fact irreducible and smooth. Therefore ${\mathcal M}_0$ can be identified with  the space obtained from $\bar D\times\P^1$ by collapsing the fibers over $S^1$ to points. This space is homeomorphic to $S^4$ as explained in the introduction.  This makes clear why the assumption $\pi_1(X)=\Z$ (made in the introduction in order to avoid technical complications)  simplifies considerably the proof.
\begin{co} ${\cal M}^0=\overline{{\cal M}^\st_\emptyset}$, where ${\cal M}^\st_\emptyset$ denotes the set of stable bundles which can be written as extensions  of type $\emptyset$.
\end{co}
This follows immediately from Proposition \ref{firststr} (\ref{compldiv}).
\qed

Our next purpose is to determine the position and the shape of the locus 
$${\mathcal M}^\pst_{\{1\}}\cup {\mathcal M}^\pst_{\{2\}}$$
 of polystable extensions of types $\{1\}$ and $\{2\}$ in the moduli space.    
 We know by formula \ref{dimensions} and Remarks \ref{vanishing}, \ref{ineq}  that, under our assumptions, for every line bundle  ${\mathcal L}$ with $c_1({\mathcal L})\in e_1\cup e_2$ there exists a unique nontrivial extension  ${\mathcal E}({\mathcal L})$ of ${\mathcal K}\otimes{\mathcal L}^\vee$ by ${\mathcal L}$. We define the map 
$$\varphi_{12}:\Pic^{e_1}\map \{\hbox{Isomorphism classes of bundles on }X\}
$$
 by
 $${\mathcal L}\mapsto\left\{\begin{array}{ccc} {\mathcal E}({\mathcal L})&{\rm when} &\deg_g({\mathcal L})<\kg\\
{\ } {\mathcal L}\oplus({\mathcal K}\otimes{\mathcal L}^{\vee})&{\rm when} &\deg_g({\mathcal L})=\kg\\
{\ } {\mathcal E}({\mathcal K}\otimes{\mathcal L}^\vee)&{\rm when} &\deg_g({\mathcal L})>\kg\ .
   \end{array}\right.
 $$
 \begin{pr}\label{secondstr} Suppose that  $\deg_g({\cal K})<0$, and  \begin{enumerate}
 \item The classes $-e_1$, $-e_2$ do not contain any cycle, 
 \item The classes $\pm(e_1-e_2)$ do not contain any effective divisor. 
 \end{enumerate}
 Then  
 \begin{enumerate}
 \item $\varphi_{12}$ takes values in ${\mathcal M}^\pst(0,{\mathcal K})\setminus\Rg'$,  is holomorphic, injective and identifies $\Pic^{e_1}$ with   ${\mathcal M}^\pst_{\{1\}}\cup {\mathcal M}^\pst_{\{2\}}$.
\item For every $d\in e_1$ the  map $\varphi_{12}$ identifies $\Pic^{d}_{=\kg}$ with $\Rg''_{d}$.
 \item    For every $c\in \Tors$ the  map $\varphi_{12}$ defines an open embedding 
 $\Pic^{d(c)}\to {\mathcal D}_{c}$.
 \item The closure   of  the  component $\varphi_{12}(\Pic^{d(c)})$ in the  moduli space is precisely the  irreducible component  ${\mathcal D}^0_c$ of ${\mathcal D}_c$.
 \item ${\mathcal D}^0_c$ is obtained from $\varphi_{12}(\Pic^{d(c)})$ by adding two points ${\mathcal B}^1_c$, ${\mathcal B}^2_c$ which are  fixed under the involution $\otimes\rho$.
 \end{enumerate}
 \end{pr}
 \pf
1. We have already  mentioned above that an extension of type $\{1\}$ or $\{2\}$ cannot be written as an extension of type $\emptyset$. Since  $X$ is not an Enoki surface, we see by Corollary \ref{injectivity},  that an extension of type $\{1\}$ can be written as an extension of type $\{2\}$ if and only if both of them  are split. The condition that $-e_i$ is not represented by a cycle   guarantees that an extension of type $\{i\}$ cannot be written as a non-trivial extension of type $I_0=\{1,2\}$ (see Proposition \ref{prevresult}), whereas the condition that   $(e_j-e_i)$ is not represented by an effective divisor guarantees that an extension of type $\{i\}$ cannot be written as an extension of the same type in a different way. Therefore, the stability condition for  ${\mathcal E}({\mathcal L})$ reduces to $\deg({\mathcal L})<\kg$.  This also proves injectivity.  The holomorphy follows from the properties of the holomorphic structure on ${\mathcal M}^\pst(0,{\mathcal K})\setminus\Rg'$ established in Proposition \ref{extholstr}. 

2. This is obvious.  

3. Since  $\varphi_{12}(\Pic^{d(c)})$   contains the circle $\Rg''_{d(c)}$, it is contained in the connected component ${\mathcal M}_c$ of the moduli space which contains this circle (see  Proposition \ref{firststr}). On the other hand $\varphi_{12}(\Pic^{d(c)})$ does not intersect the locus $\varphi_\emptyset(\Pi_\emptyset)$ of type $\emptyset$-extensions, so it is contained in the complement ${\mathcal D}_c={\mathcal M}_c\setminus\varphi_\emptyset(\Pi_\emptyset^c)$. Since $\varphi_{12}$ is holomorphic and injective, the statement follows. 

4.  ${\mathcal D}^0_c$ was defined as the irreducible component of ${\mathcal D}_c$    containing the circle $\Rg''_{d(c)}$. 

5. It suffices to note that both ${\mathcal D}^0_c$ and  $\varphi_{12}(\Pic^{d(c)})$ are invariant under  $\otimes\rho$. 
\qed 

We will need a similar description of ${\mathcal M}^\pst_{\{1\}}\cup {\mathcal M}^\pst_{\{2\}}$ in the more difficult case when one of the classes $\pm(e_1-e_2)$ does contain an effective divisor. Since $X$ is not an Enoki surface, only one of these classes, say $e_2-e_1$ is represented by an effective divisor. 
\begin{lm}\label{A} Suppose that   $e_2-e_1$ is represented by an effective divisor $A$.
\begin{enumerate}
 \item $A$ is a smooth rational curve. 
\item One has ${\mathcal E}({\mathcal L})\simeq {\mathcal E}({\mathcal K}\otimes{\mathcal L}^\vee(-A))$ for every ${\mathcal L}\in\Pic^{e_1}$.
\end{enumerate}
\end{lm}
\pf 1. If $X$ had two irreducible curves, then $X$ possesses a global spherical shell by the result of Dloussky-Oeljeklaus-Toma \cite{DOT} mentioned  in Theorem \ref{previousresults} in the introduction. In this case the possible configuration of curves is  known, and in all cases $e_2-e_1$ is either not represented by an effective divisor, or it is represented   by a smooth rational curve. If $X$ has only one irreducible curve, then this curve will be $A$ (because $e_2-e_1$ is not a divisible class in $H^2(X,\Z)/\Tors$). Therefore $A$ is irreducible, so it is either a smooth rational curve or a rational curve with a node (see \cite{Na1} Lemma 2.2, Theorem 10.2). But in the latter case $A$ has arithmetic genus 1, whereas the adjunction formula gives $p_a(A)=0$; so only the first case is possible.\\

2.  We first prove  that, when ${\mathcal L}\not\simeq {\mathcal K}\otimes{\mathcal L}^\vee(-A)$,  the canonical morphism 
$s:{\mathcal K}\otimes {\mathcal L}^\vee(-A)\to {\mathcal K}\otimes {\mathcal L}^\vee$
 admits a lift   $i:{\mathcal K}\otimes {\mathcal L}^\vee(-A)\to {\mathcal E}({\mathcal L})$, or equivalently, that the canonical section $\sigma\in H^0( {\mathcal O}(A))$ has a lift in $H^0({\mathcal E}({\mathcal L})\otimes {\mathcal K}^\vee\otimes {\mathcal L}(A))$.  Consider the following diagram with exact lines and an exact column (compare with the proof of Proposition \ref{prevresult}):
$$
\begin{array}{cccccc}
&&&&&H^0(({\mathcal K}^\vee\otimes {\mathcal L}^{\otimes 2})(A))\\
&&&&&\downarrow\\
&&&&&H^0(({\mathcal K}^\vee\otimes {\mathcal L}^{\otimes 2})(A)_A)\\
&&&&&\ \ \downarrow u \\
\map&H^0({\mathcal E}({\mathcal L})\otimes {\mathcal K}^\vee\otimes {\mathcal L})&\map& H^0({\mathcal O})&\textmap{\partial}& H^1({\mathcal K}^\vee\otimes {\mathcal L}^{\otimes 2})\\
&\downarrow&&\ \ \downarrow a&&\ \ \downarrow v\\
\map&H^0({\mathcal E}({\mathcal L})\otimes {\mathcal K}^\vee\otimes {\mathcal L}(A))&\map&H^0( {\mathcal O}(A))&\textmap {\partial_A}&H^1({\mathcal K}^\vee\otimes {\mathcal L}^{\otimes 2}(A))
\end{array}
$$
The section $\sigma=a(1)\in H^0( {\mathcal O}(A))$ has a lift in $H^0({\mathcal E}({\mathcal L})\otimes {\mathcal K}^\vee\otimes {\mathcal L}(A))$  if and only if $\partial_A\circ a (1)=v(\varepsilon)$ vanishes, where $\varepsilon\in H^1({\mathcal K}^\vee\otimes {\mathcal L}^{\otimes 2})\setminus\{0\}$ is the invariant of the extension defining ${\mathcal E}({\mathcal L})$. But when ${\mathcal L}\not\simeq {\mathcal K}\otimes{\mathcal L}^\vee(-A)$, one has $h^0(({\mathcal K}^\vee\otimes {\mathcal L}^{\otimes 2})(A))=0$, $h^0(({\mathcal K}^\vee\otimes {\mathcal L}^{\otimes 2})(A)_A)=h^1({\mathcal K}^\vee\otimes {\mathcal L}^{\otimes 2})=1$, so $u$ is an isomorphism and $v$ vanishes. This proves the existence of a lift $i$ of $s$. Using Proposition 4.8, 5. in \cite{Te3} it is easy to see that $i$ is a bundle embedding. The corresponding extension of ${\cal L}(A)$ by ${\mathcal K}\otimes {\mathcal L}^\vee(-A)$ with central term ${\cal E}({\cal L})$ cannot be trivial, because there exists no non-trivial morphism ${\cal L}(A)\to {\cal E}({\cal L})$.
\qed

Using Lemma \ref{A} one obtains:
\begin{pr}\label{secondstra} Suppose that $\deg_g({\cal K})<0$, and
\begin{enumerate} 
\item The classes $-e_1$, $-e_2$ do not contain any cycle,
\item  The class  $e_2-e_1$ is represented by an effective divisor $A$.
\end{enumerate}
Put $a=c_1({\cal O}(A))\in H^2(X,\Z)$, $\ag:=\deg_g({\cal O}(A))\in\R$. Then 
\begin{enumerate}
\item For every $d\in e_1$ the map  
$$\varphi_d:[\Pic^d]_{\leq \kg}^{\geq \kg-\ag}\to \{\hbox{Isomorphism classes of bundles on }X\}\ ,$$
defined by 
$${\mathcal L}\mapsto\left\{\begin{array}{ccc}
{\ } {\mathcal L}\oplus({\mathcal K}\otimes{\mathcal L}^{\vee})&{\rm when} &\deg_g({\mathcal L})=\kg\\
{\ }{\mathcal E}({\mathcal L})&{\rm when} &\kg-\ag<\deg_g({\mathcal L})<\kg\\
({\mathcal K}\otimes{\mathcal L}^\vee(-A))\oplus{\mathcal L}(A)&{\rm when} &\deg_g({\mathcal L})=\kg-\ag\ ,
   \end{array}\right.
$$
maps continuously  $[\Pic^d]_{\leq \kg}^{\geq \kg-\ag}$ into  ${\mathcal M}^\pst_{\{1\}}$.
\item The involution $d\mapsto k-a-d$ on $H^2(X,\Z)$  is identity on the class $e_1\subset H^2(X,\Z)$.
The involution  $\iota:\Pic^{e_1}\to \Pic^{e_1}$ given by ${\mathcal L}\mapsto {\mathcal K}\otimes{\mathcal L}^\vee(-A)$ leaves invariant every connected component of $\Pic^{e_1}$.
\item For $d\in e_1$ let $D^d$ be the disk 
$$D^d:=\qmod{[\Pic^d]_{\leq \kg}^{\geq \kg-\ag}}{\iota}
$$
and $C^d$ the punctured 2-sphere $C^d:=D^d\cup_{u_d}\Pic^{k-d}_{\leq\kg}$ where $u_d$ is the natural isomorphism between the boundaries acting by  $u_d({\mathcal L})={\mathcal K}\otimes{\mathcal L}^\vee$ for ${\mathcal L}\in\Pic^d_{=\kg}$ (see the picture below).  The map 
$$\psi_{12}:\coprod_{d\in e_1}C^d\to \{\hbox{Isomorphism classes of bundles on }X\}$$
given by
$${\mathcal L}\mapsto\left\{\begin{array}{ccc} {\mathcal E}({\mathcal L})&{\rm when} &{\mathcal L}\in [\Pic^{e_1}]_{<\kg}^{>\kg-\ag}\cup\Pic^{e_2}_{<\kg}\\
{\mathcal L}\oplus({\mathcal K}\otimes{\mathcal L}^{\vee})&{\rm when} &{\mathcal L}\in \Pic^{e_1}_{=\kg}
   \end{array}\right.
$$
maps injectively  $\coprod_{d\in e_1} C^d$ onto ${\mathcal M}^\pst_{\{1\}}\cup {\mathcal M}^\pst_{\{2\}}$.
 \begin{figure}[!h]
\vspace{-11mm}
%\hspace{2cm}
\centering
\scalebox{0.6}
{\includegraphics{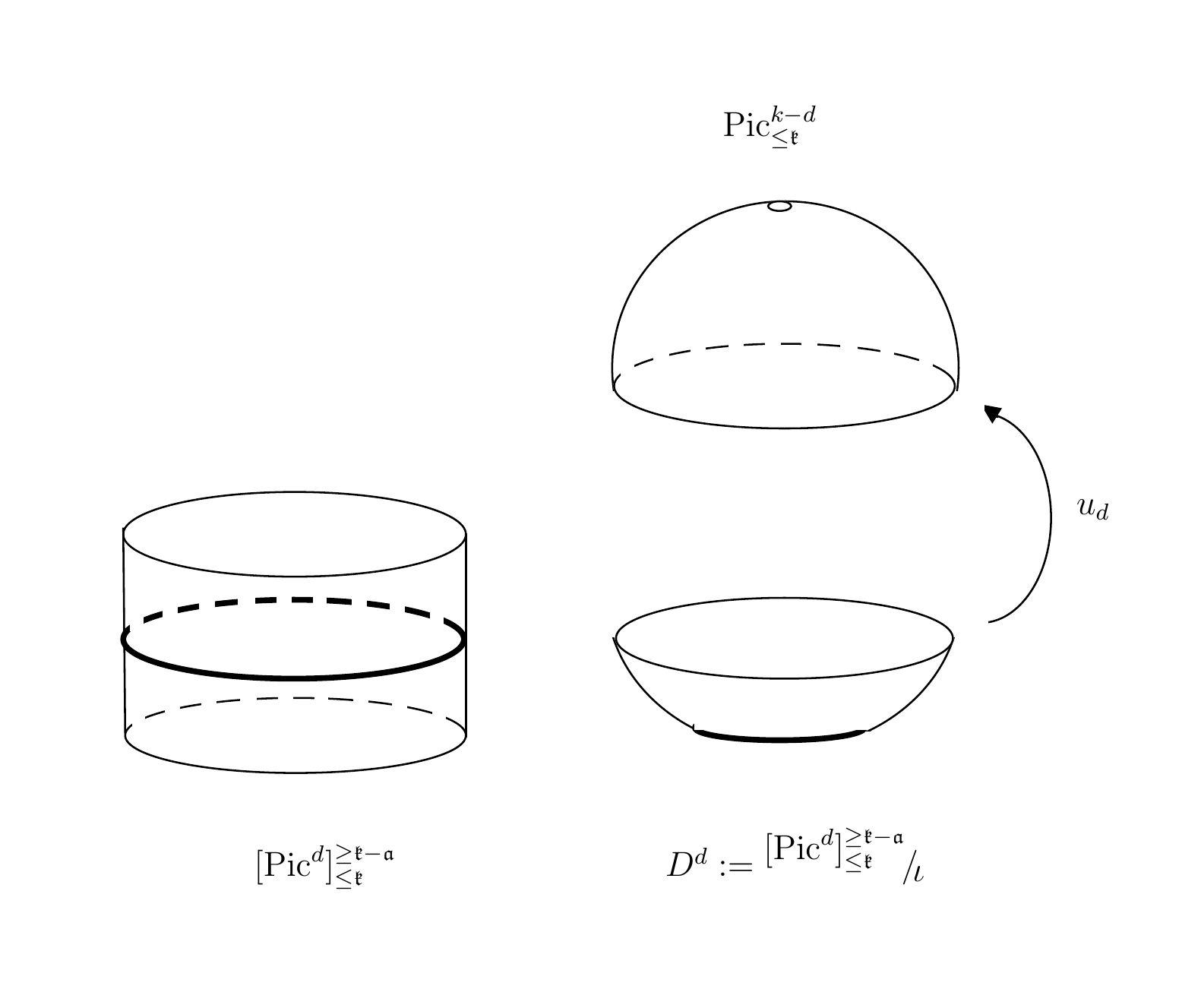}}
\label{divisor}
\vspace{-11mm}
\end{figure}

\item The map $\psi_{12}$ induces, for every $c\in \Tors$, an open holomorphic  embedding $C^{d(c)}\to {\mathcal D}^0_c$. The complement ${\mathcal D}^0_c\setminus \psi_{12}(C^{d(c)})$ is a point ${\mathcal B}_c$ which is fixed under the involution $\otimes\rho$.
\end{enumerate}
\end{pr}

\pf We denote by $a$ the Chern class of ${\mathcal O}(A)$ and by $\ag$ the degree of  ${\mathcal O}(A)$. \\

1. Let  ${\mathcal L}\in\Pic^{e_1}$. It suffices to note that a  subbundle of   ${\mathcal E}({\mathcal L})$ is isomorphic either  to ${\mathcal L}$ or   to ${\mathcal K}\otimes {\mathcal L}^\vee(-A)$. This follows from Corollary \ref{injectivity}  as in the proof of Proposition \ref{secondstr}. Continuity at points ${\cal L}\in [\Pic^d]_{\leq \kg}^{> \kg-\ag}$ is obvious, whereas continuity at points ${\cal L}\in \Pic^d_{=\kg-\ag} $ follows easily from Lemma \ref{A}, 2.\\

2. Statement (1) shows that the circles $\Rg''_d$ and $\Rg''_{k-a-d}$ belong to the same component of the moduli space. Therefore, by Proposition \ref{firststr}, \ref{bij}. one  has $d=k-a-d$. The second statement follows from the first.\\

3., 4. are proved using similar arguments as in the proof of Proposition \ref{secondstr}.
\qed

The following simple lemma concerns the geometry of the tree of rational curves ${\cal D}_c$. Recall that we denoted by ${\cal D}^0_c$ the irreducible component which contains the circle ${\cal R}''_{d(c)}$ of reductions.
\begin{lm}\label{newrem} Let ${\cal D}^1_c$ be an irreducible component of the tree ${\cal D}_c$ which is different from ${\cal D}^0_c$ (if such a component exists).  
\begin{enumerate}
\item When  the hypothesis of Proposition \ref{secondstr} holds, then ${\cal D}^1_c\cap {\cal D}^0_c$ is either empty or coincide with one of  the two singletons $\{{\cal B}^i_c\}$. 
\item When  the hypothesis of Proposition \ref{secondstra} holds, then ${\cal D}^1_c\cap {\cal D}^0_c$ is either empty or the singleton $\{{\cal B}_c\}$.
\end{enumerate}
\end{lm}
\pf Suppose that ${\cal D}^1_c\cap {\cal D}^0_c=\{{\cal F}_0\}$, where ${\cal F}_0\in {\cal D}^0_c$ does not coincide with any of the bundles  $\{{\cal B}^i_c\}$ (respectively ${\cal B}_c$).  According to  Propositions \ref{secondstr} , \ref{secondstra}, the bundle ${\cal F}_0$ is a  (possibly split) extension of type  $\{1\}$  or $\{2\}$. Therefore the set of  points of ${\cal D}^1_c$ which correspond to filtrable bundles is non-empty, but finite:   it consists of  ${\cal F}_0$ and possibly some of the (finitely many)    bundles ${\cal A}_{{\cal R}}$. Here we use essentially the fact that, for $c\ne c'$, the projective lines  ${\cal D}^0_c$, ${\cal D}^0_{c'}$ belong to different connected components of the moduli space, so ${\cal D}^1_c$ cannot be another component of the filtrable locus  ${\cal M}^{\pst}_{\{1\}}\cup {\cal M}^{\pst}_{\{2\}}$. 

If  ${\cal F}_0$ is a non-split extension (and hence stable), we get already a contradiction, because there cannot exist a family of rank 2 simple bundles on $X$ parameterized by a closed Riemann surface which contains both filtrable and non-filtrable bundles (see Corollary 5.3 in \cite{Te2}). If ${\cal F}_0$ is a split extension ${\cal L}\oplus ({\cal K}\otimes {\cal L}^{-1})$ with ${\cal L}\in \Pic^{e_1}$, $\deg_g({\cal L})=\kg$, we have to take into account the way in which the holomorphic structure of ${\cal M}^\pst(0,{\cal K})$ around $\Rg''$  has been defined (see section \ref{general}). The inclusion ${\cal D}^1_c\hookrightarrow  {\cal M}^\pst(0,{\cal K})$ is (by definition) holomorphic if,   replacing ${\cal F}_0$ by the nontrivial extension ${\cal E}({\cal L})$, we obtain a holomorphic curve  in ${\cal M}^\st_{g_t}(0,{\cal K})$ for small $t>0$ (see   Lemma \ref{holmap}). This  reduces the problem to the case when ${\cal F}_0$ is not split. \qed
\begin{pr}\label{aposteriori} Suppose that $\pi_1(X)=\Z$. Then the subspace ${\cal M}^0$ of ${\cal M}^\pst(0,{\cal K})$ contains a single component  ${\cal M}_0=\varphi_\emptyset(\Pi_\emptyset^0)\cup{\cal D}_0$. Moreover, one has ${\cal D}_0={\cal D}_0^0$, so ${\cal D}_0$ is an irreducible smooth rational curve, and ${\cal M}_0$ can be identified with the space obtained from   $\bar D^\bullet\times\P^1$ by   collapsing  to a point each fiber over  $\partial \bar D$. 
\end{pr}
\pf By the hypothesis and the coefficients formula we have $H^2(X,\Z)\simeq\Z$, so $\Tors=\{0\}$. Using (\ref{unionofcomponents}) we get the first statement. For the second, suppose that ${\cal D}_0$ had more than one connected component.  ${\cal D}_0$ is connected, so there exists an irreducible component ${\cal D}^1_0\ne{\cal D}_0^0$ which intersects ${\cal D}^0_0$. By Lemma \ref{newrem} the intersection point will be a fixed point of the involution $\otimes\rho$. Since the tree ${\cal D}_0$ is obtained by successive blow ups applied to a smooth ruled surface, it cannot contain triple crossings, so ${\cal D}_0^1$ will also be invariant under  this involution. Therefore $\otimes\rho$ would have at least three fixed points, which contradicts Proposition \ref{fixedpoints}.
\qed

We come back to the case of a general minimal class VII surface with $b_2=2$ which is not an Enoki surface.  
\begin{pr}\label{Y}
Suppose    that  none of the classes $-e_1$, $-e_2$, $-e_1-e_2$ is  represented by a cycle. Choose a Gauduchon metric $g$ on $X$ such that $\deg_g({\mathcal K})<0$. Then 
\begin{enumerate}
\item The bundles ${\mathcal A}_{\mathcal R}$ are stable, and the map $\Tors_2(\Pic(X))\to {\mathcal M}^\st(0,{\mathcal K})$ defined by ${\mathcal R}\mapsto{\mathcal A}_{\mathcal R}$ is injective.
\item The bundles  ${\cal A}_{\cal R}$ do not belong to ${\mathcal M}^0$.
\item The subspace ${\mathcal M}^\pst(0,{\mathcal K})\setminus{\mathcal M}^0$ is a smooth compact complex surface whose only filtrable points are ${\mathcal A}_{\mathcal R}$, where ${\mathcal R}^{\otimes 2}={\mathcal O}$.
\end{enumerate}
\end{pr}
\pf 1. The stability of ${\mathcal A}_{\mathcal R}$ was stated in Lemma \ref{topman}.  The injectivity of the map ${\mathcal R}\mapsto{\mathcal A}_{\mathcal R}$ is a direct consequence of   Proposition \ref{prevresult}. Indeed, if ${\mathcal A}_{\mathcal R}\simeq {\mathcal A}_{\mathcal R'}$ with ${\cal R}\not\simeq {\cal R}'$ then ${\cal A}\simeq   {\cal A}_{{\cal R}\otimes {\cal R}'}$,  so ${\cal A}$ can be written as an extension with kernel ${\cal K}\otimes [{\cal R}\otimes {\cal R}']\not\simeq {\cal K}$. Therefore, by Proposition \ref{prevresult} there would  exist a cycle in the class $-e_1-e_2$.

2.  By Corollary  \ref{injectivity} and Proposition \ref{prevresult} we see easily that, under our assumptions, ${\mathcal A}_{\mathcal R}$ cannot be isomorphic to a type $\emptyset$ extension.  Therefore, if ${\mathcal A}_{\mathcal R}$ belonged to ${\mathcal M}^0$, it will belong to one of the divisors ${\mathcal D}_c$, $c\in \Tors$. We claim that it  belongs precisely to the irreducible component  ${\mathcal D}_c^0$ which contains the circle $\Rg''_{d(c)}$. Indeed, if it belonged to another component, say  ${\mathcal D}_c^1$, this component would be contained in the stable locus (use Lemma \ref{newrem}) and would contain both filtrable and non-filtrable points (because the filtrable locus of type $\emptyset$, $\{1\}$ or $\{2\}$ is contained in $\varphi_\emptyset(\Pi_\emptyset)\cup (\cup_{c\in\Tors}{\mathcal D}^0_c)$). This would contradict Corollary  5.3 in \cite{Te2}. Therefore ${\mathcal A}_{\mathcal R}\in {\mathcal D}^0_c$, as claimed. 
On the other hand, since  ${\mathcal A}_{\mathcal R}$ is not isomorphic to a   type $\{1\}$ or a type $\{2\}$ extension, we see by Propositions  \ref{secondstr},  \ref{secondstra} that it must be  isomorphic to either one of the two points ${\mathcal B}^i_c$, or with the point ${\mathcal B}_c$. But these points are fixed under the involution $\otimes\rho$  whereas, by 1.,   ${\mathcal A}_{\mathcal R}$  is not fixed under this involution.\\

3. This follows directly from 2.
\qed

 We denote by ${\cal M}^\s(0,{\cal K})$  the moduli space of simple bundles ${\cal E}$ on $X$ with $c_2({\cal E})=0$ and $\det({\cal E})={\cal K}$.
\begin{thry} \label{mapf} Suppose that $X$ has no cycle. Then there exists  a connected, compact smooth complex surface $Y$ and an open embedding $f:Y\hookrightarrow {\cal M}^\s(0,{\cal K})$, $y\mapsto {\cal E}_y$ with  the properties
\begin{enumerate}
\item  $H^2({\cal E}nd_0({\cal E}_y))=0$ for any $y\in Y$. 
\item The set of filtrable bundles in $f(Y)$ contains the bundle ${\cal A}$ and is contained in the finite set $\{{\mathcal A}_{\mathcal R}|\ {\cal R}^{\otimes 2}={\cal O}\}$.
\end{enumerate}
\end{thry}
\pf The claim follows directly  form Proposition \ref{Y} and Lemma \ref{degree} .  
\qed
\vspace{1mm}\\
{\bf Remark:} (the case $\deg_g({\cal K})>0$) The existence of Gauduchon metrics $g$ with  $\deg_g({\cal K})<0$ simplifies considerably the proof, but is it is not absolutely necessary. In the case $\deg_g({\cal K})>0$, one can  still prove Theorem \ref{mapf} using ideas similar to those developed in the case $b_2=1$ \cite{Te2}.
\section{Universal families}

The goal of this section is to prove that the open embedding  $f:Y\to {\cal M}^\s(0,{\cal K})$ whose existence is given by Theorem \ref{mapf} (under the assumption that $X$ had no cycle) is induced by a universal family ${\cal F}$ on $Y\times X$, whose determinant is isomorphic to $p_X^*({\cal K})\otimes  p_Y^*({\cal N})$, where ${\cal N}$ is a line bundle on $Y$. In the next section we will apply the Grothendieck-Riemann-Roch theorem to this family, and we will see that its Chern classes must satisfy a set of very restrictive relations.

 Let $E$ be a rank $2$-bundle on a compact complex  surface $X$ and let ${\mathcal L}$ be a fixed holomorphic structure on the determinant line bundle $L$. We denote by ${\mathcal M}^{\s}(E,{\mathcal L})$ the moduli space of simple holomorphic structures on $E$ which induce ${\mathcal L}$ on $L$, modulo the complex gauge group ${\cal G}^\C:=\Gamma(X,\SL(E))$\footnote{We emphasize  here that we consider holomorphic structures which induce precisely ${\mathcal L}$ on  $L=\det(E)$, not only    a structure which is isomorphic to ${\mathcal L}$.}.

Let $Y$ be any complex manifold and $\iota:Y\to  {\mathcal M}^{\s}(E,{\mathcal L})$ be a holomorphic map. Put  $\iota(y)={\mathcal E}_y$. A {\it universal family for the pair} $(\iota, {\mathcal N})$ (where ${\mathcal N}$ is a holomorphic line bundle on $Y$) is a holomorphic rank 2-bundle ${\mathcal F}$ on $Y\times X$, together with an isomorphism $\det({\mathcal F})\textmap{\simeq} p_Y^*({\mathcal N})\otimes p_X^*({\mathcal L})$ such that $\resto{{\mathcal F}}{\{y\}\times X}\simeq {\mathcal E}_y$ for all $y\in Y$.
\begin{pr}\label{holclass}
 Let ${\mathcal T}$ be a holomorphic line bundle  on $X$ such that 
 \begin{enumerate}
 \item $\chi(E\otimes{\mathcal T})=-1$ 
 \item $h^0({\mathcal E}_y\otimes{\mathcal T})=h^2({\mathcal E}_y\otimes{\mathcal T})=0$ for every $y\in Y$.
\end{enumerate}
Then 
\begin{enumerate}
\item For two representatives ${\mathcal E}_y^1$, ${\mathcal E}_y^2$ of the isomorphism class $\iota(y)$,  the two lines $H^1({\mathcal E}_y^i\otimes{\mathcal T})^{\otimes 2}$,  $i=1$, $2$ and the two planes ${\mathcal E}^i_y(x)\otimes  H^1({\mathcal E}_y^i\otimes{\mathcal T})$ $i=1$, $2$ can be canonically identified,  for any $y\in Y$, $x\in X$. 
\item The assignments  
$$y\mapsto H^1({\mathcal E}_y\otimes{\mathcal T})^{\otimes 2}\ ,\ (y,x)\mapsto {\mathcal E}_y(x)\otimes  H^1({\mathcal E}_y\otimes{\mathcal T})
$$
descend to a holomorphic line bundle ${\mathcal N}^{\mathcal T}_\iota$ on $Y$ and a universal family ${\mathcal F}^{\mathcal T}_\iota$ for the pair $(\iota, {\mathcal N}^{\mathcal T}_\iota)$ respectively.
\end{enumerate}
\end{pr}

\pf The determinant line bundles of all the ${\cal E}_y$s coincide (not only are isomorphic!)  with ${\cal L}$.  Two different holomorphic $\SL(2)$-isomorphisms ${\cal E}_y^1\to {\cal E}_y^2$ differ by composition with $-\id_{{\cal E}_y^1}$, which operates trivially on the tensor products  $H^1({\mathcal E}_y^i\otimes{\mathcal T})^{\otimes 2}$, ${\mathcal E}^i_y(x)\otimes  H^1({\mathcal E}_y^i\otimes{\mathcal T})$. \\

For the second statement  suppose for simplicity that  ${\cal E}_y$ is {\it regular} (i.e. it holds $H^2({\mathcal E}nd_0({\mathcal E}_y))=0$) for all $y\in Y$. This  case is sufficient for our purposes.  We denote by ${\cal H}^\s_{\reg}(E,{\cal L})$ the space of regular simple holomorphic structures (integrable semiconnections) on $E$ which induce the fixed holomorphic structure ${\cal L}$ on $\det(E)$, and by ${\cal M}^\s_\reg(E,{\cal L})$ the corresponding open part of the moduli space ${\cal M}^\s(E,{\cal L})$. After suitable Sobolev completions ${\cal H}^\s_{\reg}(E,{\cal L})$ becomes a Banach complex manifold. On the bundle $p_X^*(E)$ over the product ${\cal H}^\s_\reg(E,{\cal L})\times X$ we introduce the tautological holomorphic structure ${\mathscr E}$ (which is trivial in the ${\cal H}^\s_\reg(E,{\cal L})$-directions).   This holomorphic structure is   ${\cal G}^\C$-invariant but, unfortunately, it does not descend  to a holomorphic bundle on  ${\cal M}^\s_\reg(E,{\cal L})\times X$, because  the center $\{\pm 1\}$ of   ${\cal G}^\C$ operates non-trivially on $p_X^*(E)$ but trivially on its base  ${\cal H}^\s_\reg(E,{\cal L})\times X$. However one can factorize  ${\mathscr E}$ by the based gauge group group ${\cal G}^\C_{x_0}$ and get a bundle ${\mathscr F}$ on $\tilde {\cal M}^\s_\reg(E,{\cal L})\times X$, where $\tilde {\cal M}^\s_\reg(E,{\cal L}):={\cal H}^\s_\reg(E,{\cal L})/{\cal G}^\C_{x_0}$ is a  holomorphic  principal $\mathrm{PSL}(2,\C)$-bundle over ${\cal M}^\s_\reg(E,{\cal L})$. Denote  by  ${\cal M}^\s_0(E,{\cal L})$, $\tilde {\cal M}^\s_0(E,{\cal L})$ the open subspaces of ${\cal M}^\s_\reg(E,{\cal L})$ and  $\tilde {\cal M}^\s_\reg(E,{\cal L})$ consisting of isomorphism classes $[{\cal E]}$ for which $h^0({\cal E}\otimes {\cal T})=h^2({\cal E}\otimes {\cal T})=0$, and denote by $ \pi$ and $p$ the projections of $\tilde {\cal M}^\s_0(E,{\cal L})\times X$ on the two factors. By the Grauert local triviality theorem and the assumption, the sheaf $R^1(\pi)_*({\mathscr F}\otimes   p^*({\cal T}))$  is a line bundle on $\tilde {\cal M}^\s_0(E,{\cal L})$, which we denote by ${\mathscr R}$. The bundles ${\mathscr R}^{\otimes 2}$, ${\mathscr F}\otimes {\mathscr R}$ descend to ${\cal M}^\s_0(E,{\cal L})$ and ${\cal M}^\s_0(E,{\cal L})\times X$ respectively,  because the $\SL(2,\C)$-actions   on these bundles induce  well-defined $\mathrm{PSL}(2,\C)$-actions which lift the free $\mathrm{PSL}(2,\C)$-actions on $\tilde {\cal M}^\s_0(E,{\cal L})$ and $\tilde {\cal M}^\s_0(E,{\cal L})\times X$. We denote by ${\cal N}^{\cal T}$, ${\cal F}^{\cal T}$ the obtained bundles. By the hypothesis (and our simplifying regularity assumption) the map $\iota:Y\to {\cal M}^\s(E,{\cal L})$ takes values in ${\cal M}^\s_0(E,{\cal L})$. It suffices to put
${\cal N}^{\cal T}_\iota:=\iota^*({\cal N}^{\cal T})\ ,\ {\cal F}^{\cal T}_\iota:=(f\times \id_X)^*({\cal F}^{\cal T})\ .
$
\qed
\begin{re} Let $X$ be a class VII surface, and choose $E$ such that $c_2(E)=0$, $c_1(E)=c_1({\mathcal K}_X)$. Then any holomorphic line bundle ${\mathcal T}$    on $X$ with $c_1({\mathcal T})^2=-1$  satisfies the condition $\chi(E\otimes{\mathcal T})=-1$. Therefore, Proposition  \ref{holclass} applies as soon as
$$h^0({\mathcal T}\otimes {\mathcal E}_y)=h^2({\mathcal T}\otimes {\mathcal E}_y)=0\hbox{ for all }y\in Y\ .$$
\end{re}

This remark applies to the embedding  $f:Y\hookrightarrow {\cal M}^\s(0,Y)$ obtained in Theorem \ref{mapf}. Indeed, since the only filtrable points on $Y$ have the form ${\cal A}_{\cal R}$ it is easy to see that  $h^0({\mathcal T}\otimes {\mathcal E}_y)=h^2({\mathcal T}\otimes {\mathcal E}_y)=0$  for all $y\in Y$. Therefore
\begin{co} \label{familyf} For any holomorphic line bundle ${\mathcal T}$on $X$ with $c_1({\mathcal T})^2=-1$ there exists a universal family ${\cal F}^{\cal T}_f$ for the pair $(f,{\cal N}^{\cal T}_f)$.
\end{co}

\section{Grothendieck-Riemann-Roch computations}
\label{indexsection}

Let ${\cal F}$ be  a universal family ${\cal F}$ for the obtained map $f:Y\hookrightarrow {\cal M}^\s$ (see Corollary \ref{familyf}). In this section we will see that,   applying the  Grothendieck-Riemann-Roch theorem   to  the bundles ${\cal F}$ and ${\cal E}nd_0({\cal F})$  and the proper morphism $Y\times X\to Y$, we will obtain important information about the Chern classes of  ${\cal F}$ and also about the Chern classes of $Y$ itself.

Let   $X$ be a class  VII  surface with $b_2=2$, $Y$ an arbitrary compact complex surface, and ${\mathcal L}$  a holomorphic line bundle on $Y$. Throughout this section we will  consider   a holomorphic rank 2 bundle ${\mathcal F}$ on $Y\times X$ with the following properties

\begin{enumerate}
\item $\det({\mathcal F})\simeq p_Y^*({\mathcal N})\otimes p_X^*({\mathcal K})$, where  ${\mathcal N}$ is a line bundle on $Y$.
\item $c_2(\resto{{\mathcal F}}{\{y\}\times X})=0$, for all $y\in Y$.
\end{enumerate}
The K\"unneth decompositions of the Chern classes $c_1({\mathcal F})$, $c_2({\mathcal F})$ in rational cohomology have the form:
$$c_1({\mathcal F})=\eta\otimes 1+1\otimes k\ ,\ c_2({\mathcal F})=c\otimes 1+ s\otimes t+\sum_i \nu_i\otimes e_i+ \sigma\otimes \theta\ ,
$$
where 
\begin{itemize}
\item[(a)] $\eta=c_1({\mathcal N})$, $k=c_1({\mathcal K}_X)$, 
\item [(b)] $c:=c_2(\resto{{\mathcal F}}{Y\times\{x\}})$ for any $x\in X$,
\item [(c)] $t$, $\theta$ are generators of $H^1(X,\Z)$ and  $H^3(X,\Z)/\Tors$ respectively such that 
$$\langle \theta\cup t, [X]\rangle=-\langle t\cup \theta, [X]\rangle=1\ ,$$ 
\item [(d)] $s\in H^3(Y,\Z)/\Tors$, $\sigma\in H^1(Y,\Z)$, and 
\item [(e)] $(e_1,e_2)$ is a basis of $H^2(X,\Z)/\Tors$ such that $e_i^2=-1$ and $k=e_1+e_2$. 
\item [(f)] $\nu_i\in H^2(Y,\Z)/\Tors$.
\end{itemize}
It is convenient to write formally ${\mathcal F}$ as $\Fg\otimes p_Y^*(\Mg)$, where $\Mg$ is a formal line bundle on $Y$ of Chern class $\frac{\eta}{2}$.
The Chern classes of the formal rank 2-bundle $\Fg$ will be 
$$c_1(\Fg)=1\otimes k\ ,\ c_2(\Fg)=c_2({\mathcal F})- \frac{1}{2} \eta\otimes  k-\frac{1}{4} \eta^2\otimes  1=$$
$$=\frac{1}{4}\left\{ U\otimes 1+ S\otimes t+ 2\sum_i M_i\otimes e_i+4\sigma\otimes\theta\right\}\ .
$$
where 
$$U:= 4c-\eta^2\ ,\ M_i=2\nu_i-\eta\ ,\ S=4s\ . $$
Let ${\mathcal T}$ be a holomorphic line bundle of Chern class $\tau$ on $X$. Our next purpose is to compute  the Chern character
$ch \left(p_{Y!}({\mathcal F}\otimes p_X^*({\mathcal T}) )\right)$.
\def\td{{\rm td}}
Using the multiplicative property of the Todd class and of the Chern character, the Grothendieck-Riemann-Roch theorem gives
$$ch (p_{Y!}({\mathcal F}\otimes p_X^*({\mathcal T}))=(p_Y)_*\left[ch({\mathcal F}\otimes p_X^*({\mathcal T})\cup p_X^*({\rm td}(X))\right]=$$
$$=(p_Y)_*\left[(ch({\mathcal F})\cup (p_X)^*(ch({\mathcal T})\cup  {\rm td}(X))\right]\ .$$

\begin{re} The non-K\"ahlerian version of the Grothendieck-Riemann-Roch theorem \cite{OTT} for proper analytic maps gives an identity in the Hodge algebra, not in rational cohomology. In this section we will use a cohomological version of this Grothendieck-Riemann-Roch theorem, which can be deduced  from the index theorem for families of coupled $Spin^c$-Dirac operators. Our fiber $X$ is  not K\"ahlerian, so the operator ${1}/{\sqrt{2}}(\bar\partial+\bar\partial^*)$ does not coincide with the corresponding coupled $Spin^c$-Dirac operator; however the difference is an operator of order 0, so the usual index formula applies. We will need this formula only for sheaves ${\cal E}$ for which  $R^i((p_Y)_*({\cal E}))$ are locally free.
\end{re} 

In our case, one has
$$\td(X)=1-\frac{k}{2}\ , \ ch({\mathcal T})\cup  {\rm td}(X)=1+\frac{1}{2}(2\tau- k)+\frac{1}{2}( \tau^2-\tau k)[X]\ .
$$

We get 
$$p_{Y*}\left[ch({\mathcal F})\cup  p_X^*(ch({\mathcal T})\cup  {\rm td}(X))\right]= p_{Y*}\left[ch(\Fg)\cup p_Y^*ch(\Mg)\cup p_X^*(ch({\mathcal T})\cup  {\rm td}(X))\right] 
$$
$$=ch(\Mg)\cup(p_Y)_*\left[ch(\Fg)\cup \left(1+\frac{1}{2}(2\tau- k)+\frac{1}{2}( \tau^2-\tau k)[X]\right)\right]
$$
Writing formally  
$$ch(\Sg):=(p_Y)_*\left[ch(\Fg)\cup \left(1+\frac{1}{2}(2\tau- k)+\frac{1}{2}( \tau^2-\tau k)[X]\right)\right]\  ,
$$
one will have $ch (p_{Y!}({\mathcal F}\otimes p_X^*({\mathcal T})))=ch(\Mg)ch(\Sg)$, where   $ch(\Sg)$ is given by 
\begin{equation}
\begin{array}{ccl}
ch_0(\Sg)&=&\ \ \tau^2\ ,\\ \\
 ch_1(\Sg)&=& -\frac{1}{2} \sum_i \tau_i M_i\ ,\  \tau_i:=\tau e_i\ ,\\ \\
ch_2(\Sg)&=&-\frac{1}{24}(1+3\tau^2) U-\frac{1}{48} \sum_iM_i^2+\frac{1}{6} s\cup\sigma \ .
\end{array}
\end{equation}
These formulae are obtained as follows. Putting $c_i:=c_i(\Fg)$, $i=1$, 2, we have  
$$ch(\Fg)=2+c_1 +\frac{1}{2}(c_1^2-2c_2)+\frac{1}{6}(c_1^3-3c_1c_2)+\frac{1}{24}(c_1^4+2c_2^2-4c_1^2c_2)\ .
$$
The known formulae for the Chern classes of $\Fg$ give: 
$$c_1^2=-2\otimes[X]\ ,\ c_1^3=0\ ,\ c_1c_2=\frac{1}{4}(U\otimes k-2\sum_i M_i\otimes[X])\ ,\  
$$
$$c_1^2c_2=-\frac{1}{2} U\otimes[X]\ ,\ c_2^2=\frac{1}{4}\left( -\sum_i M_i^2\otimes [X]+ 2(S\cup\sigma)\otimes[X]\right) \ .
$$
In the computation of $c_1c_2$ we used  that $t\cup H^2(X,\Z)=0$ in $H^3(X,\Z)/\Tors$; this follows from Poincar\'e duality. Therefore
$$ch(\Fg)=2(1\otimes 1)+1\otimes k+\frac{1}{2}\left[-2\otimes[X]-\frac{1}{2}[U\otimes 1+S\otimes t+
2\sum_iM_i\otimes e_i+4\sigma\otimes\theta ]\right]-$$
$$\frac{1}{8}\left[U\otimes k-2\sum_i M_i\otimes[X]\right]+\frac{1}{12}\left[\frac{1}{4}[-\sum_i M_i^2\otimes[X]+2(S\cup\sigma)\otimes[X]]+U\otimes[X]\right]  .
$$
Regrouping with respect to the  degree  of the $X$ component, this reads:
$$ch(\Fg)=\left\{2(1\otimes 1)-\frac{1}{4} U\otimes 1\right\}+\left\{-\frac{1}{4} S\otimes t\right\} +\left\{1\otimes k-\frac{1}{2}\sum_i M_i\otimes e_i-\frac{1}{8} U\otimes k\right\}
$$
$$+\{-\sigma\otimes\theta\}+\left\{-1+\frac{1}{4}\sum_i M_i-\frac{1}{48}\sum_i M_i^2+\frac{1}{24} (S\cup\sigma)+\frac{1}{12}U\right\}\otimes[X]\ .
$$
Taking the product with $\left(1+\frac{1}{2}(2\tau- k)+\frac{1}{2}( \tau^2-\tau k)[X]\right)$ and projecting on $Y$ we get as claimed $ch_0(\Sg)= \tau^2$,  
 $ch_1(\Sg)= -\frac{1}{2} \sum_i \tau_i M_i$ and 
$$ch_2(\Sg)=-\frac{1}{8}(\tau^2-\tau k)U-\frac{1}{16}(2\tau-k)kU-\frac{1}{48}\sum M_i^2+\frac{1}{24} S\cup\sigma+\frac{1}{12}U=$$
$$U(\frac{1}{12}-\frac{1}{8}\tau^2+\frac{1}{16}k^2)-\frac{1}{48}\sum M_i^2+\frac{1}{24} S\cup\sigma=U(-\frac{1}{24}-\frac{1}{8}\tau^2)-\frac{1}{48}\sum M_i^2+\frac{1}{6} s\cup\sigma.
$$
\begin{pr}\label{strange} For any universal family ${\cal F}$ for the map $f:Y\to {\cal M}^\s$ given by Theorem  \ref{mapf} one has the identities $U=-M_j^2$, $s\cup\sigma=0$.
\end{pr}
 \pf Chose ${\cal T}={\cal O}$. We know the $H^i({\cal E}_y)=0$ for all $y\in Y$, hence in this case $p_{Y!}({\mathcal F}\otimes p_X^*({\mathcal T}))=0$. Therefore (since $ch(\Mg)$ is invertible) we  have $ch(\Sg)=0$, which gives
\begin{equation}\label{oo}
-\frac{1}{24}  U-\frac{1}{48} \sum_iM_i^2+\frac{1}{6} s\cup\sigma=0
 \end{equation}
 
Choose now ${\cal T}$ such that $\tau=e_j$. In this case $h^0({\cal T}\otimes {\cal E}_y)=h^2({\cal T}\otimes {\cal E}_y)=0$ and $h^1({\cal T}\otimes {\cal E}_y)=1$. Therefore  $p_{Y!}({\mathcal F}\otimes p_X^*({\mathcal T}))$ can be written as $-{\cal H}$ for a line bundle ${\cal H}$ on $Y$. Writing $h=c_1({\cal H})$ we get
$$ch(\Sg)=-ch({\cal H})ch(\Mg)^{-1}=-\exp(h-\frac{\eta}{2})\ .
$$
This gives $ch_2(\Sg)=-\frac{1}{2}ch_1(\Sg)^2$ (which is equivalent to the vanishing of the second Chern class  of the line bundle ${\cal H}$), so we have
$$-\frac{1}{24}  U-\frac{1}{48} \sum_iM_i^2+\frac{1}{6} s\cup\sigma=-\frac{1}{8} M_j^2-\frac{1}{8} U \ .
$$
Combined with (\ref{oo}), this proves the claimed formulae. 
 \qed

For the endomorphism bundle ${\cal E}nd_0({\cal F})$ one has 
$$c_2({\mathcal E}nd_0({\mathcal F}))=4c_2({\mathcal F})-c_1({\mathcal F})^2=$$
$$4(c\otimes 1+ s\otimes t+\sum_i \nu_i\otimes e_i+ \sigma\otimes \theta)-(\eta^2\otimes 1+2\eta\otimes k-2 \otimes [X])
$$
$$=(4c-\eta^2)\otimes 1+4 s\otimes t+4\sigma\otimes\theta+ 2\sum_i(2\nu_i-\eta)e_i+2\otimes[X]
$$
$$=U\otimes 1+S\otimes t+2\sum_i  M_i\otimes e_i+4\sigma\otimes\theta+2\otimes[X]\ .
$$

Since ${\mathcal E}nd_0({\mathcal F})$ is isomorphic to its dual, its odd Chern classes vanish. The Chern character of ${\mathcal E}nd_0({\mathcal F})$ is $ch({\mathcal E}nd_0({\mathcal F}))=3-c_2({\mathcal E}nd_0({\mathcal F}))+\frac{1}{12}c_2({\mathcal E}nd_0({\mathcal F}))^2$. We get
$$(c_2({\mathcal E}nd_0({\mathcal F}))^2=-4\sum_iM_i^2\otimes [X]+4U\otimes [X]+8S\sigma\otimes[X]\ ,
$$
$$ch({\mathcal E}nd_0({\mathcal F}))=3-\big[U\otimes 1+S\otimes t+2\sum_i  M_i\otimes e_i+4\sigma\otimes\theta+2\otimes[X] \big]+$$
$$
+\frac{1}{12}\big[-4\sum_iM_i^2\otimes [X]+4U\otimes [X]+8S\sigma\otimes[X] \big]=
$$
$$=(3-U\otimes 1) -S\otimes t -2 \sum_i M_i\otimes e_i -4\sigma\otimes\theta+\frac{1}{3}(-\sum_i M_i^2+U +8 s\sigma-6)\otimes[X],
$$

$$[p_Y]_*(ch({\mathcal E}nd_0({\mathcal F}))\cup \td(X))=[p_Y]_*\big[ch({\mathcal E}nd_0({\mathcal F}))\cup (1-\frac{k}{2})\big]=
$$
$$=\frac{1}{3}(-\sum_i M_i^2+U +8 s\sigma-6)-\sum _i M_i\ . $$

Therefore, setting ${\mathcal U}:=(p_Y)_!({\mathcal E}nd_0({\mathcal F}))$, we have
$$ch_0({\mathcal U})=-2\ ,\ ch_1({\mathcal U})=-\sum_i M_i\ ,\ ch_2({\mathcal U})=\frac{1}{3}(U-\sum_i M_i^2+ 8s\sigma)\ .
$$
Since the image of $Y$ by the obtained {\it open} embedding $f:Y\hookrightarrow {\cal M}^\st(0,{\cal K})$ is contained in the regular part of the moduli space, we get
$$ch(T^{1,0}_Y)=ch\{R^1(p_{Y*})({\cal E}nd_0({\cal F}))\}=-ch\{p_{Y!}({\cal E}nd_0({\cal F}))\}=-ch({\cal U})\ .
$$
In particular
$$c_1^{\Q}(T^{1,0}_Y)=M_1+M_2=2(\nu_1+\nu_2-\eta)\in 2\im [H^2(Y,\Z)\to H^2(Y,\Q)]\ .
$$
In other words
\begin{pr}\label{even} The Chern class $c_1(Y)$ is even modulo torsion, so  the intersection form of $Y$ is even. In particular  $Y$ is minimal.
\end{pr}

\section{End of the proof}

Using the results proved in the previous sections, we will show now that the  assumption ``$X$ has no cycle" leads to a contradiction.

\begin{thry} Any minimal class VII surface with $b_2=2$ has a cycle representing one  of the classes 0, $-e_1$, $-e_2$, $-e_1-e_2$.
\end{thry}
\pf If $X$ is an Enoki surface, it possesses a homologically trivial cycle and is a GSS surface. Suppose now that  $X$ is not an Enoki surface. If $X$ had no cycle we have proved that  there exist a smooth compact complex surface $Y$ and an embedding $f:Y\hookrightarrow {\cal M}^\s$ whose image is contained in the regular locus  and contains a finite non-empty set of  filtrable   points. 
Moreover, we know that $f$ admits a universal family ${\cal F}\to Y\times X$ whose determinant line bundle has the form $p_Y^*({\cal N})\otimes p_X^*({\cal K})$ for a line bundle ${\cal N}$ on $Y$. By the results  in \cite{Te2}, we see that $Y$ cannot be a union of curves, because if it was, one could find a curve which passes through a filtrable point, normalize it if necessary, and get a family which contradicts  Corollary 5.3 in \cite{Te2}. Therefore $a(Y)=0$. Since the intersection form of $Y$ is even by Proposition \ref{even}, we are left with the following possibilities:
\begin{enumerate}

\item $Y$ is a class VII surface with $a(Y)=0$ and $b_2(Y)=0$.
\item $Y$ is a K3 surface with $a(Y)=0$, or  a torus with  $a(Y)=0$,
\end{enumerate}{\ }\\
\vspace{2mm}  
{\it Case 1:  $Y$ is a class VII surface with $a(Y)=0$ and $b_2(Y)=0$.} 

Suppose  that $Y$ is a class VII surface with  $b_2(Y)=0$. Using the formula $U=-M_j^2$ proved in  Proposition \ref{strange}, we see that the class $U=4c-\eta^2 
\in H^4(Y,\Z)$  vanishes. In other words, the bundles ${\cal F}^x$ on $Y$ have trivial characteristic number $\Delta({\cal F}^x)=4c_2({\cal F}^x)-c_1({\cal F}^x)^2$.  
Choose a Gauduchon metric on $Y$ in order to give sense to stability. By Theorem 1.3 in \cite{Te5}, the set
$$X^\st:=\{x\in X|\ {\cal F}^x\hbox{ is stable}\}
$$
is {\it Zariski open} in $X$. Note that, in our non-K\"ahlerian framework, this statement is not obvious. The point is that in our case the parameter space  of the holomorphic family  $({\cal F}^x)_{x\in X}$ is compact. Stability with respect to a Gauduchon metric is  always an open condition, but in general not Zariski open \cite{Te5}.
\\\\
{\it Case 1a:  The family $({\cal F}^x)_{x\in X}$ is generically stable ($X^\st\ne\emptyset$)}.\\

In this case we get a map $\emptyset\ne X^\st\to {\cal M}^\st(0,{\cal N})$. Since $\Delta({\cal F}^x)=0$, it is easy to see that the moduli space ${\cal M}^\st(0,{\cal N})$ is 0-dimensional. Indeed, the expected dimension of the moduli space vanishes, whereas the non-regular points are non-trivial line bundle extensions (use the same argument   as in  Proposition 3.7 \cite{Te2}). Using the Riemann-Roch theorem and the methods explained in section  \ref{b2=2}, we see that the set of line bundles ${\cal L}\in\Pic(Y)$ for which $h^1({\cal L}^{\otimes 2}\otimes {\cal N}^\vee)\ne 0$ is discrete, and one always has $h^1({\cal L}^{\otimes 2}\otimes {\cal N}^\vee)\leq 1$, so the space of line bundle extensions in ${\cal M}^\st(0,{\cal N})$ is also discrete.

Therefore the map $X^\st\to {\cal M}^\st(0,{\cal N})$ is constant; let ${\cal F}_0$ be this constant. We use   the same argument as in the proof of Corollary 4.2 in \cite{Te5}:  The sheaf ${\cal L}:=[p_X]_*(p_Y^*({\cal F}_0^\vee)\otimes {\cal F})$ on $X$ has rank 1, because it is a line bundle on $X^\st$. One obtains a tautological morphism $p_X^*({\cal L})\otimes p_Y^*({\cal F}_0)\to {\cal F}$, which is a bundle isomorphism on $X^\st\times Y$. Its restriction  to a fiber  $X\times\{y\}$ is a morphism  ${\cal L}\otimes {\cal O}^{\oplus 2}_X\simeq {\cal L}\otimes {\cal F}_0(y)\to {\cal F}_y$,
which is a bundle embedding on $X^\st$. This would imply that all our bundles   ${\cal F}_y$ are filtrable, which is not the case.\\ 
\\
{\it Case 1b: The family $({\cal F}^x)_{x\in X}$ is not generically stable ($X^\st=\emptyset$).}\\

In particular, in this case all the bundles ${\cal F}^x$ are filtrable. The idea is to show that the maximal destabilizing line subbundle  of ${\cal F}^x$ is  independent of $x$ for generic $x\in X$. Using the methods introduced in \cite{Te5}, consider the Brill-Noether locus of the family:
$$BN_X({\cal F}):=\{(x,{\cal U})\in X\times\Pic(Y)|\ H^0({\cal U}^\vee\otimes {\cal F}^x)\ne 0\}\ ,
$$
(which is   a closed analytic set of $X\times\Pic(Y)$) and its compact subsets
$$BN_X({\cal F})_{\geq d}:=\{(x,{\cal U})\in BN_X({\cal F})|\ \deg_g({\cal U})\geq d\}\ . $$
We denote by $p_\Pic$ the projection of $ X\times\Pic(Y)$ on $\Pic(Y)$. Since $X$ is compact, it is easy to see -- by the open mapping theorem -- that  $\deg_g\circ p_\Pic$ is locally constant on $BN_X({\cal F})$.  The reason is that $\deg_g$ is pluriharmonic on $\Pic(Y)$  and  the sets $BN_X({\cal F})_{\geq d}$ are compact (see Remark 2.13 \cite{Te5} for details).

Denote by ${\cal C}$ the set of  irreducible  components of  $BN_X({\cal F})$. For any $C\in {\cal C}$ define $d_C\in\R$ by   $\deg_g\circ p_\Pic (C)=\{d_C\}$,  and note that $C$ is closed in $BN_X({\cal F})_{\geq d_C}$, so it is compact. The projections on $X$ of all these components (which are analytic subsets of $X$) cover $X$, so there exists    $C\in{\cal C}$ with $p_X(C)=X$. Choose an irreducible component $C_0$ with this property such that  $d_{C_0}$ is maximal. Such a component exists because -- for any $d\in\R$ --  the set $\{C\in{\cal C}|\ d_C\geq d\}$ is finite.

The set
$$Z:=p_X\big[\union_{\{C\in {\cal C}|\ d_C> d_{C_0}\}} C\big]
$$
is a finite union of   analytic subsets of dimension $\leq 1$, and for any $x\in X\setminus Z$, one obviously has   
\begin{equation}\label{degmax}
\degmax_g({\cal F}^x)=d_{C_0}\ ,
\end{equation}
where we used the notation $\degmax_g({\cal F}^x):=\sup\{\deg_g({\cal L})|\ H^0({\cal L}^\vee\otimes {\cal F}^x)\ne 0\}$. Note that, for a filtrable rank 2-bundle ${\cal F}$, the invariant  $\degmax_g({\cal F})$ is well-defined and the supremum in the definition of this invariant is  attained at some line bundle (see Lemma 4 in \cite{Bu1}, and Definition 2.6, Remark 2.7 in \cite{Te5}). Formula (\ref{degmax}) shows that, for every  $(x,{\cal L})\in C_0$ and $x\in X\setminus Z$,  the line bundle ${\cal L}$ is a ``maximal destabilizing line bundle" of ${\cal F}^x$.

In our case $\Pic(Y)$ can be identified with a finite union of copies of $\C^*$, and with respect to suitable identifications $\Pic^c(Y)\simeq\C^*$,  the restriction of   $\deg_g$ to   a component $\Pic^c(Y)$ has the form $\zeta\to  \ln|\zeta|$.   Since $\deg_g\circ p_\Pic$ is locally constant on $BN_X({\cal F})$, it follows that $p_\Pic$ is locally constant on $BN_X({\cal F})$, too. Let ${\cal L}_0\in\Pic(Y)$ be the line bundle which corresponds to $C_0$. We have obviously $h^0({\cal L}_0^\vee\otimes {\cal F}^x)>0$ for all $x\in X$. Using the fact that the  bundles ${\cal F}^x$ are non-stable, it is easy to see that  $h^0({\cal L}_0^\vee\otimes {\cal F}^x)\leq 2$ for any $x\in X\setminus Z$, and equality occurs if and only if ${\cal F}^x\simeq {\cal L}_0\oplus {\cal L}_0$ (see Lemma 4.1 in \cite{Te5}). Let $U\subset (X\setminus Z)$ the open Zariski subset where the map $x\mapsto h^0({\cal L}_0^\vee\otimes {\cal F}^x)$ takes its minimal value.
The sheaf
$${\cal T}:=[p_X]_*(p_Y^*({\cal L}_0^\vee)\otimes {\cal F})$$
has rank 1 or 2 and is locally free on $U$.  The obvious morphism
$p_X^*({\cal T})\otimes p_Y^*({\cal L}_0)\to {\cal F}$
is a bundle embedding on $U\times Y$. Restricting this morphism to fibers $\{y\}\times X$, we 
get a morphism ${\cal T}\to {\cal F}_y$ which is a bundle embedding on $U$. When ${\cal T}$ has rank 2, it follows that the bundles  ${\cal F}_y$ have all the same type (filtrable or non-filtrable) as the reflexivization of ${\cal T}$. When ${\cal T}$ has rank 1,  the bundles ${\cal F}_y$ contain all a rank 1 subsheaf, hence they are all filtrable.   But our manifold $Y$ contains both filtrable and non-filtrable points.
  \qed
{\ }\vspace{1mm}\\
{\it Case 2. $Y$ is a K3 surface with $a(Y)=0$, or  a torus with  $a(Y)=0$.}\\

This case is ruled out by Corollary  1.6 in \cite{Te5}, which states:\\

{\it  Let $(Y,g)$ be a compact K\"ahler manifold, $X$ a surface with $b_1(X)$ odd and $a(X)=0$. Let ${\cal E}\to X\times Y$ be an arbitrary family of rank 2-bundles on $Y$  parameterized by $X$. Then there exist a locally free  sheaf ${\cal T}_0$  on $X$ of rank 1 or 2, a non-empty open Zariski open set $U\subset X$ and, for every $y\in Y$, a morphism $e_y:{\cal T}_0\to {\cal E}_y$ which is a bundle embedding on $U$.
}\\

Therefore, such a family cannot contain both filtrable and non-filtrable bundles. For completeness, we explain below briefly the outline of the proof of this statement in the special case when $Y$ is a K3 surface or a bidimensional torus. This assumption implies that any  moduli space of simple oriented bundles on $Y$ is smooth. Indeed, since ${\cal K}_Y$ is trivial, we have $h^2({\cal E}nd_0({\cal {\cal S}}))=h^0({\cal E}nd_0({\cal {\cal S}}))=0$ for any simple bundle ${\cal S}$ on $Y$, by Serre duality.

Suppose that the family  ${\cal E}$ is generically stable, i.e. $X^\st\ne \emptyset$. ${\cal E}$ induces a holomorphic map $X^\st\textmap{f} {\cal M}^\st (E)$, where $E$ denotes a differentiable rank 2-bundle over $Y$ and ${\cal M}^\st (E)$ stands for the moduli space of stable holomorphic structures on $E$. In our special case we know that $\det({\cal E})\simeq p_X^*({\cal K})\otimes p_Y^*({\cal N})$ by Corollary \ref{familyf}, so $f$ factorizes through a morphism $X^\st\to {\cal M}^\st_{\cal N}(E)$, but we do not need this property.

The pull-back of the Petersson-Weil K\"ahler form of ${\cal M}^\st(E)$ is a  form $\eta\in A^{1,1}_\R(X^\st)$, which is closed and positive (in the non-strict sense). By Theorem 1.4   \cite{Te5}, this form extends as a closed positive current $pw({\cal E})$ on $X$. Let $R$ be the residual part of this current with respect to the Siu decomposition. When $f$ has rank 2  at a point $x_0\in X^\st$, $R$ is smooth and strictly positive  at $x_0$, so $[R]^2>0$ by the self-intersection inequality given by Theorem 5.3 \cite{Te5}. This contradicts the assumption $b_1(X)$ odd (see the signature theorem, Theorem 2.14 \cite{BHPV}). 

When $f$ has generically rank 1, note first that the level sets $C_{{\cal S}}\subset X^\st$, ${\cal S}\in\im(f)$  of  $f$ are all 1-dimensional, by the semicontinuity theorem for the fiber dimension.  The Brill-Noether locus $B_{\cal S}:=\{x\in X|\ h^0({\cal S}^\vee\otimes {\cal E}^x)>0\}$ is a closed analytic set of $X$. One has  $C_S=X^\st\cap B_{\cal S}$  (recall that for two non-isomorphic stable bundles ${\cal S}_i$ of the same rank and degree, one has $h^0({\cal S}_1^\vee\otimes{\cal S}_2)=0$). For any ${\cal S}\in\im(f)$ choose an irreducible component
$D_{\cal S}$ of $B_{\cal S}$ which intersects $X^\st$. Two such curves  are distinct, because they are disjoint on $X^\st$.  We conclude that $X$ contains infinitely many curves, which contradicts $a(X)=0$. When $f$ is constant, we proceed as in {\it Case 1a} above. 
The case $X^\st= \emptyset$  is treated as {\it Case 1b}, using the fact that $\Pic(Y)$ is K\"ahlerian.  
\qed

\vspace{3mm} 
{\small
Author's address: \vspace{2mm}\\
Andrei Teleman, LATP, CMI,   Universit\'e de Provence,  39  Rue F.
Joliot-Curie, 13453 Marseille Cedex 13, France,  e-mail:
teleman@cmi.univ-mrs.fr. }

\end{document}